\documentclass[UTF-8,reqno]{amsart}
\usepackage{enumerate}
\usepackage{amssymb,url,color, booktabs}
\usepackage{mathrsfs}
\usepackage{amsmath}
\usepackage{color}

\usepackage[colorlinks=true]{hyperref}
\hypersetup{
    linkcolor=blue,          
    citecolor=red,        
    filecolor=blue,      
    urlcolor=cyan
}

\usepackage[nobysame]{amsrefs}
\BibSpec{article}{%
+{}{\PrintAuthors} {author}
+{,}{ \textrm} {title}
+{.}{ \textit} {journal}
+{,}{ \textbf} {volume}
+{}{ \parenthesize} {date}
+{,}{ } {pages}
+{.}{ arXiv:} {eprint}
+{.}{} {transition}
}
\BibSpec{book}{%
+{}{\PrintAuthors} {author}
+{,}{ \textit} {title}
+{.}{ \textrm} {series} 
+{,}{ Vol.} {volume}
+{.}{ } {publisher}
+{,}{ } {date}
+{.}{} {transition}
}
\usepackage{color}
\usepackage[colorlinks=true]{hyperref}
\hypersetup{
    linkcolor=blue,          
    citecolor=red,        
    filecolor=blue,      
    urlcolor=cyan
}

\numberwithin{equation}{section}

\newcommand{\be}{\begin{eqnarray}}
\newcommand{\ee}{\end{eqnarray}}
\newcommand{\ce}{\begin{eqnarray*}}
\newcommand{\de}{\end{eqnarray*}}
\newtheorem{theorem}{Theorem}[section]
\newtheorem{lemma}[theorem]{Lemma}
\newtheorem{remark}[theorem]{Remark}
\newtheorem{definition}[theorem]{Definition}
\newtheorem{proposition}[theorem]{Proposition}
\newtheorem{Examples}[theorem]{Example}
\newtheorem{corollary}[theorem]{Corollary}

\def\nor{|\!|\!|}
\def\wt{\widetilde}

\def\eps{\varepsilon}

\def\e{\mathrm{e}}

\def\p{\partial}

\def\[{{\Big[}}
\def\]{{\Big]}}
\def\<{{\langle}}
\def\>{{\rangle}}
\def\({{\Big(}}
\def\){{\Big)}}

\def\bx{{\mathbf{x}}}
\def\tr{\mathrm {tr}}

\def\dif{{\mathord{{\rm d}}}}

\def\min{{\mathord{{\rm min}}}}

\def\no{\nonumber}
\def\={&\!\!=\!\!&}

\def\bB{{\mathbf B}}
\def\bC{{\mathbf C}}

\def\cA{{\mathcal A}}
\def\cB{{\mathcal B}}

\def\cF{{\mathcal F}}

\def\cI{{\mathcal I}}

\def\cL{{\mathcal L}}

\def\cP{{\mathcal P}}

\def\cR{{\mathcal R}}
\def\cS{{\mathcal S}}
\def\cT{{\mathcal T}}

\def\mB{{\mathbb B}}

\def\mE{{\mathbb E}}

\def\mG{{\mathbb G}}

\def\mI{{\mathbb I}}

\def\mL{{\mathbb L}}
\def\mM{{\mathbb M}}
\def\mN{{\mathbb N}}

\def\mP{{\mathbb P}}

\def\mR{{\mathbb R}}
\def\mS{{\mathbb S}}

\def\bB{{\mathbf B}}

\def\1{{\mathbf{1}}}

\def\sF{{\mathscr F}}

\def\sI{{\mathscr I}}

\def\sL{{\mathscr L}}
\def\sM{{\mathscr M}}

\def\geq{\geqslant}
\def\leq{\leqslant}

\def\div{\mathord{{\rm div}}}

\def\eps{\varepsilon}

\def\e{\mathrm{e}}

\def\p{\partial}

\def\[{{\Big[}}
\def\]{{\Big]}}
\def\<{{\langle}}
\def\>{{\rangle}}
\def\({{\Big(}}
\def\){{\Big)}}

\def\bx{{\mathbf{x}}}
\def\tr{\mathrm {tr}}

\def\dif{{\mathord{{\rm d}}}}

\def\min{{\mathord{{\rm min}}}}

\def\no{\nonumber}
\def\={&\!\!=\!\!&}
\def\bt{\begin{theorem}}
\def\et{\end{theorem}}
\def\bl{\begin{lemma}}
\def\el{\end{lemma}}
\def\br{\begin{remark}}
\def\er{\end{remark}}
\def\bx{\begin{Examples}}
\def\ex{\end{Examples}}
\def\bd{\begin{definition}}
\def\ed{\end{definition}}
\def\bp{\begin{proposition}}
\def\ep{\end{proposition}}
\def\bc{\begin{corollary}}
\def\ec{\end{corollary}}

\def\geq{\geqslant}
\def\leq{\leqslant}

\def\div{\mathord{{\rm div}}}

\def\<{\langle} \def\>{\rangle}

\allowdisplaybreaks

\begin{document}

\title{Cauchy Problem of Stochastic Kinetic Equations}
\date{}

\author{Xiaolong Zhang,\ \ Xicheng Zhang}

\address{Xiaolong Zhang:
School of Mathematics and Statistics, Wuhan University,
Wuhan, Hubei 430072, P.R.China\\
Email: zhangxl.math@whu.edu.cn
 }

\address{Xicheng Zhang:
School of Mathematics and Statistics, Wuhan University,
Wuhan, Hubei 430072, P.R.China\\
Email: XichengZhang@gmail.com
 }

\thanks{
This work is supported by NNSFC grant of China (No. 11731009) and the DFG through the CRC 1283
``Taming uncertainty and profiting from randomness and low regularity in analysis, stochastics and their applications''.
}

\begin{abstract}
In this paper we establish the optimal regularity estimates for the Cauchy problem of stochastic kinetic equations 
with random coefficients in anisotropic Besov spaces.
As applications, we study the nonlinear filtering problem for a degenerate diffusion process, and obtain the existence and
regularity of conditional probability densities under few assumptions. 
Moreover, we also show the well-posedness for a class of super-linear growth 
stochastic kinetic equations driven by velocity-time white noises, as well as 
a kinetic version of Parabolic Anderson Model with measure as initial values.


\end{abstract}

\maketitle \rm

\tableofcontents

\section{Introduction}
Let $\{(W^k_t)_{t\geq 0},\ k\in\mathbb{N}\}$ be a sequence of independent one-dimensional standard Brownian motions 
defined on some stochastic basis $(\Omega,\sF,\mathbb{P};(\sF_t)_{t\geq0})$.
Let $\cP$ be the $(\cF_t)$-predictable $\sigma$-algebra over $\mR_+\times\Omega$, and
$\mM^d_{sym}$ the set of all $d\times d$-symmetric positive-definite matrices.
Given $\cP\times\cB(\mR^{2d})$-measurable processes 
$$
(a,b,\sigma):\mathbb{R}_+\times\Omega\times\mathbb{R}^{2d}\to(\mM^d_{sym}, \mathbb{R}^d,\mathbb{R}^d\otimes\ell^2),
$$
we introduce the following operators: for a function $u$ of $(x,v)\in\mR^{2d}$,
$$
\sL_v u:=\tr(a\cdot\nabla^2_v u)+b\cdot\nabla_v u,\ \ \sM^k_v u:=\sigma^k\cdot\nabla_v u,\ \ k\in\mN,
$$
where $\nabla_v$ (resp. $\nabla_x$) stands for the gradient operator that acts only on the variable $v$ (resp. $x$).
In this paper we consider the following stochastic kinetic equation (SKE for short) of It\^o's type:
\begin{align}\label{spde2}
\dif u=\big[\sL_v u+v\cdot\nabla_x u+f\big]\dif t+\big[\sM^k_vu+g^k\big]\dif W^k_t,
\end{align}
as well as its adjoint form
\begin{align}\label{spde22}
\dif u=\big[\sL^*_vu-v\cdot\nabla_x u+f\big]\dif t+\big[(\sM^k_v)^*u+g^k\big]\dif W^k_t,
\end{align}
subject to the initial condition $u(0,\omega,x,v)=u_0(\omega,x,v)$ being $\sF_0$-measurable,
where the asterisk stands for the adjoint operator, the unknown $u$ is a function of $(t,\omega,x,v)$, and $(x,v)$ represents the position and velocity,
the nonhomogenous or free terms
$$
(f, g):\mathbb{R}_+ \times\Omega\times\mathbb{R}^{2d}\to (\mR, \ell^2)
$$
are $\cP\times\cB(\mR^{2d})$-measurable processes.
Here and below we use the usual Einstein convention for summation: An index appearing in a product
will be summed automatically. Notice that SPDEs \eqref{spde2} and \eqref{spde22} are highly degenerate in $x$-direction.

\smallskip

When $\sigma$, $g$ are zero and $a,b,f$ do not depend on $\omega$, PDEs \eqref{spde2} and \eqref{spde22} are  degenerate deterministic 
equations and called kinetic Fokker-Planck-Kolmogorov equations in the literature 
since they are just the associated forward and backward Kolmogorov equations of the following SDE
\begin{align}\label{SDE0}
\dif X_t=V_t\dif t,\ \ \dif V_t=b(t,X_t,V_t)\dif t+\sqrt{2a}(t,X_t,V_t)\dif B_t,
\end{align}
where $B$ is a $d$-dimensional standard Brownian motion.
Some backgrounds of nonrandom PDEs \eqref{spde2} and \eqref{spde22} are referred to \cite{Li, So, Vi} and references therein. 
We mention that the two-sides estimates of the distribution density of SDE \eqref{SDE0} (also called heat kernel of $\sL_v+v\cdot\nabla_x$)
were studied by Delarue and Menozzi \cite{De-Me}.
The Schauder estimates for deterministic kinetic equation \eqref{spde2} were established in \cite{Lo, Lu, Pr},
and the maximal $L^p$-regularity estimates in $(x,v)$ were obtained in \cite{Bo} (see also \cite{chen+zhang}, \cite{Im-Si} 
and \cite{sen} for nonlocal versions). 
Moreover, through studying PDE \eqref{spde2} with rough drift $b$, 
the strong and weak well-posedness of SDE \eqref{SDE0} with irregular drift $b$ was also studied in \cite{Ch, Zh9, Wa-Zh}.

\smallskip

On the other hand, SKE \eqref{spde2} has a close connection with the filtering problem associated with SDE \eqref{SDE0}.
In fact, let us consider the following simple stochastic Langevin equation
$$
\dif X_t=V_t\dif t,\quad \dif V_t=\dif B_t-\dif W_t,
$$
where $B$ and $W$ are two independent Brownian motions. In this model, $W$ is regarded as the observable noise, 
$B$ is the hidden noise, $(X,V)$ stands for the position and velocity of particles, whose distribution needs to be predicted. 
Let $\pi(t,\omega,x,v)$ be the conditional probability density of $\mP((X_t,V_t)\in \cdot|\sF^W_t)$, i.e.,
for any $\varphi\in L^\infty(\mR^{2d})$,
$$
\int_{\mR^{2d}}\pi(t,\omega,x,v)\varphi(x,v)\dif x\dif v=\mE(\varphi(X_t,V_t)|\sF^W_t)(\omega),
$$
where $\sF^W_t$ is the $\sigma$-algebra generated by $\{W_s,s\leq t\}$.
It is well known that (see \cite{roz})
\begin{eqnarray}
(\partial_t+v\cdot\nabla_x)\pi=\Delta_v\pi\dif t+\nabla_{v}\pi \dif W_t.\label{spde0}
\end{eqnarray}
This is just a model equation of SKE \eqref{spde2}. 
Note that the nonlinear filtering problem of SDE \eqref{SDE0} is naturally related to SKE \eqref{spde22} (see Section 6 below).

\smallskip

The study of SPDEs has a long history since the earlier work of Pardoux \cite{Pa} about the SPDEs in Hilbert spaces, 
Krylov and Rozovskii \cite{Kr-Ro} about the SPDEs in the framework of Gelfand triple, and Walsh \cite{Wa} about the study of stochastic wave equations. 
Nowadays, the theory of SPDEs has been greatly developed, and there are vast literatures, see, for examples, 
the monographs \cite{Da-Za, roz, wei+roc, kry99}, etc.. 
Here we only mention parts of the related work.
In the case of smooth coefficients, Rozovsky and Lototsky \cite{roz} systematically developed a complete linear $L^2$-theory about
parabolic SPDEs, especially, addressed the application in nonlinear filtering problems (see also \cite{zxc} for a non-smooth, degenerate $L^2$-theory). 
A complete $L^p$-theory about semi-linear SPDEs was established by Krylov \cite{kry99}.
Here by a theory, according to Krylov's words \cite{kry99}, we mean not only results that, for $f, g^k$ belonging to a space $\mS$,
the solution belongs to some stochastic spaces $\mS'$, but also that  every element of $\mS'$ 
can be obtained as a solution for certain $f, g^k$ belonging to the same $\mS$. 
In other words, we have a bijection $F: \mS\to\mS'$ between $f, g^k$ and the solutions.
 In \cite{mik},
Mikulevicius obtained the Schauder estimates for nondegenerate SPDEs when the leading coefficient $a$
is nonrandom and $\sigma=0$. The full Schauder theory for nondegenerate second order SPDEs 
with random coefficients was established recently by Du and Liu \cite{du+liu}.
The main difficulty, when we consider the random leading coefficients, is that the direct Duhamel formulation 
is no longer applicable due to the non-adaptedness of the integrands in the stochastic integrals since we have random kernels in this case.
In \cite{du+liu}, the authors adopt Campanato's energy characterization of H\"older spaces to overcome this difficulty.

\smallskip

The aim of this paper is to establish the optimal regularity estimates for SKE \eqref{spde2}
in anisotropic Besov spaces under some H\"older regularity assumptions on $a,b,\sigma$.
Our optimal regularity estimates are not only for nondegenerate velocity component $v$, but also for degenerate 
position component $x$. As discussed above, when $a$ is nonrandom, by Duhamel's formulation, and using completely 
the same argument developed in \cite{sen}, one can establish a satisfactory Schauder theory.
However, for random leading coefficient $a$, the method adopted by Du and Liu \cite{du+liu}, if it is not impossible,
seems hard to be used for SKE \eqref{spde2}.
We mention that when $d=1$,  Pascucci and Pesce \cite{APAPesce2020} obtained the Schauder estimates for SKE \eqref{spde2}
by using It\^o-Wentzell's formula to reduce SKE \eqref{spde2} into a deterministic PDE
with random coefficients. Thus the price to pay is that the coefficients have to be at least $C^3$-differentiable in $x,v$.
Moreover, we are also interested in the following nonlinear SKE with super-linear growth coefficient in $\mR^2$:
\begin{align}\label{EZ1}
\dif u=[\Delta_v u+v\cdot\nabla_x u]\dif t+|u|^{\gamma+1} \dif B,\ \gamma\in[0,\tfrac18),
\end{align}
where $\dif B$ stands for the velocity-time white noise.
In particular, it includes a kinetic version of continuous Parabolic Anderson Model (abbreviated as PAM) when $\gamma=0$ and $u\geq 0$.
We referred to \cite{Ba-Ch}  and references therein for some background about PAM.
In the nondegenerate case, i.e., $u$ does not depend on $x$, 
the global well-posedness to the above super-linear growth SKEs with Dirichlet boundary conditions for $\gamma\in[0,\frac12)$
was first established by Mueller \cite{Mu}. His proof is based on the large deviation estimates and only for positive solutions.
While in \cite[Section 8.4]{kry99}, Krylov provides a quite different and more direct proof, and still for positive solutions,
where the key observation is that  for any bounded stopping time $\tau$,
$$
\mE\int_{\mR}u(\tau,v)\dif v\leq\int_{\mR}u(0,v)\dif v.
$$
Here we aim to show the existence and uniqueness of weak solutions 
and the optimal regularity  of $u$ in $t,x,v$ to nonlinear SKE \eqref{EZ1}. As above, the key point is to show
$$
\mE\int_{\mR^2}|u(\tau,x,v)|\dif x\dif v\leq\int_{\mR^2}|u(0,x,v)|\dif x\dif v.
$$
See Lemma \ref{78} below.

\smallskip

Now we describe the strategies of proving the optimal Besov regularity to SKE \eqref{spde2}. 
First of all, we study the model equation
\begin{align}\label{spde0}
\dif u=[\Delta_v u+v\cdot\nabla_x u+f]\dif t+g^k\dif W^k_t.
\end{align}
Using some tools from \cite{sen}, we establish the optimal regularity estimates of the solution
in anisotropic Besov spaces  with respect to the time and spatial-velocity variables by Duhamel's formula. 
We would like to emphasize that the semigroup $P_t$ associated with $\Delta_v+v\cdot\nabla_x$ 
behaves unlike the Gaussian heat semigroup. One has to consider the transport semigroup $\Gamma_tf(x,v):=f(x+tv,v)$.
The non-commutativity $\Gamma_tP_t\not=\Gamma_tP_t$ brings us many difficulties.
Moreover, although $P_t$ is a strongly continuous semigroup in $L^p$-space, 
there seems not exist a good characterization for the domain of $\Delta_v+v\cdot\nabla_x$ in $L^p$-space. 
Thus, SKE \eqref{spde0} does not fall into the abstract framework studied in \cite{Da-Za}. We have to carefully handle 
the anisotropicity caused by degenerate term $v\cdot\nabla_x$.
Next, by a generalized 
It\^o-Wentzell's formula (\cite{kry99, Krylov11}),
we transform SKE \eqref{spde2} with random but constant $a,\sigma$ into model equation \eqref{spde0}, and then
obtain the regularity estimates for \eqref{spde2}, where the key point is that we need to work in Besov spaces 
with  finite integrability exponent.
Finally, we shall use the freezing coefficient argument to derive the optimal regularity estimates for variable random coefficients $a,\sigma$.
In order to make the perturbation term can be absorbed by the freezing term, we introduce a new localized
anisotropic Besov norm $\nor\cdot\nor_{\wt\bB^s_{p;\theta}}$ in \eqref{Se0} below. Unlike $L^\infty\not\subset L^p$ for $p<\infty$,
one obvious advantage of using this localized norm is that for any $s>0$,
$$
\nor\cdot\nor_{\wt\bB^s_{p;\theta}}\leq C \nor\cdot\nor_{\wt\bB^s_{p';\theta}},\ \ p'\geq p.
$$
Thus our initial value can be a constant. Such a norm was also used in \cite{Zh-Zh21}. 
It is noticed that if the coefficients $a,b,\sigma$ and $f,g$, 
initial value $u_0$ do not depend on the position variable $x$, then 
SKEs \eqref{spde2} and \eqref{spde22} naturally reduce to the classical SPDEs studied in \cite{kry99}.
To the authors' knowledge, even in this classical case, our Besov estimates are also new.

\smallskip

This paper is organized as follows: In Section 2, we introduce some anisotropic function spaces, and prepare some useful results for later use.  
In Section 3, we establish the optimal regularity estimates for random constant coefficients case, 
namely, $a$ and $\sigma^k$ are independent of space and velocity variables. In Sections 4 and 5, we prove our main results 
(see Theorems \ref{Th41},\ref{es},\ref{Th45}). 
In Section 6, we apply our results to the filtering problem of degenerate diffusion processes, and show the existence and regularity of
conditional density processes, which satisfies a nonlinear SPDE in classical sense under some assumptions.
In Section 7, we also show a well-posedness result for a class of nonlinear SKEs driven by velocity-time white noises, and obtain
the optimal Besov regularity, which seems to be new even in the classical nondegenerate case (cf. \cite{Wa, kry99}).

\smallskip

Throughout this paper, we use the following conventions: The letter $C=C(\cdots)$ denotes an unimportant constant, 
whose value may change in different places, and which is increasing with  respect to its arguments.
We use $A\asymp B$ and $A\lesssim B$  to denote $C^{-1}B\leq A\leq CB$ and
$A\leq C B$, respectively, for some unimportant constant $C\geq 1$. As usual, we use $:=$ or $=:$ as a way of definition, and for any $a,b\in\mR$,
$$
a\wedge b:=\min(a,b),\ \ a\vee b:=\max(a,b),\ \ a^+:=a\vee 0,\ \ a^-:=-(a\wedge 0).
$$
Moreover, for $p\in[1,\infty]$, we use $\ell^p$ to denote the usual space of a sequence of real numbers that is $p$-order summable.

\section{Preliminaries}
\label{prelimi}
\subsection{Vector-valued anisotropic Besov and H\"older spaces}
In this subsection we introduce the vector-valued anisotropic H\"older and Besov spaces and their basic properties for later use.
Let $m=(m_1,\cdots,m_n)\in\mathbb{N}^n$ with $m_1+\cdots+m_n=N$ and $\theta=(\theta_1,\cdots,\theta_n)\in [1,\infty)^n$
be fixed. We denote $\theta\cdot m=\theta_1m_1+\cdots+\theta_nm_n$.
For $x=(x_1,\cdots,x_n)$ and $y=(y_1,\cdots,y_n)\in\mathbb{R}^{m_1}\times\cdots\times\mathbb{R}^{m_n}$,
we introduce the following distance in $\mathbb{R}^N$ by
$$
|x-y|_\theta:=\sum^{n}_{i=1}|x_i-y_i|^{1/\theta_i},\quad x_i,y_i\in\mathbb{R}^{m_i}.
$$

Let $\mB$ be a Banach space. For $p\in[1,\infty]$, let $\mL^p_x(\mB):=L^p(\mR^N,\dif x;\mB)$
be the usual vector-valued $L^p$-space over $\mR^N$. For $h\in\mR^N$ and a map $f:\mR^N\to\mB$, the first order difference operator is defined by
$$
\delta^{(1)}_hf(x):=f(x+h)-f(x),
$$
and for $M\in\mN$, the $M$-order difference operator is defined recursively by
$$
\delta^{(M)}_hf(x):= \delta^{(1)}_h\delta^{(M-1)}_hf(x).
$$
By induction, it is easy to see that
\begin{align}\label{Def8}
\delta^{(M)}_hf(x)=\sum^{M}_{k=0}(-1)^{M-k}{M\choose k} f(x+kh),\quad  h\in\mathbb{R}^{N},
\end{align}
where $\binom M k$ is the binomial coefficient.

\begin{definition}[Vector-valued anisotropic H\"older spaces]
For $s\in(0,\infty)$, the $\mB$-valued anisotropic H\"older space is defined by norm
$$
\|f\|_{\mathbf{C}^{s}_{\theta}(\mB)}:=\|f\|_{\mL^\infty_x(\mB)}+[f]_{\mathbf{C}^{s}_{\theta}(\mB)}<\infty,
$$
where
$$
[f]_{\mathbf{C}^{s}_{\theta}(\mB)}
:=\sup_{h}\|\delta^{([s]+1)}_hf\|_{\mL^\infty_x(\mB)}/|h|^s_\theta.
$$
Here $[s]$ denotes the greatest integer less than $s$.
\end{definition}

In order to introduce the anisotropic Besov space, we need a symmetric 
nonnegative $C^{\infty}-$function $\phi^\theta_0$ on $\mathbb{R}^N$ with
$$
\phi^\theta_0(\xi)=1\ \mathrm{for}\ \xi\in B^\theta_1\ \mathrm{and}\ \phi^\theta_0(\xi)=0\ \mathrm{for}\ \xi\notin B^\theta_2,
$$
where for $r>0$ and $y\in\mR^N$,
$$
B^\theta_r(y):=\{x\in\mathbb{R}^N:|x-y|_\theta\leq r\},\ \ B^\theta_r:=B^\theta_r(0).
$$
For $\xi=(\xi_1,\cdots,\xi_n)\in\mathbb{R}^{m_1}\times\cdots\times\mathbb{R}^{m_n}$ and $j\in\mathbb{N}$, we define
$$
\phi^\theta_j(\xi):=\phi^\theta_0(2^{-\theta j}\xi)-\phi^\theta_0(2^{-\theta(j-1)}\xi),
$$
where
$$
2^{-\theta j}\xi:=(2^{-\theta_1 j}\xi_1,\cdots,2^{-\theta_n j}\xi_n).
$$
By the very definition, one sees that for $j\in\mathbb{N}$, $\phi^\theta_j(\xi)=\phi^\theta_1(2^{-\theta(j-1)}\xi)\geq 0$ and
$$
\mathrm{supp}\,\phi^\theta_j\subset B^\theta_{2^{j+1}}\setminus B^\theta_{2^{j-1}},\quad\sum^n_{j=0}\phi^\theta_j(\xi)=\phi^\theta_0(2^{-n}\xi)\to 1,\quad n\to\infty.
$$
For $f\in \mL^1_x(\mB)$, let $\hat f$ be the Fourier transform of $f$ defined by
$$
\hat f(\xi):=(2\pi)^{-d/2}\int_{\mR^N} \e^{-{\rm i}\xi\cdot x}f(x)\dif x, \quad\xi\in\mR^N,
$$
and $\check f$ the Fourier  inverse transform of $f$ defined by
$$
\check f(x):=(2\pi)^{-d/2}\int_{\mR^N} \e^{{\rm i}\xi\cdot x}f(\xi)\dif\xi, \quad x\in\mR^N.
$$

Let $\cS$ be the space of all Schwartz functions on $\mR^N$ and $\cS'(\mB)$ the space of all continuous linear operators from
$\cS$ to $\mB$, called vector-valued distribution space. 
Note that $\cS':=\cS'(\mR)$ is just the tempered distribution space.
 In a natural way, $\mL^1_x(\mB)\subset\cS'(\mB)$, and we can extend the Fourier 
transform to the element in $\cS'(\mB)$ by duality. 
For given $j\in\mathbb{N}_0$, the block operator  $\mathcal{R}^\theta_j$ is defined on $\cS'(\mB)$ by (see \cite{fn})
\begin{align}\label{Ph0}
\mathcal{R}^\theta_jf(x):=(\phi^\theta_j\hat{f})\check{\ }(x)=\check{\phi}^\theta_j*f(x)
=2^{\theta\cdot m(j-1)}\int_{\mathbb{R}^N}\check{\phi}^\theta_1(2^{\theta(j-1)}y)f(x-y)\dif y,
\end{align}
where the convolution is understood in the distributional sense. In particular, by the symmetry of $\phi^\theta_j$, we have
$$
\<\cR^\theta_j f,g\>=\< f,\cR^\theta_jg\>,\ \ f\in\cS', \ g\in\cS,
$$
where $\<\cdot,\cdot\>$ stands for the dual pair between $\cS'$ and $\cS$.  

\begin{definition}[Vector-valued anisotropic Besov spaces]\label{bs}
For $s\in\mR$ and $p\in[1,\infty]$, the $\mB$-valued anisotropic Besov space is defined by
$$
\mathbf{B}^{s}_{p;\theta}(\mB):=\left\{f\in \cS'(\mB): \|f\|_{\mathbf{B}^{s}_{p;\theta}(\mB)}
:=\sup_{j\geq0}\big(2^{sj}\|\mathcal{R}^\theta_{j}f\|_{\mL^p_x(\mB)}\big)<\infty\right\}.
$$
When $\theta=(1,\cdots,1)$, we shall simply write 
$$
\mathbf{B}^{s}_{p}(\mB):=\mathbf{B}^{s}_{p;\theta}(\mB),\ \cR_j:=\cR^\theta_j,\ \phi_j:=\phi^\theta_j.
$$
When $\mB=\mR$, we shall simply write $\mathbf{B}^{s}_{p;\theta}:=\mathbf{B}^{s}_{p;\theta}(\mR)$.
\end{definition}

\br\rm
In the literature, there are usually two subscripts $p,q$ in the definition of Besov spaces, where
$p$ stands for the integrability of spatial variables and $q$ denotes the $\ell^q$-norm of frequency index $j$.
Since we only use the Besov space of $q=\infty$, in our definition, we take $q=\infty$ and drop it for simplicity.
\er

\br\rm
In application below, we usually take $\mB=\mL^p_\omega:=L^p(\Omega,\sF,\mP)$.
In this case, by Fubini's theorem, for any $s<s'$ and $p\in[1,\infty)$, one has
\begin{align}\label{CT1}
\mathbf{B}^{s'}_{p;\theta}(\mL^p_\omega)\subset \mL^p_\omega(\mathbf{B}^{s}_{p;\theta})\subset \mathbf{B}^s_{p;\theta}(\mL^p_\omega).
\end{align}
Indeed, by definition and Fubini's theorem,
\begin{align*}
\|f\|^p_{\mL^p_\omega(\mathbf{B}^{s}_{p;\theta})}&=\mE \|f(\omega,\cdot)\|^p_{\mathbf{B}^{s}_{p;\theta}}
=\mE \left(\sup_{j\geq 0}2^{sjp}\|\cR^\theta_jf(\omega,\cdot)\|^p_{\mL^p_x}\right)\\
&\leq\sum_{j\geq 0}2^{sj p}\mE \left(\|\cR^\theta_jf(\omega,\cdot)\|^p_{\mL^p_x}\right)
=\sum_{j\geq 0}2^{sjp}\|\cR^\theta_jf\|^p_{\mL^p_x(\mL^p_\omega)}\\
&\leq\sum_{j\geq 0}2^{(s-s')jp}\|f\|^p_{\mathbf{B}^{s'}_{p;\theta}(\mL^p_\omega)}
=1/(1-2^{(s-s')p})\|f\|^p_{\mathbf{B}^{s'}_{p;\theta}(\mL^p_\omega)},
\end{align*}
and
\begin{align*}
\|f\|^p_{\mathbf{B}^{s}_{p;\theta}(\mL^p_\omega)}&=\sup_{j\geq 0}2^{sjp}\|\cR^\theta_jf\|^p_{\mL^p_x(\mL^p_\omega)}
\leq\mE \left(\sup_{j\geq 0}2^{sjp}\|\cR^\theta_jf(\omega,\cdot)\|^p_{\mL^p_x}\right)=\|f\|^p_{\mL^p_\omega(\mathbf{B}^{s}_{p;\theta})}.
\end{align*}
\er

We recall the following result whose proof is completely the same as in 
\cite{fn}*{p.52, Lemma 2.1}. We omit the details.
\bl
[Bernstein's type inequalities]
For any $1\leq p\leq q\leq\infty$, $k\in\mN_0$ and $i=1,\cdots,n$, there is a constant
$C=C(\theta,m,p,q,k,i)>0$
 such that for all $j\geq 0$,
\begin{align}\label{Ber}
\|\nabla^k_{x_i}\cR^\theta_j f\|_{\mL^q_x(\mB)}\lesssim_C 2^{j(\theta_ik +(\theta\cdot m)(\frac1p-\frac1q))}\|\cR^\theta_j f\|_{\mL^p_x(\mB)},
\end{align}
where $\nabla^k_{x_i}$ denotes the $k$-order gradient with respect to $x_i$.
\el

As easy consequences of the above lemma, we have
\bl
Let $1\leq p\leq q\leq\infty$, $s\in\mR$ and $k\in\mN_0$, $i=1,\cdots,n$.
\begin{enumerate}[(i)]
\item If $0<\theta_ik+(\theta\cdot m)(\frac1p-\frac1q)<s$, then for some $C=C(\theta,m,p,q,s,k,i)>0$,
\begin{align}\label{Ber1}
\|\nabla^k_{x_i} f\|_{\mL^q_x(\mB)}\lesssim_C\|f\|_{\bB^{s}_{p;\theta}(\mB)},
\end{align}
and for $s'=s+\theta_ik +(\theta\cdot m)(\frac1p-\frac1q)$,
\begin{align}\label{Ber0}
\|\nabla^k_{x_i}f\|_{\bB^{s}_{q;\theta}(\mB)}\lesssim_C \|f\|_{\bB^{s'}_{p;\theta}(\mB)}.
\end{align}
\item For any $s_1<s<s_2$, there is a $C=C(\theta,m,p,s,s_1,s_2)>0$ such that
\begin{align}\label{DQ9}
\|f\|_{\bB^{s}_{p;\theta}(\mB)}\lesssim_C \|f\|^{(s-s_1)/(s_2-s_1)}_{\bB^{s_2}_{p;\theta}(\mB)}\|f\|^{(s_2-s)/(s_2-s_1)}_{\bB^{s_1}_{p;\theta}(\mB)}.
\end{align}
\end{enumerate}
\el
\begin{proof}
(i) Noting that $f=\sum_j\cR^\theta_j f$, by Bernstein's inequality \eqref{Ber}, we have
\begin{align*}
\|\nabla^k_{x_i}f\|_{\mL^q_x(\mB)}&\leq\sum_{j\geq 0}\|\nabla^k_{x_i}\cR^\theta_j f\|_{\mL^q_x(\mB)}
\leq\sum_{j\geq 0}2^{j(\theta_i k+(\theta\cdot m)(\frac1p-\frac1q))}\|\cR^\theta_j f\|_{\mL^p_x(\mB)}\\
&\lesssim\sum_{j\geq 0}2^{j (\theta_i k+(\theta\cdot m)(\frac1p-\frac1q))-js}\|f\|_{\bB^s_{p;\theta}(\mB)}\lesssim\|f\|_{\bB^s_{p;\theta}(\mB)}.
\end{align*}
Thus we get \eqref{Ber1}. Estimate \eqref{Ber0} is direct by \eqref{Ber} and $\cR^\theta_j\nabla^k_{x_i}=\nabla^k_{x_i}\cR^\theta_j$.

(ii) It follows by definition.
\end{proof}

The following lemma is elementary.
\bl\label{Le27}
For any $0<\beta<\alpha<\infty$, it holds that for all $\lambda>0$,
$$
\sum_{j\geq 0}(\lambda2^{\alpha j}\wedge 1)2^{-\beta j}\leq 
\left(\frac{2^\beta}{\ln 2}\int^\infty_0(r^\alpha\wedge 1)r^{-\beta-1}\dif r\right)\lambda^{\frac\beta\alpha}.
$$
\el
\begin{proof}
Note that
\begin{align*}
\sum_{j\geq 0}(\lambda2^{\alpha j}\wedge 1)2^{-\beta j}\leq 2^\beta\int^\infty_0(\lambda2^{\alpha s}\wedge 1)2^{-\beta s}\dif s
=\frac{2^\beta\lambda^{\frac\beta\alpha}}{\ln 2}\int^\infty_{\lambda^{1/\alpha}}(r^\alpha\wedge 1)r^{-\beta-1}\dif r,
\end{align*}
which gives the desired estimate.
\end{proof}

Now we show the following important characterization of anisotropic Besov space $\bB^s_{p;\theta}(\mB)$
(cf. \cite{Triebel3}).
\bl
For any $s\in(0,\infty)$ and $p\in[1,\infty]$,
there exists a constant $C=C(\theta,m,p,s)\geq1$ such that for all $f\in\bB^s_{p;\theta}(\mB)$,
\begin{align}\label{FG1}
\|f\|_{\bB^s_{p;\theta}(\mB)}\asymp_C \sup_{h\in\mR^N}\Big(|h|^{-s}_\theta\|\delta^{([s]+1)}_hf\|_{\mL^p_x(\mB)}\Big)+\|f\|_{\mL^p_x(\mB)}.
\end{align}
In particular, $\bC^s_\theta(\mB)=\bB^s_{\infty;\theta}(\mB)$ and
\begin{align}\label{28}
\|f \|_{\bB^{s}_{p;\theta}(\mB)}\asymp_C\|f \|_{\bB^{s/\theta_1}_{p;x_1}(\mB)}+\cdots+\|f \|_{\bB^{s/\theta_n}_{p;x_n}(\mB)},
\end{align}
where for $i=1,\cdots,n$,
$$
\|f \|_{\bB^s_{p;x_i}(\mB)}:=\sup_{j\geq0}\left[2^{sj}\left(\int_{\mR^N}\|\cR^{x_i}_j f(x)\|^p_{\mB}\dif x\right)^{\frac1p}\right],
$$
and for $x=(x_1,\cdots,x_n)$,
\begin{align}\label{DD9}
\cR^{x_i}_j f(x):=\mathcal{R}_{j}f(\cdots, x_{i-1},\cdot,x_{i+1},\cdots)(x_i).
\end{align}
\el
\begin{proof}
(i) For simplicity, we set $M:=[s]+1$. We first prove that for any $h\in\mR^N$,
$$
\|\delta^{(M)}_{h}f\|_{\mL^p_x(\mB)}\lesssim |h|^{s}_\theta\|f\|_{\bB^s_{p;\theta}(\mB)}.
$$
By definition, without loss of generality, it suffices to prove that for $h_1\in\mR^{m_1}$,
\begin{align}\label{DX102}
\|\delta^{(M)}_{{\bf h}_1}f\|_{\mL^p_x(\mB)}\lesssim |h_1|^{s/\theta_1}\|f\|_{\bB^s_{p;\theta}(\mB)},
\end{align}
where ${\bf h}_1:=(h_1,0,\cdots,0)\in\mR^N$.
Since $\cR^\theta_j\delta^{(M)}_{{\bf h}_1}f=\delta^{(M)}_{{\bf h}_1}\cR^\theta_jf$, by \eqref{Def8}
with Taylor's expansion and Bernstein's inequality \eqref{Ber}, we have
\begin{align*}
\|\cR^\theta_j\delta^{(M)}_{{\bf h}_1}f\|_{\mL^p_x(\mB)}\lesssim |h_1|^{M}\|\nabla^{M}_{x_1}\cR^\theta_j f\|_{\mL^p_x(\mB)}
\lesssim |h_1|^{M} 2^{(M\theta_1-s)j}\|f\|_{\bB^s_{p;\theta}(\mB)},
\end{align*}
and also,
$$
\|\cR^\theta_j\delta^{(M)}_{{\bf h}_1}f\|_{\mL^p_x(\mB)}\leq \sum_{k=0}^M{M\choose k} \|\cR^\theta_j f\|_{\mL^p_x(\mB)}
\leq 2^M 2^{-sj}\|f\|_{\bB^s_{p;\theta}(\mB)}.
$$
Hence,
$$
\|\delta^{(M)}_{{\bf h}_1}f\|_{\mL^p_x(\mB)}\leq\sum_{j\geq 0}\|\cR^\theta_j\delta^{(M)}_{{\bf h}_1}f\|_{\mL^p_x(\mB)}
\lesssim\sum_{j\geq 0}\Big[(|h_1|^{M}2^{(M\theta_1-s)j})\wedge 2^{-sj}\Big]\|f\|_{\bB^s_{p;\theta}(\mB)},
$$
which gives \eqref{DX102} by Lemma \ref{Le27}.

(ii) For $j\geq 1$, since $\int_{\mathbb{R}^N} \check{\phi}^\theta_j(h)\dif h=(2\pi)^{d/2}\phi^\theta_j(0)=0$, by \eqref{Def8} and the change of variable, we have
\begin{align*}
\int_{\mathbb{R}^N} \check{\phi}^\theta_j(h)\delta^{(M)}_hf(x)\dif h
&=\sum^{M}_{k=0}(-1)^{M-k}{M\choose k}\int_{\mathbb{R}^N} \check{\phi}^\theta_j(h)f(x+kh)\dif h\\
&=\sum^{M}_{k=1}(-1)^{M-k}{M\choose k}\int_{\mathbb{R}^N} \check{\phi}^\theta_j(h)f(x+kh)\dif h\\
&=\sum^{M}_{k=1}(-1)^{M-k}{M\choose k}\int_{\mathbb{R}^N} \phi^\theta_{j}(k\cdot)^{\check{ }}(h)f(x+h)\dif h.
\end{align*}
In particular, if we define for $j\in\mN_0$,
$$
\phi^{\theta, M}_j(\xi):=(-1)^{M+1}\sum^{M}_{k=1}(-1)^{M-k}{M\choose k} \phi^{\theta}_{j}(k \xi),
$$
then
$$
(-1)^{M+1}\int_{\mathbb{R}^N} \check{\phi}^\theta_j(h)\delta^{(M)}_hf(x)\dif h=[\phi^{\theta, M}_j]^{\check{}}*f(x)=:\mathcal{R}^{\theta, M}_jf(x),
$$
and for $j\geq 1$,
\begin{align}
\|\cR^{\theta, M}_j f\|_{\mL^p_x(\mB)}&\leq\int_{\mR^N} |\check{\phi}^\theta_j(h)|\,\|\delta^{(M)}_hf\|_{\mL^p_x(\mB)}\dif h\no\\
&\leq\left(\int_{\mR^N} |\check{\phi}^\theta_j(h)||h|^s_\theta\dif h\right)\sup_h\frac{\|\delta^{(M)}_hf\|_{\mL^p_x(\mB)}}{|h|^s_\theta}\no\\
&\lesssim2^{-sj}\sup_h\Big(|h|^{-s}_\theta\|\delta^{(M)}_hf\|_{\mL^p_x(\mB)}\Big).\label{Ea0}
\end{align}
On the other hand, noting that
$$
\phi^{\theta, M}_j(\xi)= \phi^{\theta, M}_0(2^{-\theta j}\xi)-\phi^{\theta, M}_0(2^{-\theta(j-1)}\xi)
$$
and
$$
\phi^{\theta, M}_0(\xi)=1\mbox{ for $\xi\in B^\theta_{1/M}$  and }
\phi^{\theta, M}_0(\xi)=0\mbox{ for $\xi\notin B^\theta_{2}$},
$$
we have
$$
\mathrm{supp}\ \phi^{\theta, M}_j\subset B^{\theta}_{2^{j+1}} \setminus B^{\theta}_{(2^{j-1})/M}.
$$
Hence, for any $i,j\in\mN_0$ with $|j-i|>\log_2 M+2=:\gamma$,
$$
\mathcal{R}^{\theta,M}_{j} \mathcal{R}^{\theta}_{i}f(x)=0.
$$
Moreover, noting that for any $\xi\in\mathbb{R}^N$,
\begin{align*}
\sum_{j\geq0}\phi^{\theta, M}_j(\xi)
&= (-1)^{M+1}\sum^{M}_{k=1} \sum_{j\geq0}(-1)^{M-k}{M\choose k} \phi^{\theta}_{j}(k\xi)\\
&= (-1)^{M+1}\sum^{M}_{k=1}(-1)^{M-k}{M\choose k}=1,
\end{align*}
we have
\begin{align*}
\cR^{\theta}_i f=\sum_{j\geq 0}\cR^{\theta}_i\cR^{\theta, M}_j f=\sum_{|j-i|\leq\gamma}\cR^{\theta}_i\cR^{\theta, M}_j f.
\end{align*}
Therefore, for $i\geq\gamma+1$, by \eqref{Ea0},
\begin{align*}
 \|\cR^{\theta}_i f\|_{\mL^p_x(\mB)}
 &\leq \sum_{|i-j|\leq \gamma}\| \cR^{\theta}_i\cR^{\theta, M}_j f\|_{\mL^p_x(\mB)}
 \lesssim \sum_{|i-j|\leq \gamma}\|\cR^{\theta, M}_j f\|_{\mL^p_x(\mB)}\\
&\lesssim \sum_{|i-j|\leq \gamma} 2^{-sj} \sup_h\Big(|h|^{-s}_\theta\|\delta^{(M)}_hf\|_{\mL^p_x(\mB)}\Big) \\
&\lesssim  2^{-si} \sup_h\Big(|h|^{-s}_\theta\|\delta^{(M)}_hf\|_{\mL^p_x(\mB)}\Big).
\end{align*}
For $i<\gamma+1$, we always have
$$
 \|\cR^{\theta}_i f\|_{\mL^p_x(\mB)}\lesssim \|f\|_{\mL^p_x(\mB)}.
$$
Thus we obtain another side estimate and complete the proof.
\end{proof}

For a Banach space $\mB$, let $\cL_\mB$ be the Banach space of all bounded linear operators from $\mB$  to $\mB$.
\bl
For any $s\in(0,2)$ and $p\in[1,\infty]$, there exists a constant
$C=C(\theta,m,p,s)>0$
 such that for
all $\cT\in \bC^s_\theta(\cL_\mB)$ and $g\in \bB^s_{p;\theta}(\mB)$,
\begin{align}\label{24}
\|\cT g\|_{\bB^{s}_{p;\theta}(\mB)}\lesssim_C
\left\{
\begin{aligned}
&\|\cT\|_{\mathbb{L}^\infty_x(\cL_\mB)}\|g\|_{\bB^s_{p;\theta}(\mB)}+\|\cT\|_{\bC^s_\theta(\cL_\mB)}\|g\|_{\mL^p_x(\mB)},&s\in(0,1),\\
&\|\cT\|_{\mathbb{L}^\infty_x(\cL_\mB)}\|g\|_{\bB^s_{p;\theta}(\mB)}+\|\cT\|_{\bC^s_\theta(\cL_\mB)}\|g\|_{\bB^{s/2}_{p;\theta}(\mB)},& s\in[1,2),\\
\end{aligned}
\right.
\end{align}
where $(\cT g)(x):=\cT(x)g(x)$ for $x\in\mR^N$.
\el
\begin{proof}
For $s\in(0,1)$, by definition and \eqref{FG1}, we have
\begin{align*}
\|\delta^{(1)}_h(\cT g)\|_{\mL^p_x(\mB)}
&\leq \|\cT\|_{\mathbb{L}^\infty_x(\cL_\mB)}\|\delta^{(1)}_h g\|_{\mL^p_x(\mB)}
+\|\delta^{(1)}_h \cT\|_{\mathbb{L}^\infty_x(\cL_\mB)}\|g\|_{\mL^p_x(\mB)}\\
&\lesssim |h|^s_\theta\|\cT\|_{\mathbb{L}^\infty_x(\cL_\mB)}\|g\|_{\bB^s_{p;\theta}(\mB)}
+|h|^s_\theta\|\cT\|_{\bC^s_\theta (\cL_\mB)}\|g\|_{\mL^p_x(\mB)}.
\end{align*}
For $s\in[1,2)$, noting that
\begin{align*}
\delta^{(2)}_h(\cT g)=\cT(\cdot+2h)\delta^{(2)}_h g+\delta^{(2)}_h \cT g(\cdot+h)+\delta^{(1)}_{2h}\cT\,\delta^{(1)}_h g,
\end{align*}
by \eqref{FG1}, we have
\begin{align*}
\|\delta^{(2)}_h(\cT g)\|_{\mL^p_x(\mB)}
&\leq \|\cT\|_{\mathbb{L}^\infty_x(\cL_\mB)}\|\delta^{(2)}_h g\|_{\mL^p_x(\mB)}
+\|\delta^{(2)}_h \cT\|_{\mathbb{L}^\infty_x(\cL_\mB)}\|g\|_{\mL^p_x(\mB)}\\
&\quad+\|\delta^{(1)}_{2h}\cT\| _{\mathbb{L}^\infty_x(\cL_\mB)}\|\delta^{(1)}_{h} g\|_{\mL^p_x(\mB)}\\
&\lesssim |h|^s_\theta\|\cT\|_{\mathbb{L}^\infty_x(\cL_\mB)}\|g\|_{\bB^s_{p;\theta}(\mB)}
+|h|^s_\theta\|\cT\|_{\bC^s_\theta (\cL_\mB)}\|g\|_{\mL^p_x(\mB)}\\
&\quad+|2h|^{s/2}_\theta|h|^{s/2}_\theta\|\cT\|_{\bC^{s/2}_\theta (\cL_\mB)}\|g\|_{\bB^{s/2}_{p;\theta}(\mB)}.
\end{align*}
The estimate \eqref{24} follows by $\|\cT g\|_{\mL^p_x(\mB)}\leq\|\cT\|_{\mathbb{L}^\infty_x(\cL_\mB)}\|g\|_{\mL^p_x(\mB)}$ and \eqref{FG1} again.
\end{proof}
\br\rm
For any $s\in\mR$ and $s'>|s|$, we in fact have (cf. \cite{fn})
\begin{align}\label{Ea1}
\|\cT g\|_{\bB^{s}_{p;\theta}(\mB)}\lesssim_C \|\cT\|_{\bC^{s'}_\theta(\cL_\mB)}\|g\|_{\bB^{s}_{p;\theta}(\mB)},
\end{align}
where $C=C(\theta,m,p,s,s')>0$.
\er
The following commutator estimates will be used to improve the regularity in $x$.
\bl\label{Le28}
Let $p\in[1,\infty]$, $\gamma\in(0,1)$ and $s\in(0,2)$. For any $i=1,\cdots,n$, there is a constant
$C=C(\theta,m,p,s,\gamma,i)>0$
such that for all $\cT\in \bC^\gamma_{x_i}(\bC^s_\theta(\cL_\mB))$, $g\in \bB^s_{p;\theta}(\mB)$ and $j\geq 0$,
\begin{align}
\|[\cR^{x_i}_j,\cT]g\|_{\bB^s_{p;\theta}(\mB)}\lesssim_C
2^{-\gamma j}[\cT]_{\bC^\gamma_{x_i}(\bC^s_\theta(\cL_\mB))}\|g\|_{\bB^s_{p;\theta}(\mB)},\label{1}
\end{align}
where $[\cR^{x_i}_j,\cT]g:=\cR^{x_i}_j(\cT g)-\cT(\cR^{x_i}_jg)$ and $\cR^{x_i}_j$ is defined by \eqref{DD9}. Moreover,
\begin{align}
\|[\cR^{x_i}_j,\cT]g\|_{\mathbb{L}^p_x(\mB)}\lesssim_C 2^{-\gamma j}[\cT]_{\bC^\gamma_{x_i} (\cL_\mB)}\|g\|_{\mathbb{L}^p_x(\mB)}.\label{2}
\end{align}
\el
\begin{proof}
Without loss of generality, we assume $i=1$. For simplicity, we write
$$
x^*_1:=(x_2,\cdots,x_n),\ \ \tau^{(x_1)}_h g(x):=g(x_1-h,x_1^*),\ h\in\mR^{m_1},
$$
and
$$
\delta^{(x_1)}_{h}\cT(x):=\cT(x_1-h,x^*_1)-\cT(x_1,x^*_1).
$$
By definition we have
\begin{align*}
[\cR^{x_1}_j,\cT]g(x)=\int_{\mR^{m_1}}\check{\phi}_j(h)\,\delta^{(x_1)}_{h}\cT(x)\,\tau^{(x_1)}_h\! g(x)\dif h.
\end{align*}
Note that by \eqref{24},
\begin{align*}
\|\delta^{(x_1)}_{h}\cT\,\tau^{(x_1)}_h\! g\|_{\bB^s_{p;\theta}(\mB)}
&\lesssim 
\|\delta^{(x_1)}_{h}\cT\|_{\bC^s_\theta (\cL_\mB)}\|g\|_{\bB^{s}_{p;\theta}(\mB)}\lesssim|h|^\gamma\|\cT\|_{\bC^\gamma_{x_1}(\bC^s_\theta (\cL_\mB))}\|g\|_{\bB^s_{p;\theta}(\mB)}.
\end{align*}
Hence,
\begin{align*}
\|[\cR^{x_1}_j,\cT]g\|_{\bB^s_{p;\theta}(\mB)}&\lesssim \left(\int_{\mR^{m_1}}\check{\phi}_j(h)|h|^\gamma\dif h\right)
[\cT]_{\bC^\gamma_{x_1}(\bC^s_\theta (\cL_\mB))}\|g\|_{\bB^s_{p;\theta}(\mB)}\\
&\lesssim 2^{-\gamma j}[\cT]_{\bC^\gamma_{x_1}(\bC^s_\theta (\cL_\mB))}\|g\|_{\bB^s_{p;\theta}(\mB)}.
\end{align*}
Thus we get \eqref{1}. Similarly, we can prove \eqref{2}.
\end{proof}

Let $\chi:\mR^N\to [0,1]$ be a smooth function with $\chi=1$ on $B^\theta_1$ and $\chi=0$ on the complement of $B^\theta_2$.
For $\delta>0$ and $x_0\in\mR^N$, define
\begin{align}\label{Cu1}
\chi^\delta_{x_0}(x):=\chi((x-x_0)/\delta).
\end{align}
We also introduce the following localized Besov space (see \cite{Zh-Zh21}).
\bl
\label{3gf}
Let $p\in[1,\infty]$ and $s\in\mR$.
For any fixed $\delta,\delta'>0$, there is a constant
$C=C(\theta,m,p,s,\delta,\delta')\geq 1$ such that
\begin{align}\label{Se0}
\nor f\nor_{\wt\bB^{s}_{p;\theta}(\mB)}:=\sup_{x_0} \|\chi^\delta_{x_0} f\|_{\mathbf{B}^{s}_{p;\theta}(\mB)}
\asymp_C \sup_{x_0} \|\chi^{\delta'}_{x_0} f\|_{\mathbf{B}^{s}_{p;\theta}(\mB)}.
\end{align}
In particular, the following localized Besov space is independent of the choice of $\delta$:
$$
\wt\bB^{s}_{p;\theta}(\mB):=\big\{f\in\cS'(\mB): \nor f\nor_{\wt\bB^{s}_{p;\theta}(\mB)}<\infty\big\}.
$$
\el
\begin{proof}
Without loss of generality, we assume $\delta<\delta'$.
For fixed $x_0\in\mR^N$, note the support of $\chi^{\delta'}_{x_0}$ is contained in $B^\theta_{2\delta'}(x_0)$.
By a finite covering technique, one can find a number $M$ independent of $x_0$ and points $\{x_i, i=1,\cdots,M\}$ such that
$$
B^\theta_{2\delta'}(x_0)\subset\cup_{i=1}^M B^\theta_{\delta}(x_i).
$$
Let $\varphi_i, i=1,\cdots,M$ be the partition of unity associated with $\{B^\theta_{\delta}(x_i), i=1,\cdots,M\}$ so that
$$
\mbox{supp}(\varphi_i)\subset B^\theta_{\delta}(x_i),\ \ \sum_{i=1}^M\varphi_i(x)=1\ \mbox{ on } \cup_{i=1}^M B^\theta_{\delta}(x_i).
$$
Since $\chi^{\delta}_{x_i}=1$ on $B^\theta_\delta(x_i)$, by \eqref{Ea1} we have
\begin{align*}
\|\chi^{\delta'}_{x_0}f\|_{\mathbf{B}^{s}_{p;\theta}(\mB)}&\leq\sum_{i=1}^M\|\chi^{\delta'}_{x_0}\varphi_i f\|_{\mathbf{B}^{s}_{p;\theta}(\mB)}
=\sum_{i=1}^M\| \chi^{\delta'}_{x_0}\varphi_i \chi^{\delta}_{x_i} f\|_{\mathbf{B}^{s}_{p;\theta}(\mB)}\\
&\lesssim\sum_{i=1}^M \| \chi^{\delta'}_{x_0} \varphi_i\|_{\mathbf{C}^{|s|+1}_{\theta}}\| \chi^{\delta}_{x_i}f\|_{\mathbf{B}^{s}_{p;\theta}(\mB)}\lesssim \sup_{x_0} \|\chi^\delta_{x_0} f\|_{\mathbf{B}^{s}_{p;\theta}(\mB)}.
\end{align*}
The proof is complete.
\end{proof}
\br\rm
By \eqref{Ber0} and the definition, we have
$$
\nor \nabla^k_{x_i}f\nor_{\wt\bB^{s}_{p;\theta}(\mB)}\lesssim \nor f\nor_{\wt\bB^{s+k\theta_i}_{p;\theta}(\mB)},
$$
and  for $s>0$ and $1\leq p\leq p'\leq \infty$,
$$
\nor f\nor_{\wt\bB^{s}_{p;\theta}(\mB)}\lesssim \nor f\nor_{\wt\bB^{s}_{p';\theta}(\mB)}.
$$
\er
\subsection{It\^o-Wentzell's formula for distribution-valued processes}
In this subsection we recall the generalized It\^o-Wentzell's formula of distribution-valued It\^o's processes in \cite{Krylov11}.
Let $\mathbf{D}$ be the set of all $\cS'$-valued function
$u(t,\omega,\cdot): \mR_+\times\Omega\to\cS'$ with that,
for any $\varphi\in C^{\infty}_{0}(\mathbb{R}^N)$,
$$
(t,\omega)\mapsto \<u(t,\omega,\cdot),\varphi\>
$$
is predictable.
For $p\geq 1$, we denote by $\mathbf{D}^p$ the subset of $\mathbf{D}$ consisting of all $u$ such that for any
$\varphi\in C^{\infty}_0(\mR^N)$ and $T,R\in\mathbb{R}_+$,
\begin{align}
\int^T_0\sup_{|x|\leq R}|\<u(t,\omega,\cdot),\varphi(\cdot-x)\>|^p\dif t<\infty\quad\quad(\mathrm{a.s.}-\omega).\label{iw1}
\end{align}
In the same way, we define $\mathbf{D}^p(\ell^2)$ by $u=(u_1,\cdots)$ with $u_i\in\mathbf D$
and replacing $|\cdot|$ in \eqref{iw1} by $\|\cdot\|_{\ell^2}$.

\begin{definition}
[Generalized It\^o's process]
Let $u, f\in\mathbf{D}$, $g\in\mathbf{D}(\ell^2)$ and $\tau$ be a stopping time. One says that the equality
\begin{eqnarray}
\dif u(t,x)=f(t,x)\dif t+g^k(t,x)\dif W^k_t,\quad t\leq \tau, \label{iw2}
\end{eqnarray}
holds in the distribution sense if $\1_{(0,\tau]}f\in\mathbf{D}^1$, $\1_{(0,\tau]}g\in\mathbf{D}^2(\ell^2)$
and for any $\varphi\in C^{\infty}_0(\mR^N)$, with probability one, for all $t\in\mathbb{R}_+$,
\begin{eqnarray}
\<u(t\wedge\tau), \varphi\>=\<u(0),\varphi\>+\int_0^{t\wedge\tau}\<f(s),\varphi\>\dif s+\int_0^{t\wedge\tau}\<g^k(s),\varphi\>\dif W^k_s.
\end{eqnarray}
\end{definition}

Let $X_t$ be an $\mathbb{R}^N-$valued stochastic process given by
$$
X^i_t=\int^t_0b^i_s\dif s+\int^t_0\sigma^{ik}_s\dif W^k_s,\ i=1,\cdots,N,
$$
where $b_t$, $\sigma^{\cdot k}_t$ are predictable $\mathbb{R}^N-$valued processes
with
$$
\int^T_0(|b_t|+\|\sigma_t\|_{\mR^N\otimes\ell^2}^2)\dif t<\infty,\ a.s., \ \forall T>0.
$$
The following generalized It\^o-Wentzell's formula are proven in
\cite{Krylov11}*{Theorem 1.1}.
\begin{lemma}[It\^o-Wentzell's formula]\label{IW}
Let $f,u\in\mathbf{D}$, $g\in\mathbf{D}(\ell^2)$  and $\tau$ be a stopping time,
$X$ an It\^o's process. Suppose \eqref{iw2} holds. Then $w(t,x)=u(t,x+X_t)$ satisfies
\begin{align}
\begin{split}
\dif w(t,x)&=\Big[\tfrac{1}{2}\sigma^{ik}_t\sigma^{jk}_t\p_i\p_jw(t,x)+b_t\cdot\nabla w(t,x)\Big]\dif t\nonumber\\
&\quad+\Big[f(t,x+X_t)+\p_ig^k_t(x+X_t)\sigma^{ik}_t\Big]\dif t\nonumber\\
&\quad+\Big[g^k_t(x+X_t)+\sigma^{ik}_t\p_iw(t,x)\Big]\dif W^k_t,\quad t\leq\tau,
\end{split}
\end{align}
in the distribution sense.
\end{lemma}

\section{Model equations}

Let $W_t$ be a $d$-dimensional standard Brownian motion. Define
$$
(X_t, V_t):=\left(\sqrt{2}\int^t_0 W_s\dif s, \sqrt{2}W_t\right)
$$
and the kinetic semigroup
\begin{align}\label{Ki}
P_t f(x,v):=\mE f(x+tv+X_t, v+V_t)=(\Gamma_tp_t)*(\Gamma_tf)(x,v),
\end{align}
where $p_t$ is the density of $(X_t, V_t)$ given by
$$
p_t(x,v)=\left(\frac{2\pi t^4}{3}\right)^{-\frac{d}{2}}\exp\left(-\frac{3|x|^2+|3x-2tv|^2}{4t^3}\right),
$$
and 
\begin{align}\label{Ga1}
\Gamma_t f(x,v):=f(x+tv,v).
\end{align}
It is easy to see that for any $\lambda>0$,
$$
p_{\lambda t}(x,v)=\lambda^{-2d}p_t(\lambda^{-3/2}x,\lambda^{-1/2}v)
$$
and for $f\in\cS$, by It\^o's formula,
\begin{align}\label{KE}
\p_tP_t f=[\Delta_v+v\cdot\nabla_x]P_tf,
\end{align}
which also holds for $f\in\cS'(\mB)$ in the distributional sense by duality.

In the remainder of this paper, we take the parameters in Subsection 2.1 as 
$$
N=2d,\ n=2,\ m_1=m_2=d,\ \theta=(3,1).
$$
For notational simplicity, we use $z=(x,v)$ to denote a generic point in $\mR^{2d}$, and 
 for a Banach space $\mB$, and for $T>0$ and $p\in[1,\infty]$, we write
$$
\mL^p_T(\mB):=L^p([0,T];\mB),\ \ \mL^p_{z,\omega}(\mB):=\mL^p_z(\mL^p_\omega(\mB))=\mL^p_\omega(\mL^p_z(\mB)).
$$
The operator $\Gamma_t$ in \eqref{Ga1} will play a crucial role below. Note that it is in general not a bounded operator in $\bB^{s}_{p;\theta}$, 
but we obviously have
$$
\|\Gamma_tf\|_{\mL^p_z}=\|f\|_{\mL^p_z}.
$$

\subsection{Estimates of kinetic semigroup $P_t$}
In this subsection we show some basic estimates about $P_t$ in anisotropic Besov spaces.
First of all we have the following estimate about the heat kernel $p_t$
that is similar to \cite{sen}*{Lemma 5.1}. 
\bl
\label{ies}
For any $l\geq 0$, there is a constant $C=C(d,l)>0$ such that
\begin{align}\label{C2}
\|\mathcal{R}^\theta_j\Gamma_{t}p_t\|_{\mL^1_z}\lesssim_C (t 4^{j})^{-l},\ \ \forall j\in\mN,\ t>0.
\end{align}
\el
\begin{proof}
Note that
\begin{align*}
\Gamma_t p_t(x,v)=t^{-2d}p_1(t^{-\frac32}x+ t^{-\frac12} v,t^{-\frac12}v).
\end{align*}
Let $\hbar:=t^{-\frac12}2^{-j}$. Denote the left hand side of \eqref{C2} by $\sI$.
By the change of variable, we have
\begin{align*}
\sI&=\int_{\mR^{2d}}\left|\int_{\mR^{2d}}
\check\phi^\theta_1(\bar x,\bar v)p_1(x-\hbar^{3}\bar x+v-\hbar\bar v,v-\hbar\bar v)\dif \bar x\dif\bar v\right|\dif x\dif v.
\end{align*}
We define the following operators: for $m\in\mN$,
$$
\Lambda^m f(x,v):=(\Delta_x+\Delta_v^3)^mf(x,v),
$$
and for $\hat f\in\cS$ with $0\notin$ support of $\hat f$,
$$
\widehat{\Lambda^{-m} f}(\xi,\eta):=(|\xi|^2+|\eta|^{6})^{-m}\hat f(\xi,\eta).
$$
Let
$$
H(\bar x,\bar v,x,v):=p_1(x-\hbar^{3}\bar x+v-\hbar\bar v,v-\hbar\bar v).
$$
Since $\phi^{\rm \theta}_1(\xi,\eta)\subset B^{\theta}_4\setminus B^{\theta}_1$, for $m\in\mN$, we have
\begin{align*}
\sI&=\int_{\mR^{2d}}\left|\int_{\mR^{2d}}
\Lambda^{-m}\check\phi^\theta_1(\bar x,\bar v)\Lambda^m H(\cdot,\cdot,x,v)(\bar x,\bar v)\dif \bar x\dif\bar v\right|\dif x\dif v\\
&\leq\int_{\mR^{2d}}|\Lambda^{-m}\check\phi^\theta_1(\bar x,\bar v)|\dif \bar x\dif\bar v
\int_{\mR^{2d}}\left|\Lambda^m H(\cdot,\cdot,x,v)\right|(\bar x,\bar v)\dif x\dif v.
\end{align*}
By the chain rule and elementary calculations,  we have
\begin{align*}
\sup_{\bar x,\bar v}\int_{\mR^{2d}}\left|\Lambda^mH(\cdot,\cdot,x,v)\right|(\bar x,\bar v)\dif x\dif v \lesssim {\hbar}^{6m}.
\end{align*}
Since $\Lambda^{-m}\check\phi^\theta_1$ is a Schwartz function, we thus have
$$
\sI\lesssim {\hbar}^{6m}\int_{\mR^{2d}}|\Lambda^{-m}\check\phi^\theta_1(\bar x,\bar v)|\dif \bar x\dif\bar v \lesssim {\hbar}^{6m}.
$$
Moreover, it is easy to see that
$$
\sI\leq \int_{\mR^{2d}}|\check\phi^\theta_1(\bar x,\bar v)|\dif \bar x\dif\bar v\lesssim 1.
$$
The desired estimate now follows by the above two estimates.
\end{proof}
By this lemma, we can show the following crucial estimates.
\bl
(i) For any $p\in[1,\infty]$, $\alpha\in\mR$, $l\geq 0$ and $T\geq 1$, there is a constant $C=C(T,d,\alpha,l)>0$ 
such that for all $t\in[0,T]$ and $j\in\mN_0$, 
\begin{align}\label{AZ1}
\|\cR^\theta_j P_{t}f\|_{\mL^p_z(\mB)}\lesssim_C2^{-j\alpha} h_{l,\alpha}(t 4^{j})\|f\|_{\bB^{\alpha}_{p;\theta}(\mB)},
\end{align}
where
$$
h_{l,\alpha}(t):=(1\wedge t^{-l})\big(1+t\big)^{|\alpha|}.
$$
(ii) Let $p\in[1,\infty]$, $\alpha\in\mR$, $r\in[1,\infty)$, $q\in[r,\infty]$ and  $\kappa,\kappa'\in(0,1]$ with $\frac{\kappa'}{r}>\frac{\kappa}q$
when $q>r$ and $\kappa'\geq\kappa$ when $q=r$.
For any $T\geq 1$, there is a constant
$C=C(T,d,\alpha,r,q,\kappa,\kappa')>0$ such that for all $j\in\mN_0$ and $t\in[0,T]$,
\begin{align}\label{C4}
\begin{split}
&\left(\int^t_0(t-s)^{\kappa'-1}\|\cR^\theta_j P_{t-s}f(s)\|^r_{\mL^p_z(\mB)}\dif s\right)^{1/r}\\
&\qquad\lesssim_C 2^{-(\alpha+2(\frac{\kappa'}r-\frac{\kappa} q))j}
\left(\int^t_0(t-s)^{\kappa-1}\|f(s)\|^q_{\bB^{\alpha}_{p;\theta}(\mB)}\dif s\right)^{1/q}.
\end{split}
\end{align}
\el
\begin{proof}
(i) By definition \eqref{Ki}, we have (see \cite[Lemm 6.7]{sen}),
$$
\cR^\theta_j P_{t}f=\cR^\theta_j\Gamma_t p_t*\Gamma_t f
=\sum_{\ell\in{\Theta}^t_j}(\cR^\theta_j\Gamma_t p_t)*(\Gamma_t\cR^\theta_\ell f),\ j\in\mN_0,
$$
where
$$
\Theta^t_0:=\Big\{\ell\in\mathbb{N}_0\,\big{|}\,2^\ell\leq 2^4(1+t)\Big\}
$$
and for $j\geq 1$,
$$
\Theta^t_j:=\Big\{\ell\in\mathbb{N}_0\,\big{|}\,2^\ell\leq2^4(2^j+t2^{3j}),\ 2^j\leq2^4(2^\ell+t2^{3\ell})\Big\}.$$
If $j\geq 1$, then by Young's inequality, we have for any $l\geq 0$ and $\alpha\not=0$,
\begin{align*}
\|\cR^\theta_j P_{t}f\|_{\mL^p_z(\mB)}
&\leq\|\cR^\theta_j\Gamma_t p_t\|_{\mL^1_z}\sum_{\ell\in{\Theta}^t_j}\|\Gamma_t\cR^\theta_\ell f\|_{\mL^p_z(\mB)}\no\\
&\stackrel{\eqref{C2}}{\lesssim} (t 4^{j})^{-l}\sum_{\ell\in{\Theta}^t_j}\|\cR^\theta_\ell f\|_{\mL^p_z(\mB)}
\lesssim (t 4^{j})^{-l}\sum_{\ell\in{\Theta}^t_j}2^{-\alpha\ell}\|f\|_{\bB^{\alpha}_{p;\theta}(\mB)}\no\\
&\lesssim (t 4^{j})^{-l} 2^{-j\alpha}\big(1+t 4^j\big)^{|\alpha|}\|f\|_{\bB^{\alpha}_{p;\theta}(\mB)},
\end{align*}
where the last step is due to \cite[Lemma 6.7]{sen}. By the arbitrariness of $l\geq 0$,
$$
\|\cR^\theta_j P_{t}f\|_{\mL^p_z(\mB)}\leq 2^{-j\alpha}(1\wedge (t 4^{j})^{-l}) \big(1+t 4^j\big)^{|\alpha|}\|f\|_{\bB^{\alpha}_{p;\theta}(\mB)}.
$$
If $j=0$, then we similarly have for $\alpha\not=0$,
\begin{align*}
\|\cR^\theta_0 P_{t}f\|_{\mL^p_z(\mB)}
&\leq\|\cR^\theta_0\Gamma_t p_t\|_{\mL^1_z}\sum_{\ell\in{\Theta}^t_0}\|\Gamma_t\cR^\theta_\ell f\|_{\mL^p_z(\mB)}
\lesssim \sum_{\ell\in{\Theta}^t_0}2^{-\alpha\ell}\|f\|_{\bB^{\alpha}_{p;\theta}(\mB)}\\
&\lesssim  \big(1+t\big)^{|\alpha|}\|f\|_{\bB^{\alpha}_{p;\theta}(\mB)}\leq T^l
(1\wedge t^{-l})\big(1+t\big)^{|\alpha|}\|f\|_{\bB^{\alpha}_{p;\theta}(\mB)}.
\end{align*}
Thus we obtain \eqref{AZ1} for $\alpha\not=0$. For $\alpha=0$, it follows by interpolation theorem.

(ii) Denote the left hand side of \eqref{C4} by $\sI$. For $q\in(r,\infty]$,
let $q'\in[r,\infty)$ with $\frac1{q'}+\frac 1 q=\frac1r$. By \eqref{AZ1}
and H\"older's inequality, we have
\begin{align*}
\sI&\lesssim 2^{-j\alpha}\left(\int^t_0(t-s)^{\kappa'-1}h^r_{l,\alpha}((t-s)4^j)\|f(s)\|^r_{\bB^{\alpha}_{p;\theta}(\mB)}\dif s\right)^{1/r}\\
&\lesssim 2^{-j\alpha}\left(\int^t_0h^{q'}_{l,\alpha}((t-s)4^j)(t-s)^{\frac{q'(1-\kappa)}{q}-\frac{q'(1-\kappa')}{r} }\dif s\right)^{1/q'}\\
&\quad\times\left(\int^t_0\|f(s)\|^q_{\bB^{\alpha}_{p;\theta}(\mB)}(t-s)^{\kappa-1}\dif s\right)^{1/q}.
\end{align*}
For the first integral denoted by $\sI_0$, by the change of variable we have
\begin{align*}
\sI_0&=4^{-(\frac1{q'}+\frac{1-\kappa}{q}-\frac{1-\kappa'}r)j}\left(\int^{t4^j}_0h^{q'}_{l,\alpha}(s)s^{\frac{q'(1-\kappa)}{q}-\frac{q'(1-\kappa')}{r}}\dif s\right)^{1/q'}\\
&\leq 4^{-(\frac{\kappa'}r-\frac\kappa q)j}\left(\int^\infty_0(1\wedge s^{-lq'})(1+s)^{|\alpha|q'}s^{\frac{q'(1-\kappa)}{q}-\frac{q'(1-\kappa')}{r}}\dif s\right)^{1/q'}
\end{align*}
In particular, for $l$ large enough, the last integral is finite. Thus,
$$
\sI\lesssim 2^{-j(\alpha+2(\frac{\kappa'}r-\frac\kappa q))}\left(\int^t_0\|f(s)\|^q_{\bB^{\alpha}_{p;\theta}(\mB)}(t-s)^{\kappa-1}\dif s\right)^{1/q}.
$$
For $q=r$, it is similar. The proof is complete.
\end{proof}

Next we show the following estimate about the kinetic semigroup $P_t$.
\bl
For any $\alpha\in\mR$, $\kappa\geq 0$ and $T>0$, there exists a constant $C=C(T, d,\alpha,\kappa)>0$ such that for all $t\in(0,T]$,
\begin{align}\label{HH7}
\|P_t f\|_{\bB^{\alpha+\kappa}_{p;\theta}(\mB)}\lesssim_C t^{-\frac\kappa 2}\|f\|_{\bB^\alpha_{p;\theta}(\mB)},\ \ \ 
\nor P_t f\nor_{\wt\bB^{\alpha+\kappa}_{p;\theta}(\mB)}\lesssim_C t^{-\frac\kappa 2}\nor f\nor_{\wt\bB^\alpha_{p;\theta}(\mB)}.
\end{align}
\el
\begin{proof}
Let $T>0$.
By  \eqref{AZ1}, we have for any $j\in\mN_0$ and $l\geq|\alpha|+\frac\kappa2$,
\begin{align}
&\|\cR^\theta_j P_tf\|_{\mL^p_z(\mB)}
\lesssim  2^{-\alpha j} h_{l,\alpha}(t 4^{j})\|f\|_{\bB^{\alpha}_{p;\theta}(\mB)}\no\\
&\qquad=2^{-(\alpha+\kappa)j}t^{-\frac\kappa 2}\|f\|_{\bB^{\alpha}_{p;\theta}(\mB)} (t4^j)^{\frac\kappa2} h_{l,\alpha}(t 4^{j})\no\\
&\qquad\leq 2^{-(\alpha+\kappa)j}t^{-\frac\kappa 2}\|f\|_{\bB^{\alpha}_{p;\theta}(\mB)}\sup_{s>0}s^{\frac\kappa2} h_{l,\alpha}(s),\label{ED8}
\end{align}
where $\sup_{s>0}s^{\frac\kappa2} h_{l,\alpha}(s)<\infty$ for $l\geq |\alpha|+\frac\kappa2$ and $\kappa\geq 0$. Thus by definition,
$$
\|P_t f\|_{\bB^{\alpha+\kappa}_{p;\theta}(\mB)}=\sup_{j\geq 0}2^{(\alpha+\kappa)j}\|\cR^\theta_j P_tf\|_{\mL^p_z(\mB)}\lesssim 
t^{-\frac\kappa 2}\|f\|_{\bB^{\alpha}_{p;\theta}(\mB)},
$$
which gives the first estimate in \eqref{HH7}. 

Next we prove the second estimate in \eqref{HH7}.
Let $\chi$ be as in \eqref{Cu1}. Fix $\delta>0$, $z_0:=(x_0,v_0)\in\mR^d\times\mR^d$, and for $t>0$, let
\begin{align}\label{CH6}
z_{0,t}:=(x_0-tv_0, v_0),\ \chi_{z_{0}}^\delta(t,z):=\chi((z-z_{0,t})/\delta)
\end{align}
and define
\begin{align}\label{CH60}
u(t,z):=P_tf(z),\ \ u^\delta_{z_0}(t,z):=\chi^\delta_{z_0}(t,z)u(t,z).
\end{align}
By definition and \eqref{KE}, it is easy to see that
$$
\p_t u^\delta_{z_0}=\Delta_v u^\delta_{z_0}+v\cdot\nabla_x u^\delta_{z_0}+F^\delta_{z_0},
$$
where 
$$
F^\delta_{z_0}:=-2\nabla_v u\cdot \nabla_v\chi^\delta_{z_0}-((v-v_0)\cdot\nabla_x\chi^\delta_{z_0}+\Delta_v\chi^\delta_{z_0})u.
$$
By Duhamel's formula, we have
\begin{align*}
u^\delta_{z_0}(t)=P_t(u^\delta_{z_0}(0))+\int^t_0P_{t-s}F^\delta_{z_0}(s)\dif s=:I_1(t)+I_2(t).
\end{align*}
For $I_1(t)$, since $u^\delta_{z_0}(0,z)=f(z)\chi^\delta_{z_0}(0,z)$, by \eqref{ED8}, we clearly have
$$
\|\cR^\theta_j I_1(t)\|_{\mL^p_z(\mB)}\lesssim 2^{-(\alpha+\kappa)j}t^{-\frac\kappa 2}\|u^\delta_{z_0}(0)\|_{\bB^{\alpha}_{p;\theta}(\mB)}
\lesssim 2^{-(\alpha+\kappa)j}t^{-\frac\kappa 2}\nor f\nor_{\wt\bB^{\alpha}_{p;\theta}(\mB)}.
$$
For $I_2(t)$, noting that by $\chi^{2\delta}_{z_0}=1$ on the support of $\chi^{\delta}_{z_0}$,
\begin{align*}
F^\delta_{z_0}=-2\nabla_v u^{2\delta}_{z_0}\cdot \nabla_v\chi^\delta_{z_0}-((v-v_0)\cdot\nabla_x\chi^\delta_{z_0}+\Delta_v\chi^\delta_{z_0})u^{2\delta}_{z_0},
\end{align*}
by \eqref{C4} with $q=r=1$ and \eqref{Ea1}, we have
\begin{align*}
\|\cR^\theta_j I_2(t)\|_{\mL^p_z(\mB)}
&\lesssim 2^{-(\alpha+\kappa)j}\int^t_0(t-s)^{-\frac\kappa2}\|F^\delta_{z_0}(s)\|_{\bB^{\alpha}_{p;\theta}(\mB)}\dif s\\
&\lesssim 2^{-(\alpha+\kappa)j}\int^t_0(t-s)^{-\frac\kappa2}\|u^{2\delta}_{z_0}(s)\|_{\bB^{\alpha+1}_{p;\theta}(\mB)}\dif s\\
&\lesssim 2^{-(\alpha+\kappa)j}\int^t_0(t-s)^{-\frac\kappa2}\nor u(s)\nor_{\wt\bB^{\alpha+1}_{p;\theta}(\mB)}\dif s.
\end{align*}
Thus, combining the above calculations, we obtain that for any $\kappa\in[0,2)$,
\begin{align}
&\nor u(t)\nor_{\wt\bB^{\alpha+\kappa}_{p;\theta}(\mB)}=\sup_{z_0}\sup_{j\geq 0} 2^{(\alpha+\kappa)j}\|\cR^\theta_j u^\delta_{z_0}(t)\|_{\mL^p_z(\mB)}\no\\
&\qquad\lesssim t^{-\frac\kappa 2}\nor f\nor_{\wt\bB^{\alpha}_{p;\theta}(\mB)}
+\int^t_0(t-s)^{-\frac\kappa2}\nor u(s)\nor_{\wt\bB^{\alpha+1}_{p;\theta}(\mB)}\dif s.\label{ED9}
\end{align}
In particular, letting $\kappa=1$, we get
\begin{align*}
\nor u(t)\nor_{\wt\bB^{\alpha+1}_{p;\theta}(\mB)}
\lesssim t^{-\frac12}\nor f\nor_{\wt\bB^{\alpha}_{p;\theta}(\mB)}
+\int^t_0(t-s)^{-\frac12}\nor u(s)\nor_{\wt\bB^{\alpha+1}_{p;\theta}(\mB)}\dif s,
\end{align*}
which implies by Gronwall's inequality of Volterra's type (see \cite{Zh10}) that
$$
\nor u(t)\nor_{\wt\bB^{\alpha+1}_{p;\theta}(\mB)}
\lesssim t^{-\frac12}\nor f\nor_{\wt\bB^{\alpha}_{p;\theta}(\mB)}.
$$
Substituting this into \eqref{ED9}, we obtain the second estimate in \eqref{HH7} for $\kappa\in(0,2)$.
For general $\kappa\geq 2$, it follows by the semigroup property of $P_t$.
\end{proof}

The following lemma will be used to show the regularity in time variable.
\bl
(i) For any $p\in[1,\infty]$, $\alpha\in(0,1)$ and $R>0$, there is a constant $C=C(d,p,\alpha,R)>0$ such that for all $t>0$,
\begin{align}\label{DX12}
\|\Gamma_t(\chi^R_0 f)-\chi^R_0 f\|_{\mL^p_z(\mB)}\lesssim_C t^{\frac{\alpha}3}\|\chi^R_0 f\|_{\bB^\alpha_{p;\theta}(\mB)},
\end{align}
and for any $\gamma\in(0,\frac\alpha 3]$ and $T>0$, there is a $C=C(T,d,p,\alpha,\gamma,R)>0$ such that
\begin{align}\label{DX11}
\|\Gamma_t(\chi^R_0 f)-\chi^R_0 f\|_{\bB^\gamma_{p;\theta}(\mB)}\lesssim_C t^{\frac{\alpha-\gamma}3}\|\chi^R_0 f\|_{\bB^\alpha_{p;\theta}(\mB)},\ \ t\in[0,T].
\end{align}
(ii) Let $p\in[1,\infty)$. For any $f\in\mL^p_z(\mB)$, it holds that
\begin{align}\label{DX31}
\lim_{t\downarrow 0}\|P_t f-f\|_{\mL^p_z(\mB)}=0,
\end{align}
and for any $\gamma<\alpha$ and $f\in\bB^{\alpha}_{p;\theta}(\mB)$,
\begin{align}\label{DX32}
\lim_{t\downarrow 0}\|P_t f-f\|_{\bB^\gamma_{p;\theta}(\mB)}=0.
\end{align}
(iii) For any $p\in[1,\infty]$ and $\alpha\in(0,1)$, 
there is a constant $C=C(d,\alpha,p)>0$ such that for all $t\geq 0$,
\begin{align}\label{DB0}
\|P_{t}f-\Gamma_t f\|_{\mL^p_{z}(\mB)}&\lesssim_C t^{\frac{\alpha}2}\|f\|_{\bB^{\alpha}_{p;\theta}(\mB)},
\end{align}
and for any $\gamma\in(0,\frac{\alpha}{3}]$, there is a $C=C(d,\alpha,\gamma,p)>0$ such that for all $t\geq 0$,
\begin{align}\label{DB1}
\|P_{t}f-\Gamma_t f\|_{\bB^{\gamma}_{p;\theta}(\mB)}&\lesssim_C t^{\frac{\alpha-\gamma}2}\|f\|_{\bB^{\alpha}_{p;\theta}(\mB)}.
\end{align}
\el
\begin{proof}
(i) Write $f_R:=\chi^R_0 f$. Let $\cR^x_j$ be defined as in \eqref{DD9}. Note that
\begin{align*}
\|\Gamma_t\cR^x_j f_R(x,v)-\cR^x_j f_R(x,v)\|_\mB\leq t|v|\int^1_0\|\nabla_x\cR^x_j f_R\|_\mB(x+stv,v)\dif s.
\end{align*}
Since $f_R$ has support in $B^\theta_{2R}$, by Bernstein's inequality \eqref{Ber}, we have
$$
\|\Gamma_t\cR^x_j f_R-\cR^x_j f_R\|_{\mL^p_z(\mB)}\leq 2R t\|\nabla_x\cR^x_j f_R\|_{\mL^p_z(\mB)}\lesssim t2^j\|\cR^x_j f_R\|_{\mL^p_z(\mB)},
$$
and also,
$$
\|\Gamma_t\cR^x_j f_R-\cR^x_j f_R\|_{\mL^p_z(\mB)}\leq 2\|\cR^x_j f_R\|_{\mL^p_z(\mB)}.
$$
Therefore, 
\begin{align*}
\|\Gamma_tf_R-f_R\|_{\mL^p_z(\mB)}
&\leq\sum_{j\geq 0}\|\Gamma_t\cR^x_j f_R-\cR^x_j f_R\|_{\mL^p_z(\mB)}
\lesssim\sum_{j\geq 0}((t2^j)\wedge 1)\|\cR^x_j f_R\|_{\mL^p_z(\mB)}\\
&\leq\sum_{j\geq 0}((t2^j)\wedge 1)2^{-\frac{\alpha}{3}j}\|f_R\|_{\bB^{\alpha/3}_{p;x}(\mB)}
\lesssim t^{\frac\alpha 3}\|f_R\|_{\bB^{\alpha/3}_{p;x}(\mB)},
\end{align*}
where the last step is due to Lemma \ref{Le27}.
For \eqref{DX11}, noting that by \eqref{FG1},
$$
\|\delta^{(1)}_h(\Gamma_tf_R)\|_{\mL^p_z(\mB)}=\|\delta^{(1)}_{\Gamma_t h}f_R\|_{\mL^p_z(\mB)}
\lesssim |\Gamma_t h|^\alpha_\theta\|f_R\|_{\bB^{\alpha}_{p;\theta}(\mB)},
$$
where for $h=(h_1,h_2)\in\mR^{2d}$, $\Gamma_t h=(h_1+th_2,h_2)$,
we have
\begin{align*}
\|\delta^{(1)}_h(\Gamma_tf_R-f_R)\|_{\mL^p_z(\mB)}
&\leq \|\delta^{(1)}_h\Gamma_tf_R\|_{\mL^p_z(\mB)}+\|\delta^{(1)}_hf_R\|_{\mL^p_z(\mB)}\\
&\lesssim(|\Gamma_t h|^\alpha_\theta+|h|^\alpha_\theta)\|f_R\|_{\bB^{\alpha}_{p;\theta}(\mB)}.
\end{align*}
Moreover, by \eqref{DX12} we also have
$$
\|\delta^{(1)}_h(\Gamma_tf_R-f_R)\|_{\mL^p_z(\mB)}\leq2\|\Gamma_tf_R-f_R\|_{\mL^p_z(\mB)}\lesssim
t^{\frac\alpha 3}\|f_R\|_{\bB^{\alpha/3}_{p;x}(\mB)}.
$$
Hence, 
\begin{align*}
\|\delta^{(1)}_h(\Gamma_tf_R-f_R)\|_{\mL^p_z(\mB)}
&\lesssim (t^{\frac\alpha 3}\wedge(|\Gamma_t h|^\alpha_\theta+|h|^\alpha_\theta))\|f_R\|_{\bB^{\alpha}_{p;\theta}(\mB)}\\
&\leq C_T t^{\frac{\alpha-\gamma} 3}|h|^\gamma_\theta\|f_R\|_{\bB^{\alpha}_{p;\theta}(\mB)}.
\end{align*}
Thus we obtain \eqref{DX11} by \eqref{FG1}.

(ii) Note that by \eqref{Ki},
\begin{align*}
\|P_t f-f\|_{\mL^p_z(\mB)}&\leq\|\Gamma_tp_t*\Gamma_t f-\Gamma_tf\|_{\mL^p_z(\mB)}+\|\Gamma_t f-f\|_{\mL^p_z(\mB)}\\
&=\|p_t* f-f\|_{\mL^p_z(\mB)}+\|\Gamma_t f-f\|_{\mL^p_z(\mB)},
\end{align*}
which in turn yields \eqref{DX31}.
Moreover, for $\gamma<\alpha$, we have
\begin{align*}
\|P_t f-f\|_{\bB^\gamma_{p;\theta}(\mB)}&=\sup_{j\geq 0}2^{\gamma j}\|\cR^\theta_j(P_t f-f)\|_{\mL^p_z(\mB)}
\leq\sum_{j\geq 0}2^{\gamma j}\|\cR^\theta_j(P_t f-f)\|_{\mL^p_z(\mB)}.
\end{align*}
Since $\|\cR^\theta_jP_t f\|_{\mL^p_z(\mB)}\lesssim 2^{-\alpha j}\|f\|_{\bB^\alpha_{p;\theta}(\mB)}$, by the dominated convergence theorem,
\begin{align*}
\lim_{t\downarrow 0}\|P_t f-f\|_{\bB^\gamma_{p;\theta}(\mB)}
&\leq\sum_{j\geq 0}2^{\gamma j}\lim_{t\downarrow 0}\|\cR^\theta_j(P_t f-f)\|_{\mL^p_z(\mB)}\stackrel{\eqref{DX31}}{=}0.
\end{align*}

(iii) Note that by \eqref{FG1},
\begin{align}
\|\Gamma_tf(\cdot+\bar x,\cdot+\bar v)-\Gamma_tf\|_{\mL^p_z(\mB)}
&=\|f(\cdot+\bar x+t\bar v,\cdot+\bar v)-f\|_{\mL^p_z(\mB)}\no\\
&\lesssim (|\bar x+t\bar v|^{\frac\alpha3}+|\bar v|^{\alpha})\|f\|_{\bB^\alpha_{p;\theta}(\mB)}.\label{DD6}
\end{align}
Thus, by definition  we have
\begin{align*}
\|P_{t}f-\Gamma_t f\|_{\mL^p_{z}(\mB)}&\leq \int_{\mR^{2d}}\Gamma_tp_t(\bar x,\bar v)
\|\Gamma_tf(\cdot-\bar x,\cdot-\bar v)-\Gamma_tf\|_{\mL^p_z(\mB)}\dif \bar x\dif\bar v\\
&\leq \int_{\mR^{2d}}\Gamma_tp_t(\bar x,\bar v)(|\bar x+t\bar v|^{\frac\alpha3}+|\bar v|^{\alpha})\|f\|_{\bB^\alpha_{p;\theta}(\mB)}\dif \bar x\dif \bar v\\
&=\int_{\mR^{2d}}p_t(\bar x,\bar v)(|\bar x|^{\frac\alpha3}+|\bar v|^{\alpha})\dif \bar x\dif \bar v\|f\|_{\bB^\alpha_{p;\theta}(\mB)}
\lesssim t^{\frac{\alpha}2}\|f\|_{\bB^{\alpha}_{p;\theta}(\mB)}.
\end{align*}
For \eqref{DB1},
set $\hbar:=t^{-\frac12}2^{-j}$. Since $\|\Gamma_tp_t\|_{\mL^1_z}=1$, we have for $j\geq 0$,
\begin{align*}
&\|\cR^\theta_j P_{t}f-\cR^\theta_j\Gamma_t f\|_{\mL^p_z(\mB)}=\|(\Gamma_t p_t)*(\cR^\theta_j\Gamma_t f)-\cR^\theta_j\Gamma_t f\|_{\mL^p_z(\mB)}\\
&\qquad\leq\int_{\mR^{2d}}\Gamma_t p_t(\bar x,\bar v)
\big\|\cR^\theta_j\Gamma_t f(\cdot-\bar x,\cdot-\bar v)-\cR^\theta_j\Gamma_t f\big\|_{\mL^p_z(\mB)}\dif \bar x\dif\bar v\\
&\qquad\leq\int_{\mR^{2d}}\Gamma_t p_t(\bar x,\bar v)
\Big(|\bar x|\|\nabla_x\cR^\theta_j\Gamma_t f\|_{\mL^p_z(\mB)}
+|\bar v|\|\nabla_v\cR^\theta_j\Gamma_tf\|_{\mL^p_z(\mB)}\Big)\dif \bar x\dif\bar v\\
&\qquad\lesssim (t^{3/2}2^{3j}+t^{1/2}2^{j})\|\cR^\theta_j\Gamma_t f\|_{\mL^p_z(\mB)}
=(\hbar^{-3}+\hbar^{-1}) \|\cR^\theta_j\Gamma_t f\|_{\mL^p_z(\mB)}.
\end{align*}
Moreover, we clearly have
$$
\|(\cR^\theta_j\Gamma_t p_t)*(\Gamma_t f)-\cR^\theta_j\Gamma_t f\|_{\mL^p_z(\mB)}\leq 2\|\cR^\theta_j\Gamma_t f\|_{\mL^p_z(\mB)}.
$$
Thus
\begin{align*}
\|\cR^\theta_j P_{t}f-\cR^\theta_j\Gamma_t f\|_{\mL^p_z(\mB)}
&\lesssim ((\hbar^{-3}+\hbar^{-1})\wedge 1)\|\cR^\theta_j\Gamma_t f\|_{\mL^p_z(\mB)}.
\end{align*}
On the other hand, for $j\geq 1$, noting that
\begin{align*}
\cR^\theta_j\Gamma_t f(x,v)=\int_{\mR^{2d}}\check\phi^\theta_j(\bar x,\bar v)
[\Gamma_t f(x-\bar x,v-\bar v)-\Gamma_t f(x,v)]\dif\bar x\dif\bar v,
\end{align*}
by \eqref{DD6}, we have
\begin{align*}
\|\cR^\theta_j\Gamma_t f\|_{\mL^p_z(\mB)}
&\leq \int_{\mR^{2d}}\check\phi^\theta_j(\bar x,\bar v)
\|\Gamma_t f(\cdot-\bar x,\cdot-\bar v)-\Gamma_t f\|_{\mL^p_z(\mB)}\dif\bar x\dif\bar v\\
&\lesssim \int_{\mR^{2d}}\check\phi^\theta_j(\bar x,\bar v)
(|\bar x+t\bar v|^{\alpha/3}+|\bar v|^{\alpha})\|f\|_{\bB^\alpha_{p;\theta}(\mB)}\dif\bar x\dif\bar v\\
&\lesssim ((t2^{-j})^{\frac\alpha 3}+2^{-\alpha j})\|f\|_{\bB^\alpha_{p;\theta}(\mB)}
=t^{\frac\alpha 2}(\hbar^{\frac\alpha 3}+\hbar^{\alpha})\|f\|_{\bB^\alpha_{p;\theta}(\mB)}.
\end{align*}
Therefore, for $j\geq 1$ and $\gamma\in(0,\tfrac\alpha3]$,
\begin{align*}
\|\cR^\theta_j P_{t}f-\cR^\theta_j\Gamma_t f\|_{\mL^p_z(\mB)}
&\lesssim t^{\frac\alpha 2} ((\hbar^{-3}+\hbar^{-1})\wedge 1)(\hbar^{\frac\alpha 3}+\hbar^{\alpha})\|f\|_{\bB^\alpha_{p;\theta}(\mB)}\\
&\leq 2t^{\frac\alpha 2}\hbar^{\gamma}\|f\|_{\bB^\alpha_{p;\theta}(\mB)}=2t^{\frac{\alpha-\gamma}2}2^{-\gamma j}\|f\|_{\bB^\alpha_{p;\theta}(\mB)},
\end{align*}
which together with
$$
\|\cR^\theta_0 P_{t}f-\cR^\theta_0\Gamma_t f\|_{\mL^p_z(\mB)}
\lesssim ((t^{\frac32}+t^{\frac12})\wedge 1)\|\cR^\theta_0\Gamma_t f\|_{\mL^p_z(\mB)}\lesssim t^{\frac{\alpha-\gamma}{2}}\|f\|_{\mL^p_z(\mB)}
$$
yields \eqref{DB1}.
\end{proof}

\subsection{Simple model equations}
We consider the following simple model equation:
\begin{align}\label{spde44}
\dif u=[\nu\Delta_v u+v\cdot\nabla_x u+f]\dif t+g^k\dif W^k_t,\ \ u(0)=u_0,
\end{align}
where $\nu>0$ and $f\in\mathbf{D}^1$, $g\in\mathbf{D}^2(\ell^2)$, and $u_0$ is an $\sF_0$-measurable $\cS'$-valued random variable.
\bd
An $\cS'$-valued predictable process $u$ is called a distribution solution of 
the above SKE if for any $\varphi\in\cS$ and $t\geq 0$,
\begin{align*}
\<u(t),\varphi\>&=\<u_0,\varphi\>+\int^t_0\<u(s),\nu\Delta_v\varphi-v\cdot\nabla_x\varphi\>\dif s\\
&\quad+\int^t_0\<f(s),\varphi\>\dif s+\int^t_0\<g^k(s),\varphi\>\dif W^k_s,\ \mP-a.s.
\end{align*}
\ed

We have the following regularity estimate for this model equation.
\bt\label{Th326}
For any $p\geq 2$, $q\in[2,\infty]$, $\beta\in\mR$, $\kappa\in(0,1]$, $\alpha\leq \beta+2(1-\frac{\kappa}q)$ and $T>0$, 
there is a constant $C=C(T, d,p,q,\kappa,\alpha,\beta,\nu)>0$ such that
for any distribution solution $u$ of SKE \eqref{spde44} and $t\in(0,T]$, 
\begin{align}\label{SZ1}
\begin{split}
&\|u(t)\|_{\bB^{\beta+2(1-\frac{\kappa}q)}_{p;\theta}(\mL^p_\omega)}
\lesssim_C
t^{\frac{\alpha-\beta}2-1+\frac{\kappa}{q}}\|u_0\|_{\bB^{\alpha}_{p;\theta}(\mL^p_\omega)}\\
&\qquad+\left(\int^t_0(t-s)^{\kappa-1}
\Big[\|f(s)\|^q_{\bB^{\beta}_{p;\theta}(\mL^p_\omega)}+\|g(s)\|^q_{\bB^{\beta+1}_{p;\theta}(\mL^p_\omega(\ell^2))}\Big]\dif s\right)^{1/q},
\end{split}
\end{align}
and also,
\begin{align}\label{SZ2}
\begin{split}
&\nor u(t)\nor_{\wt\bB^{\beta+2(1-\frac{\kappa}q)}_{p;\theta}(\mL^p_\omega)}
\lesssim_C
t^{\frac{\alpha-\beta}2-1+\frac{\kappa}{q}}\nor u_0\nor_{\wt\bB^{\alpha}_{p;\theta}(\mL^p_\omega)}\\
&\qquad+\left(\int^t_0(t-s)^{\kappa-1}
\Big[\nor f(s)\nor^q_{\wt\bB^{\beta}_{p;\theta}(\mL^p_\omega)}+\nor g(s)\nor^q_{\wt\bB^{\beta+1}_{p;\theta}(\mL^p_\omega(\ell^2))}\Big]\dif s\right)^{1/q}.
\end{split}
\end{align}
\et
\begin{proof}
Without loss of generality, we assume $\nu=1$ in equation \eqref{spde44}.
Let $P_t$ be defined by \eqref{Ki}
and $\bar u:=P_t u_0$. Since $\p_t \bar u=\Delta_v \bar u+v\cdot\nabla_x\bar u$, one sees that
$\tilde u=u-\bar u$ solves
$$
\dif \tilde u=[\Delta_v \tilde u+v\cdot\nabla_x \tilde u+f]\dif t+g^k\dif W^k_t,\ \ \tilde u(0)=0.
$$
Thus by \eqref{HH7}, we may assume $u_0\equiv 0$.

(i) By Duhamel's formula, we have in the distributional sense, 
\begin{align*}
u(t)=\int^t_0 P_{t-s}f(s)\dif s+\int^t_0 P_{t-s}g^k(s)\dif W^k_s=:I_1(t)+I_2(t).
\end{align*}
For $I_1(t)$, by \eqref{C4} with $(r,\alpha)=(1,\beta)$ and $\mB:=\mL^p_\omega$, we have
\begin{align*}
\|\cR^\theta_j I_1(t)\|_{\mL^p_{z,\omega}}
&\lesssim\int^t_0 \|\cR^\theta_j P_{t-s}f(s)\|_{\mL^p_z(\mL^p_\omega)}\dif s\\
&\lesssim 2^{-(\beta+2(1-\frac{\kappa}q))j}\left(\int^t_0(t-s)^{\kappa-1}\|f(s)\|^q_{\bB^{\beta}_{p;\theta}(\mL^p_\omega)}\dif s\right)^{1/q}.
\end{align*}
For $I_2(t)$, by BDG's inequality and Minkowskii's inequality, we have
\begin{align}
\|\cR^\theta_j I_2(t)\|_{\mL^p_{z,\omega}}^p&=\mE\left\|\int^t_0 \cR^\theta_j P_{t-s}g^k(s)\dif W^k_s\right\|^p_{\mL^p_z}\no\\
&\lesssim \mE\left(\int^t_0 \|\cR^\theta_j P_{t-s}g(s)\|_{\mL^p_z(\ell^2)}^2\dif s\right)^{p/2}\no\\
&\lesssim \left(\int^t_0 \|\cR^\theta_j P_{t-s}g(s)\|_{\mL^p_z(\mL^p_\omega(\ell^2))}^2\dif s\right)^{p/2}.\label{DC1}
\end{align}
By \eqref{C4} with $(r,\alpha)=(2,\beta+1)$ and $\mB:=\mL^p_\omega(\ell^2)$, we obtain
$$
\|\cR^\theta_j I_2(t)\|_{\mL^p_z(\mL^p_\omega)}\lesssim 2^{-(\beta+2(1-\frac{\kappa}q))j}
\left(\int^t_0(t-s)^{\kappa-1}\|g(s)\|^q_{\bB^{\beta+1}_{p;\theta}(\mL^p_\omega(\ell^2))}\dif s\right)^{1/q}.
$$
Combining the above calculations, we obtain \eqref{SZ1}.

(ii) Let $\chi^\delta_{z_0}$ be as in \eqref{CH6} and define
$$
u^\delta_{z_0}(t,\omega,z):=u(t,\omega,z)\chi^\delta_{z_0}(t,z).
$$
It is easy to see that
\begin{align}\label{spde404}
\dif u^\delta_{z_0}=[\Delta_v u^\delta_{z_0}+v\cdot\nabla_x u^\delta_{z_0}+F^\delta_{z_0}]\dif t+g^k\chi^\delta_{z_0}\dif W^k_t,
\end{align}
where 
$$
F^\delta_{z_0}:=f\chi^\delta_{z_0}-2\nabla_v u\cdot\nabla_v\chi^\delta_{z_0}-((v-v_0)\cdot\nabla_x\chi^\delta_{z_0}+\Delta_v\chi^\delta_{z_0})u.
$$
Note that by \eqref{Ea1}, 
\begin{align*}
\sup_{z_0}\|F^\delta_{z_0}(s)\|_{\bB^{\beta}_{p;\theta}(\mL^p_\omega)}
\lesssim\nor f(s)\nor_{\wt\bB^{\beta}_{p;\theta}(\mL^p_\omega)}+
\nor u(s)\nor_{\wt\bB^{1+\beta}_{p;\theta}(\mL^p_\omega)}.
\end{align*}
By \eqref{SZ1}, we have
\begin{align*}
&\nor u(t)\nor_{\wt\bB^{\beta+2(1-\frac{\kappa}q)}_{p;\theta}(\mL^p_\omega)}
=\sup_{z_0\in\mR^{2d}}\|u^\delta_{z_0}(t)\|_{\bB^{\beta+2(1-\frac{\kappa}q)}_{p;\theta}(\mL^p_\omega)}\\
&\lesssim\left(\int^t_0(t-s)^{\kappa-1}\sup_{z_0}
\Big[\|F^\delta_{z_0}(s)\|^q_{\bB^{\beta}_{p;\theta}(\mL^p_\omega)}
+\|g\chi^\delta_{z_0}(s)\|^q_{\bB^{\beta+1}_{p;\theta}(\mL^p_\omega(\ell^2))}\Big]\dif s\right)^{1/q}\\
&\lesssim\left(\int^t_0(t-s)^{\kappa-1}
\Big[\nor u(s)\nor^q_{\wt\bB^{1+\beta}_{p;\theta}(\mL^p_\omega)}+\nor f(s)\nor^q_{\wt\bB^{\beta}_{p;\theta}(\mL^p_\omega)}+\nor g(s)\nor^q_{\wt\bB^{\beta+1}_{p;\theta}(\mL^p_\omega(\ell^2))}\Big]\dif s\right)^{1/q}.
\end{align*}
In particular, taking $q=2$, $\kappa=1$ and then by Gronwall's inequality, we get
\begin{align*}
\nor u(t)\nor^2_{\wt\bB^{\beta+1}_{p;\theta}(\mL^p_\omega)}
\lesssim\int^t_0
\Big[\nor f(s)\nor^2_{\wt\bB^{\beta}_{p;\theta}(\mL^p_\omega)}+\nor g(s)\nor^2_{\wt\bB^{\beta+1}_{p;\theta}(\mL^p_\omega(\ell^2))}\Big]\dif s.
\end{align*}
Substituting this into the above inequality we obtain \eqref{SZ2}.
\end{proof}

The following stochastic convolution lemma provides the existence of a continuous version of the solution in time variable, 
and also will be used to show the boundedness of $\|u(\omega,\cdot)\|_{\mL^\infty_{T,z}}$ in the probability sense.
\bl \label{37}(Stochastic convolutions)
Let $p\geq 2$, $T\geq 1$ and  $\alpha,\beta\in\mR$ with $\alpha<\beta+1-\frac2p$.
For any $g\in{\mathbf D}(\ell^2)\cap \mL^p_T(\bB^\beta_{p;\theta}(\mL^p_\omega(\ell^2)))$, there  is a continuous version 
$\wt\mS_g$ of
$[0,T]\ni t\mapsto \mS_g(t):=\int^t_0P_{t-s}g^k(s)\dif W^k_s\in \bB^\alpha_{p;\theta}$ so that
for some $C=C(p,T,\alpha,\beta)>0$ and any stopping time $\tau$,
\begin{align}\label{323}
\mE\left(\sup_{t\in[0,T\wedge\tau]}\left\|\wt\mS_g(t)\right\|^p_{\bB^\alpha_{p;\theta}}\right)\lesssim_C
\int^T_0\|g(s\wedge\tau)\|^p_{\bB^\beta_{p;\theta}(\mL^p_\omega(\ell^2))}\dif s.
\end{align}
\el
\begin{proof}
We use the factorization method in \cite{Da-Kw-Za}. Noting that
$$
\int^t_s(t-r)^{\kappa-1}(r-s)^{-\kappa}\dif r=\frac{\pi}{\sin(\kappa\pi)},\ \ 0\leq s<t<\infty,\ \kappa\in(0,1),
$$
by stochastic Fubini's theorem (cf. \cite{Krylov11}), we have  for each $t\in(0,T)$,
$$
\mS_g(t)=\int^t_0P_{t-s}g^k(s)\dif W^k_s=\frac{\sin(\kappa \pi)}{\pi}\int^t_0(t-s)^{\kappa-1}P_{t-s}G(s)\dif s=:\wt\mS_g(t),\ \mP-a.s.
$$
where
$$
G(t):=\int^t_0(t-s)^{-\kappa}P_{t-s}g^k(s)\dif W^k_s.
$$
We show the above defined $\wt\mS_g$ has the desired property.
Since $\alpha<\beta+1-\frac2p$, one can choose $\kappa,\eps\in(0,1)$ small enough so that
\begin{align}\label{FU2}
\alpha=\beta+1-\tfrac2p-\tfrac{2\kappa}p-\eps.
\end{align}
By Minkowskii's inequality and \eqref{C4}, we have
\begin{align*}
\left\|\cR^\theta_j\wt\mS_g(t)\right\|_{\mL^p_z}
&\lesssim\int^t_0(t-s)^{\kappa-1}\|\cR^\theta_j P_{t-s}G(s)\|_{\mL^p_z}\dif s\\
&\lesssim 2^{-\alpha j} \left(\int^t_0\|G(s)\|^p_{\bB^{\alpha-2(\kappa-\frac1p)}_{p;\theta}}\dif s\right)^{1/p},
\end{align*}
which implies that for all $t\in[0,T]$,
$$
\left\|\wt\mS_g(t)\right\|^p_{\bB^{\alpha}_{p;\theta}}\lesssim\int^t_0\|G(s)\|^p_{\bB^{\alpha-2(\kappa-\frac1p)}_{p;\theta}}\dif s
\leq\int^T_0\|G(s)\|^p_{\bB^{\alpha-2(\kappa-\frac1p)}_{p;\theta}}\dif s.
$$
From this a priori estimate, as in \cite[Lemma 1]{Da-Kw-Za}, one sees that $t\mapsto \wt\mS_g(t)$ is continuous
as long as the last integral is finite a.s.. In particular, 
\begin{align}\label{FU3}
&\mE\left(\sup_{t\in[0,T]}\left\|\wt\mS_g(t)\right\|^p_{\bB^{\alpha}_{p;\theta}}\right)
\lesssim\int^T_0\mE\|G(s)\|^p_{\bB^{\alpha-2(\kappa-\frac1p)}_{p;\theta}}\dif s\no\\
&\quad\stackrel{\eqref{CT1}}{\lesssim}\int^T_0\|G(t)\|^p_{\bB^{\alpha-2(\kappa-\frac1p)+\eps}_{p;\theta}(\mL^p_\omega)}\dif t
\stackrel{\eqref{FU2}}{=}\int^T_0\|G(t)\|^p_{\bB^{\beta+1-2\kappa-\frac{2\kappa}p}_{p;\theta}(\mL^p_\omega)}\dif t.
\end{align}
On the other hand, by BDG's inequality and \eqref{C4} again, we have
\begin{align*}
\mE\|\cR^\theta_j G(t)\|^p_{\mL^p_z}
&\lesssim \mE\left(\int^t_0(t-s)^{-2\kappa}\|\cR^\theta_jP_{t-s}g(s)\|^2_{\mL^p_z(\ell^2)}\dif s\right)^{p/2}\\
&\lesssim \left(\int^t_0(t-s)^{-2\kappa}\|\cR^\theta_jP_{t-s}g(s)\|^2_{\mL^p_z(\mL^p_\omega(\ell^2))}\dif s\right)^{p/2}\\
&\lesssim 2^{-(\beta+2(\frac{1-2\kappa}2-\frac{\kappa}p))p j}\int^t_0(t-s)^{\kappa-1}\|g(s)\|^p_{\bB^\beta_{p;\theta}(\mL^p_\omega(\ell^2))}\dif s.
\end{align*}
Hence,
$$
\|G(t)\|^p_{\bB^{\beta+2(\frac{1-2\kappa}2-\frac{\kappa}p)}_{p;\theta}(\mL^p_\omega)}\lesssim 
\int^t_0(t-s)^{\kappa-1}\|g(s)\|^p_{\bB^\beta_{p;\theta}(\mL^p_\omega(\ell^2))}\dif s.
$$
Substituting this into \eqref{FU3}, we obtain 
\begin{align*}
\mE\left(\sup_{t\in[0,T]}\left\|\wt\mS_g(t)\right\|^p_{\bB^{\alpha}_{p;\theta}}\right)
&\lesssim\int^T_0\int^t_0(t-s)^{\kappa-1}\|g(s)\|^p_{\bB^\beta_{p;\theta}(\mL^p_\omega(\ell^2))}\dif s\dif t\\
&=\frac1\kappa\int^T_0(T-s)^{\kappa}\|g(s)\|^p_{\bB^\beta_{p;\theta}(\mL^p_\omega(\ell^2))}\dif s\\
&\leq\frac{T^\kappa}\kappa\int^T_0\|g(s)\|^p_{\bB^\beta_{p;\theta}(\mL^p_\omega(\ell^2))}\dif s.
\end{align*}
Finally, for any stopping time $\tau$, since $\tilde\mS_g(t\wedge\tau)=\tilde\mS_{g(\cdot\wedge\tau)}(t\wedge\tau)$ a.s., we have
\begin{align*}
\mE\left(\sup_{t\in[0,T\wedge\tau]}\left\|\wt\mS_g(t)\right\|^p_{\bB^\alpha_{p;\theta}}\right)
&\leq
\mE\left(\sup_{t\in[0,T]}\left\|\wt\mS_{g(\cdot\wedge\tau)}(t)\right\|^p_{\bB^\alpha_{p;\theta}}\right)\\
&\leq\int^T_0\|g(s\wedge\tau)\|^p_{\bB^\beta_{p;\theta}(\mL^p_\omega(\ell^2))}\dif s.
\end{align*}
The proof is complete.
\end{proof}
The following corollary is a direct consequence of Theorem \ref{Th326} and Lemma \ref{37}.
\bc\label{Cor38}
Let $p\geq 2$, $T>0$ and $\alpha\leq\beta+2-\frac2p$. Suppose that $u_0\in\bB^{\alpha}_{p;\theta}(\mL^p_\omega)$ and
$$
\sI^T_{f,g}:=\int^T_0\Big[\|f(s)\|^p_{\bB^{\beta}_{p;\theta}(\mL^p_\omega)}
+\|g(s)\|^p_{\bB^{\beta+1}_{p;\theta}(\mL^p_\omega(\ell^2))}\Big]\dif s<\infty.
$$
Then for any $\gamma<\beta+2-\frac2p$,
the unique weak solution $u$ in Theorem \ref{Th326} admits a continuous version $\tilde u(\omega,\cdot)\in C((0,T]; \bB^{\gamma}_{p;\theta})$.
Moreover, if $\gamma<\alpha$, then $\tilde u(\omega,\cdot)\in C([0,T]; \bB^{\gamma}_{p;\theta})$, and in this case,
there is a constant $C=C(p,T,d,\alpha,\gamma,\beta)>0$ such that
$$
\mE\left(\|\tilde u\|^p_{C([0,T]; \bB^{\gamma}_{p;\theta})}\right)\lesssim\|u_0\|^p_{\bB^{\alpha}_{p;\theta}(\mL^p_\omega)}+\sI^T_{f,g}.
$$
\ec

Although the above corollary asserts the continuity of $u$ in time variable, it does not say any H\"older regularity in the time variable.
The following lemma is useful.
\bl
For any $p\geq 2$, $\beta\in(0,1)$, $\gamma\in(0,\frac\beta3]$ and $T>0$, there is a constant $C=C(T, d,p,\beta,\gamma,\nu)>0$ such that
for any distribution solution $u$ of SKE \eqref{spde44} and $0<t_1<t_2\leq T$, 
\begin{align}\label{ES9}
\begin{split}
&\|u(t_2)-\Gamma_{t_2-t_1}u(t_1)\|_{\bB^{\gamma}_{p;\theta}(\mL^p_\omega)}\lesssim_C
(t_2-t_1)^{\frac{\beta-\gamma}2}\\
&\quad\times\Big(\|u(t_1)\|_{\bB^{\beta}_{p;\theta}(\mL^p_\omega)}+
\|f\|_{\mL^\infty_{t_1,t_2}(\bB^{\beta-2}_{p;\theta}(\mL^p_\omega))}+\|g\|_{\mL^\infty_{t_1,t_2}(\bB^{\beta-1}_{p;\theta}(\mL^p_\omega(\ell^2)))}\Big),
\end{split}
\end{align}
where $\mL^\infty_{t_1,t_2}(\mB):=L^\infty([t_1,t_2]; \mB)$.
\el
\begin{proof}
For $0<t_1<t_2\leq T$, noting that
\begin{align*}
u(t_2)=\int^{t_2}_{t_1} P_{t_2-s}f(s)\dif s+\int^{t_2}_{t_1} P_{t_2-s}g^k(s)\dif W^k_s+P_{t_2-t_1}u(t_1),
\end{align*}
we have
\begin{align*}
u(t_2)-\Gamma_{t_2-t_1}u(t_1)
&=\int^{t_2}_{t_1} P_{t_2-s}f(s)\dif s+\int^{t_2}_{t_1} P_{t_2-s}g^k(s)\dif W^k_s\\
&\quad+(P_{t_2-t_1}-\Gamma_{t_2-t_1})u(t_1)\\
&=:I_1+I_2+I_3.
\end{align*}
For $I_1$, by \eqref{AZ1} with $l=3+(\gamma-\beta)/2$, we have for all $j\in\mN_0$,
\begin{align*}
\|\cR^\theta_j I_1\|_{\mL^p_z(\mL^p_\omega)}
&\lesssim \int^{t_2}_{t_1} \|\cR^\theta_j P_{t_2-s}f(s)\|_{\mL^p_z(\mL^p_\omega)}\dif s\\
&\lesssim 2^{j(2-\beta)}\int^{t_2}_{t_1} h_{l,\beta-2}((t_2-s)4^j)\|f(s)\|_{\bB^{\beta-2}_{p;\theta}(\mL^p_\omega)}\dif s\\
&\lesssim 2^{-j\beta}\int^{(t_2-t_1)4^j}_{0} h_{l,\beta-2}(s)\dif s\|f\|_{\mL^\infty_{t_1,t_2}(\bB^{\beta-2}_{p;\theta}(\mL^p_\omega))}\\
&\lesssim 2^{-j\beta}\int^{(t_2-t_1)4^j}_{0} s^{\frac{\beta-\gamma}{2}-1}\dif s\|f\|_{\mL^\infty_{t_1,t_2}(\bB^{\beta-2}_{p;\theta}(\mL^p_\omega))}\\
&\lesssim 2^{-j\gamma}(t_2-t_1)^{\frac{\beta-\gamma}{2}}\|f\|_{\mL^\infty_{t_1,t_2}(\bB^{\beta-1}_{p;\theta}(\mL^p_\omega))},
\end{align*}
where we have used that for $l=3+(\gamma-\beta)/2$,
$$
h_{l,\beta-2}(s)=(1\wedge s^{-l})(1+s)^{2-\beta}\leq (1\wedge s^{-l})(1+s)^2\leq 4s^{\frac{\beta-\gamma}{2}-1}.
$$
For $I_2$, as in \eqref{DC1}, by  \eqref{AZ1} with $l=(3+\gamma-\beta)/2$, we have for all $j\in\mN_0$,
\begin{align*}
\|\cR^\theta_j I_2\|_{\mL^p_z(\mL^p_\omega)}
&\lesssim \left(\int^{t_2}_{t_1} \|\cR^\theta_j P_{t_2-s}g(s)\|_{\mL^p_z(\mL^p_\omega(\ell^2))}^2\dif s\right)^{1/2}\\
&\lesssim 2^{j(1-\beta)}\left(\int^{t_2}_{t_1} h^2_{l,\beta-1}((t_2-s)4^j)\|g(s)\|_{\bB^{\beta-1}_{p;\theta}(\mL^p_\omega(\ell^2))}^2\dif s\right)^{1/2}\\
&\lesssim 2^{-j\beta}\left(\int^{(t_2-t_1)4^j}_{0} h^2_{l,\beta-1}(s)\dif s\right)^{1/2}\|g\|_{\mL^\infty_{t_1,t_2}
(\bB^{\beta-1}_{p;\theta}(\mL^p_\omega(\ell^2)))}\\
&\lesssim 2^{-j\beta}\left(\int^{(t_2-t_1)4^j}_{0} s^{\beta-\gamma-1}\dif s\right)^{1/2}\|g\|_{\mL^\infty_{t_1,t_2}(\bB^{\beta-1}_{p;\theta}(\mL^p_\omega(\ell^2)))}\\
&\lesssim 2^{-j\gamma}(t_2-t_1)^{\frac{\beta-\gamma}{2}}\|g\|_{\mL^\infty_{t_1,t_2}(\bB^{\beta-1}_{p;\theta}(\mL^p_\omega(\ell^2)))}.
\end{align*}
For $I_3$, since $\gamma\in(0,\frac\beta3]$, by  \eqref{DB1} we have
$$
\|I_3\|_{\bB^{\gamma}_{p;\theta}(\mL^p_\omega)}\lesssim (t_2-t_1)^{\frac{\beta-\gamma}{2}}\|u(t_1)\|_{\bB^{\beta}_{p;\theta}(\mL^p_\omega)}.
$$
Combining the above estimates, we obtain \eqref{ES9}.
\end{proof}

\br\rm
In the kinetic case, it is not possible to show the H\"older continuity of $t\mapsto u(t)$
in $\bB^{\gamma}_{p;\theta}(\mL^p_\omega)$ due to the transport term of $v\cdot\nabla_x$. 
In other words, one can not show the same
estimates   for $\|P_tf-f\|_{\bB^{\gamma}_{p;\theta}}$ as in \eqref{DB0} and \eqref{DB1}. 
However, \eqref{ES9} together with \eqref{DX11} still says some H\"older continuity of 
$t\mapsto u(t,x,v)$ in time variable locally in $x,v$ (see Corollary \ref{76} below).
\er

\section{SKEs with constant random coefficients}

In this section we consider the following SKE with  constant random coefficients:
\begin{align}\label{spde3}
\dif u=[\tr(a\cdot\nabla^2_vu)+v\cdot\nabla_x u+f]\dif t+[\sigma^{k}\cdot\nabla_v u+g^k]\dif W^k_t,\ u(0)=u_0,
\end{align}
where $u_0:\Omega\to\cS'$ is $\sF_0$-measurable, $f\in{\bf D}^1$, $g\in{\bf D}^2(\ell^2)$ and 
$$
(a,\sigma):\mR_+\times\Omega\to(\mR^d\otimes\mR^d,\mR^d\otimes\ell^2)
$$ 
are predictable processes.
Suppose that for some $c_0,c_1>0$ and all $(t,\omega)\in\mR_+\times\Omega$,
\begin{align}\label{Con1}
2a(t,\omega)-\sigma^{\cdot k}\sigma^{\cdot k}(t,\omega)\geq c_0\mI
\end{align}
and
\begin{align}\label{Con2}
\|a(t,\omega)\|_{\mR^d\otimes\mR^d}+\|\sigma(t,\omega)\|_{\mR^d\otimes\ell^2}\leq c_1.
\end{align}
\bd
We call a predictable process $u(t,\omega):\mR_+\times\Omega\to\cS'$ a distribution solution of SKE \eqref{spde3}
if for any $\varphi\in C^\infty_0(\mR^{2d})$ and $t\geq 0$,
\begin{align*}
\<u(t),\varphi\>&=\<u_0,\varphi\>+
\int^t_0[\<u, \tr(a\cdot\nabla^2_v\varphi)\>-\<u, v\cdot\nabla_x \varphi\>+\<f,\varphi\>]\dif s\\
&\quad+\int^t_0[\<g^k,\varphi\>-\<u,\sigma^{k}\cdot\nabla_v\varphi\>]\dif W^k_s.
\end{align*}
\ed

We use Krylov's trick \cite{kry99} 
to give a representation for the solution of SKE \eqref{spde3} in terms of the solution of the model equation \eqref{spde44}.
Let
\begin{align}\label{Def0}
\tilde\sigma:=(2a-\sigma\sigma^*- c_0\mI)^{1/2}\in\mR^d\otimes\mR^d
\end{align}
and
\begin{align}\label{Def1}
\tilde f:=f-\sigma^{k}\cdot\nabla_v g^k.
\end{align}
Let $(B^k)_{k=1,\cdots,d}$ be another Brownian motion that is independent of $W$.
Define
$$
\cI^i_t:=\int^t_0\sigma^{ik}(s)\dif W^k_s,\quad \tilde\cI^i_t:=\int^t_0\tilde{\sigma}^{ik}(s)\dif B^k_s, \ \ i=1,\cdots,d.
$$
We introduce the following auxiliary equation
\begin{align*}
\dif \bar u=\Big[\tfrac{c_0}{2}\Delta_v\bar u+v\cdot\nabla_x\bar u+\bar f\Big]\dif t+\bar g^k\dif W^k_t,\ \bar u(0)=u_0,
\end{align*}
where
$$
\bar f(t,x,v):=\tilde f\left(t,x+\int^t_0(\cI_s+\tilde\cI_s)\dif s,v-\cI_t-\tilde\cI_t\right)
$$
and
$$
\bar g(t,x,v):=g\left(t,x+\int^t_0(\cI_s+\tilde\cI_s)\dif s,v-\cI_t-\tilde\cI_t\right).
$$
Let
$$
w(t,x,v)=\bar u\left(t,x-\int^t_0(\cI_s+\tilde\cI_s)\dif s,v+\cI_t+\tilde\cI_t\right).
$$
By It\^o-Wentzell's formula (see Lemma \ref{IW}), we have
\begin{align*}
\dif w&=\Big[\tr (\tilde a\cdot\nabla^2 w)+v\cdot\nabla_x w+\tilde f+\sigma^{k}\cdot\nabla_v g^k\Big]\dif t\\
&\quad+(g^k+\sigma^{k}\cdot\nabla_{v} w)\dif W^k_t+(\tilde{\sigma}^{k}\cdot\nabla_{v} w)\dif B^k_t,
\end{align*}
where
$$
\tilde a:=\tfrac{1}{2}(\tilde{\sigma}\tilde{\sigma}^*+{\sigma}{\sigma}^*)+c_0\mI/2.
$$
Now we define
$$
u(t,x,v):=\mE(w(t,x,v)|\sF^W_t),
$$
where $\sF^W_t$ is the $\sigma$-algebra generated by $\{W_s, s\leq t\}$ and $\sF_0$.
By \eqref{Def0} and \eqref{Def1}, one sees that
$$
\tilde a=a,\ \ \tilde f+\sigma^{k}\cdot\nabla_v g^k=f.
$$
Note that (cf. \cite[p.28 Theorem 1.15]{roz})
$$
\mE\left(\int^t_0 \xi^k(s)\dif W^k_s\Big|\sF^W_t\right)=\int^t_0 \mE\left(\xi^k(s)|\sF^W_s\right)\dif W^k_s
$$
and
$$
\mE\left(\int^t_0 \eta^k(s)\dif B^k_s\Big|\sF^W_t\right)=0.
$$
In particular, one sees that $u$ solves SKE \eqref{spde3}. Now by Theorem \ref{Th326} and Fubini's theorem, we have
\bt\label{Th31}
Let $p\geq 2$, $q\in[2,\infty]$, $\beta\in\mR$, $\kappa\in(0,1]$, $\alpha\leq \beta+2(1-\frac{\kappa}q)$ and $T>0$.
Under \eqref{Con1} and \eqref{Con2}, there is a constant $C=C(T, d,p,q,\kappa,\alpha,\beta,c_0,c_1)>0$ so that 
 for any distribution solution $u$ of SKE \eqref{spde3} and $t\in(0,T]$,
 \begin{align}\label{SZ11}
\begin{split}
&\|u(t)\|_{\bB^{\beta+2(1-\frac{\kappa}q)}_{p;\theta}(\mL^p_\omega)}
\lesssim_C
t^{\frac{\alpha-\beta}2-1+\frac{\kappa}{q}}\|u_0\|_{\bB^{\alpha}_{p;\theta}(\mL^p_\omega)}\\
&\qquad+\left(\int^t_0(t-s)^{\kappa-1}
\Big[\|f(s)\|^q_{\bB^{\beta}_{p;\theta}(\mL^p_\omega)}+\|g(s)\|^q_{\bB^{\beta+1}_{p;\theta}(\mL^p_\omega(\ell^2))}\Big]\dif s\right)^{1/q}.
\end{split}
\end{align}
\et
\begin{proof}
Note that by Fubini's theorem,
\begin{align*}
\|\cR^\theta_j u(t)\|_{\mL^p_z(\mL^p_\omega)}&=\|\cR^\theta_j u(t)\|_{\mL^p_\omega(\mL^p_z)}
\leq\|\mE(\|\cR^\theta_j w(t)\|_{\mL^p_z}|\sF^W_t)\|_{\mL^p_\omega}\\
&\leq\|\cR^\theta_j w(t)\|_{\mL^p_\omega(\mL^p_z)}=\|\cR^\theta_j\bar u(t)\|_{\mL^p_\omega(\mL^p_z)}
=\|\cR^\theta_j\bar u(t)\|_{\mL^p_z(\mL^p_\omega)}.
\end{align*}
Therefore, by Theorem \ref{Th326},
\begin{align*}
&\|u(t)\|_{\bB^{\beta+2(1-\frac{\kappa}q)}_{p; \theta}(\mL^p_\omega)}
=\sup_j \Big(2^{(\beta+2(1-\frac{\kappa}q))j}\|\cR^\theta_j u(t)\|_{\mL^p_z(\mL^p_\omega)}\Big)\\
&\qquad\leq \|\bar u(t)\|_{\bB^{\beta+2(1-\frac{\kappa}q)}_{p; \theta}(\mL^p_\omega)}
\lesssim_C t^{\frac{\alpha-\beta}2-1+\frac{\kappa}{q}}\|u_0\|_{\bB^{\alpha}_{p;\theta}(\mL^p_\omega)}\\
&\qquad\quad+\left(\int^t_0(t-s)^{\kappa-1}\Big[\|\bar f(s)\|^q_{\bB^{\beta}_{p;\theta}(\mL^p_\omega)}
+\|\bar g(s)\|^q_{\bB^{\beta+1}_{p;\theta}(\mL^p_\omega(\ell^2))}\Big]\dif s\right)^{1/q}.
\end{align*}
Finally, it is just noticed that by Fubini's theorem again,
\begin{align*}
\|\bar f(s)\|_{\bB^{\beta}_{p;\theta}(\mL^p_\omega)}
=\|\tilde f(s)\|_{\bB^{\beta}_{p;\theta}(\mL^p_\omega)}
&\lesssim\|f(s)\|_{\bB^{\beta}_{p;\theta}(\mL^p_\omega)}
+\|\nabla_v g\|_{\bB^{\beta}_{p;\theta}(\mL^p_\omega(\ell^2))}
\end{align*}
and
\begin{align*}
\|\bar g(s)\|_{\bB^{1+\beta}_{p;\theta}(\mL^p_\omega(\ell^2))}
=\|g(s)\|_{\bB^{1+\beta}_{p;\theta}(\mL^p_\omega(\ell^2))}.
\end{align*}
Thus we obtain \eqref{SZ11}. 
\end{proof}

\section{SKEs with variable random coefficients}

In this section we consider SKEs \eqref{spde2} and \eqref{spde22} with variable coefficients,
and assume the following super-parabolic conditions:
\begin{enumerate}[{\bf (H$_0$)}]
\item There are $c_0,c_1>0$ such that for all $(t,\omega,x,v)$,
$$
2a(t,\omega,x,v)-\sigma\sigma^*(t,\omega,x,v)\geq c_0\mI
$$
and
$$
\|a(t,\omega,x,v)\|_{\mR^d}+\|\sigma(t,\omega,x,v)\|_{\mR^d\otimes\ell^2}\leq c_1.
$$
\end{enumerate}

\subsection{Optimal Besov regularity estimates}
In this subsection we consider SKE \eqref{spde2} and assume {\bf (H$_0$)} and
\begin{enumerate}[{\bf (H$^\beta_1$)}]
\item We suppose that for some $\beta\in(0,1)$ and $c_2>0$,
$$
\|a\|_{\mL^\infty_{T,\omega}(\bC^{\beta}_{\theta})}+\|b\|_{\mL^\infty_{T,\omega}(\bC^{\beta}_{\theta})}
+\|\sigma\|_{\mL^\infty_{T,\omega}(\bC^{1+\beta}_{\theta}(\ell^2))}\leq c_2,\ \forall T>0,
$$
where $\mL^\infty_{T,\omega}(\mB):=L^\infty([0,T]\times\Omega,\cP,\dif t\times\mP;\mB)$.
\end{enumerate}
We first introduce the following notion about the solution of SKE \eqref{spde2}.
\bd\label{Def41}
Let $T>0$, $\beta\in(0,1)$ and $p\in[2,\infty)$.
A predictable process $u\in \mL^\infty_T(\wt\bB^{2+\beta}_{p;\theta}(\mL^p_\omega))$ is called a weak solution of 
SKE \eqref{spde2} with initial value $u_0\in \wt\bB^{2+\beta}_{p;\theta}(\mL^p_\omega)$
if for any $\varphi\in C^\infty_c(\mR^{2d})$,
\begin{align*}
\<u(t),\varphi\>&=\<u_0,\varphi\>+\int^t_0[\<\tr(a\cdot\nabla^2_vu),\varphi\>-\<u,v\cdot\nabla_x\varphi\>+\<b\cdot\nabla_v u,\varphi\>]\dif s\\
&\quad+\int^t_0\<f,\varphi\>\dif s+\int^t_0\<\sigma^{k}\cdot\nabla_v u+g^k,\varphi\>\dif W^k_s\ \ a.s., \ \forall t\in[0,T].
\end{align*}
\ed
We have the following existence and uniqueness of weak solutions for SKE \eqref{spde2}.
\bt\label{Th41}
Let $T>0$, $p\in[2,\infty)$, $\beta\in(0,1)$  and
$$
u_0\in\widetilde{\bB}^{2+\beta}_{p;\theta}(\mL^p_\omega),\ \ f\in\mL^\infty_T(\widetilde{\bB}^{\beta}_{p;\theta}(\mL^p_\omega)),\ \ g\in\mL^\infty_T(\widetilde{\bB}^{1+\beta}_{p;\theta}(\mL^p_\omega(\ell^2))).
$$
Under {\bf (H$_0$)} and {\bf (H$^\beta_1$)}, there is a unique weak solution $u$ to SKE \eqref{spde2} so that
\begin{align}\label{ES1}
\nor u\nor_{\mL^\infty_T(\widetilde{\bB}^{2+\beta}_{p;\theta}(\mL^p_\omega))}\lesssim_C
\nor u_0\nor_{\widetilde{\bB} ^{2+\beta}_{p;\theta}(\mL^p_\omega)}+
\nor f\nor_{\mL^\infty_T(\widetilde{\bB} ^{\beta}_{p;\theta}(\mL^p_\omega))}+\nor g\nor_{\mL^\infty_T(\widetilde{\bB}^{1+\beta}_{p;\theta}(\mL^p_\omega(\ell^2)))},
\end{align}
where the constant $C=C(d,p,\beta,\kappa,T,c_0,c_1,c_2)>0$.
\et

\begin{proof}
Let $\bar u:=P_t u_0$ solve $\p_t \bar u=\Delta_v \bar u+v\cdot\nabla_x\bar u$. By \eqref{HH7} we have
$$
\nor \bar u\nor_{\mL^\infty_T(\widetilde{\bB}^{2+\beta}_{p;\theta}(\mL^p_\omega))}\lesssim_C
\nor u_0\nor_{\widetilde{\bB} ^{2+\beta}_{p;\theta}(\mL^p_\omega)}.
$$
If we consider $\tilde u=u-\bar u$, then $\tilde u$ solves SKE \eqref{spde2} with 
$$
\tilde f=f+\sL_v\bar u-\Delta_v \bar u,\quad \tilde g^k=g^k+\sigma^k\cdot\nabla_v\bar u,
$$
and $\tilde f$, $\tilde g$ satisfy the same assumptions as $f,g$. Thus,
without loss of generality we may assume $u_0\equiv0$.
We divide the proof into three steps.
In the first two steps, we prove the a priori estimate \eqref{ES1} by freezing coefficient argument. In the last step, we use
the continuity method to show the existence of weak solutions.

(i) Let $\chi^\delta_{z_0}$ be as in \eqref{CH6} and define
$$
a_{z_0}(t,\omega):=a(t,\omega,z_{0,t}),\ \ \sigma_{z_0}(t,\omega):=\sigma(t,\omega,z_{0,t})
$$
and
$$
u^\delta_{z_0}(t,\omega,z):=u(t,\omega,z)\chi^\delta_{z_0}(t,z).
$$
Let $u$ be any weak solution of SKE \eqref{spde2}. By Definition \ref{Def41}, it is easy to see that
$u^\delta_{z_0}$ is a weak solution of the following freezing SKE:
$$
\dif u^\delta_{z_0}=\Big[\tr (a_{z_0}\cdot\nabla^2_v u^\delta_{z_0})+v\cdot\nabla_x u^\delta_{z_0}
+\tilde f\Big]\dif t+\Big[\sigma^{k}_{z_0}\cdot\nabla_{v} u^\delta_{z_0}+\tilde g^k\Big]\dif W^k_t
$$
with $u^\delta_{z_0}(0)=0$, where
\begin{align*}
\tilde f&:=\tr ((a-a_{z_0})\cdot\nabla^2_vu^\delta_{z_0})
-2 \tr(a\cdot(\nabla_v u\otimes\nabla_v\chi^\delta_{z_0}))+
 \chi^\delta_{z_0} b\cdot\nabla_vu\\
&\quad-\big[\tr (a\cdot\nabla^2_v\chi^\delta_{z_0}) +(v-v_0)\cdot\nabla_x\chi^\delta_{z_0}\big]u+f\chi^\delta_{z_0},
\end{align*}
and
\begin{align*}
\tilde g&:=(\sigma-\sigma_{z_0})\cdot\nabla_{v} u^\delta_{z_0}-\sigma\cdot\nabla_v\chi^\delta_{z_0} u+g\chi^\delta_{z_0}.
\end{align*}
For simplicity of notations, if we let
$$
\bar a:=(a-a_{z_0})\chi^{2\delta}_{z_0},\
\bar b:=2a\cdot\nabla_v\chi^\delta_{z_0}-\chi^\delta_{z_0} b,\ \ \bar\sigma:=(\sigma-\sigma_{z_0})\chi^{2\delta}_{z_0},
$$
and
$$
\bar c:=-\big(\tr (a\cdot\nabla^2_v\chi^\delta_{z_0}) +(v-v_0)\cdot\nabla_x\chi^\delta_{z_0}\big),\
\bar h:=\sigma\cdot\nabla_v\chi^\delta_{z_0},
$$
then due to $\chi^{2\delta}_{z_0}=1$ on the support of $\chi^\delta_{z_0}$, we can write
\begin{align}\label{FF0}
\tilde f=\tr (\bar a\cdot\nabla^2_vu^\delta_{z_0})-\bar b\cdot\nabla_v u^{2\delta}_{z_0}
+\bar c\cdot u^{2\delta}_{z_0}
+f\chi^\delta_{z_0}
\end{align}
and
\begin{align}\label{FF1}
\tilde g=\bar\sigma\cdot\nabla_{v} u^\delta_{z_0}-\bar h\,u^{2\delta}_{z_0}+g\chi^\delta_{z_0}.
\end{align}

(ii) Let $T>0$, $p\in[2,\infty)$ and $\beta\in(0,1)$. For any $q\in[2,\infty]$, 
by Theorem \ref{Th31} with $\kappa=1$, there is a constant $C>0$
such that for all $\delta\in(0,1)$, $z_0\in\mR^{2d}$ and $t\in(0,T]$,
\begin{align}\label{DU70}
\begin{split}
&\|u^\delta_{z_0}(t)\|_{\bB^{\beta+2(1-\frac1 q)}_{p; \theta}(\mL^p_\omega)}
\lesssim_C
\left(\int^t_0
\Big[\|\tilde f(s)\|^q_{\bB^{\beta}_{p;\theta}(\mL^p_\omega)}+\|\tilde g(s)\|^q_{\bB^{\beta+1}_{p;\theta}(\mL^p_\omega(\ell^2))}\Big]\dif s\right)^{1/q}.
\end{split}
\end{align}
Now let us estimate $\|\tilde f(s)\|_{\bB^{\beta}_{p;\theta}(\mL^p_\omega)}$ and 
$\|\tilde g(s)\|_{\bB^{\beta+1}_{p;\theta}(\mL^p_\omega(\ell^2))}$. Below we drop the variable $s$.
For each $(t,z)\in\mR_+\times\mR^{2d}$, we introduce a bounded linear operator:
$$
\cT^{\bar a}_{t,z}:\mL^p_\omega\to\mL^p_\omega,\quad \cT^{\bar a}_{t,z}f(\omega):=\bar a(t,z,\omega)f(\omega).
$$
Noting that by the assumptions,
$$
\|\cT^{\bar a}_{t,z}f\|_{\mL^p_\omega}\leq c_2(2\delta)^\beta\|f\|_{\mL^p_\omega},\quad
\|(\cT^{\bar a}_{t,z}-\cT^{\bar a}_{t,z'})f\|_{\mL^p_\omega}\leq c_2|z-z'|^\beta\|f\|_{\mL^p_\omega},
$$
by \eqref{24}, \eqref{DQ9} and Young's inequality, we have
\begin{align*}
\|\tr (\bar a\cdot\nabla^2_v u^\delta_{z_0})\|_{\bB^{\beta}_{p;\theta}(\mL^p_\omega)}
&\lesssim  \delta^\beta\|\nabla^2_v u^\delta_{z_0}\|_{\bB^{\beta}_{p;\theta}(\mL^p_\omega)}
+\|\nabla^2_v u^\delta_{z_0}\|_{\mL^p_z(\mL^p_\omega)}\\
&\lesssim \delta^\beta\|u^\delta_{z_0}\|_{\bB^{2+\beta}_{p;\theta}(\mL^p_\omega)}+C_\delta\|u^\delta_{z_0}\|_{\mL^p_z(\mL^p_\omega)}.
\end{align*}
Similarly, since $\|\bar b\|_{\mL^\infty_{T,\omega}(\bC^\beta_\theta)}+\|\bar c\|_{\mL^\infty_{T,\omega}(\bC^\beta_\theta)}\lesssim 1$,
by \eqref{24} again,  we have
\begin{align*}
\|\bar b\cdot\nabla_v u^{2\delta}_{z_0}\|_{\bB^{\beta}_{p; \theta}(\mL^p_\omega)}
+\big\|\bar c\cdot u^{2\delta}_{z_0}\big\|_{\bB^{\beta}_{p;\theta}(\mL^p_\omega)}
\lesssim \|u^{2\delta}_{z_0}\|_{\bB^{1+\beta}_{p;\theta}(\mL^p_\omega)}\lesssim \|u^{2\delta}_{z_0}\|_{\bB^2_{p;\theta}(\mL^p_\omega)}.
\end{align*}
Therefore, for any $\delta\in(0,1)$ and some $C_0,C_1>0$ independent of $\delta$ and $z_0$,
\begin{align}\label{DU1}
\|\tilde f\|_{\bB^{\beta}_{p;\theta}(\mL^p_\omega)}
&\leq C_0\delta^\beta\|u^{\delta}_{z_0}\|_{\bB^{2+\beta}_{p;\theta}(\mL^p_\omega)}
+C_\delta\|u^{2\delta}_{z_0}\|_{\bB^2_{p;\theta}(\mL^p_\omega)}
+C_1\|f^{\delta}_{z_0}\|_{\bB^{\beta}_{p;\theta}(\mL^p_\omega)}\no\\
&\leq C_0\delta^\beta\nor u\nor_{\widetilde{\bB} ^{2+\beta}_{p;\theta}(\mL^p_\omega)}
+C_\delta\nor u\nor_{\widetilde{\bB}^2_{p;\theta}(\mL^p_\omega)}+C_1\nor f\nor_{\widetilde{\bB}^{\beta}_{p;\theta}(\mL^p_\omega)}.
\end{align}
Next we estimate $\|\tilde g\|_{\bB^{1+\beta}_{p; \theta}(\mL^p_\omega(\ell^2))}$.
As above, by \eqref{24}, we have
\begin{align*}
\|\bar\sigma\cdot\nabla_{v} u^\delta_{z_0}\|_{\bB^{1+\beta}_{p;\theta}(\mL^p_\omega(\ell^2))}
&\lesssim\|\bar\sigma\|_{\bC^{1+\beta}_{\theta}(\mL^\infty_\omega(\ell^2))}\|\nabla_{v} u^\delta_{z_0}\|_{\bB^1_{p;\theta}(\mL^p_\omega)}+\|\bar\sigma\|_\infty\|\nabla_{v} u^\delta_{z_0}\|_{\bB^{1+\beta}_{p;\theta}(\mL^p_\omega)}\\
&\lesssim\|u^\delta_{z_0}\|_{\bB^2_{p;\theta}(\mL^p_\omega)}+\delta\|u^\delta_{z_0}\|_{\bB^{2+\beta}_{p;\theta}(\mL^p_\omega)},
\end{align*}
and so,
\begin{align}\label{DU2}
\|\tilde g\|_{\bB^{1+\beta}_{p;\theta}(\mL^p_\omega(\ell^2))}
\leq C_0\delta\nor u\nor_{\widetilde{\bB}^{2+\beta}_{p;\theta}(\mL^p_\omega)}+C_\delta\nor u\nor_{\widetilde{\bB}^2_{p;\theta}(\mL^p_\omega)}
+C_1\nor g\nor_{\widetilde{\bB}^{1+\beta}_{p;\theta}(\mL^p_\omega(\ell^2))}.
\end{align}
By \eqref{DU70}, \eqref{DU1}, \eqref{DU2} and Lemma \ref{3gf}, we obtain that for any $q\in[2,\infty]$,  $\delta\in(0,1)$ and $t\in[0,T]$,
\begin{align}
& \nor u(t) \nor_{\widetilde{\bB}^{\beta+2(1-\frac1 q)}_{p; \theta}(\mL^p_\omega)}
=\sup_{z_0\in\mR^d}\|u^\delta_{z_0}(t)\|_{\bB^{\beta+2(1-\frac1 q)}_{p; \theta}(\mL^p_\omega)}\no\\
&\quad\lesssim \left(\int^t_0 \Big[C_0\delta^\beta\nor u(s)\nor^q_{\widetilde{\bB} ^{2+\beta}_{p;\theta}(\mL^p_\omega)}
+C_\delta\nor u(s)\nor^q_{\widetilde{\bB} ^{2}_{p;\theta}(\mL^p_\omega)}\Big]\dif s\right)^{1/q}+{\bf K}_T,\label{VV}
\end{align}
where
$$
{\bf K}_T:=\nor f\nor_{\mL^\infty_T(\widetilde{\bB} ^{\beta}_{p;\theta}(\mL^p_\omega))}
+\nor g\nor_{\mL^\infty_T(\widetilde{\bB} ^{1+\beta}_{p;\theta}(\mL^p_\omega(\ell^2)))}.
$$
In particular, taking $q=\infty$ and $\delta$ being small enough, we obtain for all $t\in[0,T]$,
\begin{align}
 \nor u \nor_{\mL^\infty_t(\widetilde{\bB}^{2+\beta}_{p; \theta}(\mL^p_\omega))}
\lesssim_C
\nor u\nor_{\mL^\infty_t(\widetilde{\bB} ^{2}_{p;\theta}(\mL^p_\omega))}+{\bf K}_T.\label{vv}
\end{align}
Substituting this into  \eqref{VV}  and taking $q=\frac2\beta$, we obtain 
\begin{align*}
 \nor u \nor_{\mL^\infty_t(\widetilde{\bB}^{2}_{p; \theta}(\mL^p_\omega))}
 &\lesssim \left(\int^t_0
\nor u\nor^q_{\mL^\infty_s(\widetilde{\bB} ^{2}_{p;\theta}(\mL^p_\omega))}\dif s\right)^{1/q}+{\bf K}_T,
\end{align*}
which implies by Gronwall's inequality that
$$
 \nor u \nor_{\mL^\infty_T(\widetilde{\bB}^{2}_{p; \theta}(\mL^p_\omega))}\lesssim{\bf K}_T.
$$
Substituting this into \eqref{vv}, we obtain \eqref{ES1}.

(iii) We use the continuity method to show the existence of weak solutions.
Let
$$
a_\lambda:=(\lambda a+(1-\lambda)\mI),\ \ \sigma^k_\lambda:=\lambda\sigma^{k},\ \lambda\in[0,1].
$$
Note that $a_\lambda$, $\sigma_\lambda$ satisfy {\bf (H$_0$)} and {\bf (H$^\beta_1$)} uniformly in $\lambda\in[0,1]$. For instance,
$$
2a_\lambda-\sigma_\lambda\sigma^*_\lambda=2(1-\lambda)\mI+\lambda(2a-\lambda\sigma\sigma^*)
\geq (2(1-\lambda)+\lambda c_0)\mI\geq(2\wedge c_0)\mI.
$$
Given $(f,g)\in\mL^\infty_T(\wt\bB^{\beta}_{p;\theta}(\mL^p_\omega))\times \mL^\infty_T(\wt\bB^{1+\beta}_{p;\theta}(\mL^p_\omega))$,
we consider the following SKE:
\begin{align}\label{SPDE8}
\dif u=[\tr(a_\lambda\cdot\nabla^2_vu)+v\cdot\nabla_x u+\lambda b\cdot\nabla_v u+f]\dif t
+[\sigma^k_\lambda\cdot\nabla_v u+g^k]\dif W^k_t
\end{align}
with $u(0)=u_0$.
Suppose that for some $\lambda_0\in[0,1]$, the above SKE admits a unique weak solution
$u\in \mL^\infty_T(\widetilde{\bB}^{2+\beta}_{p;\theta}(\mL^p_\omega))$
so that \eqref{ES1} holds. We want to show that there is a $\delta$ independent of $\lambda$ such that
for any $\lambda\in[\lambda_0,\lambda_0+\delta]$,
the above SKE still has a unique weak solution $u\in \mL^\infty_T(\widetilde{\bB}^{2+\beta}_{p;\theta}(\mL^p_\omega))$.
Indeed, given $w\in\mL^\infty_T(\wt\bB^{2+\beta}_{p;\theta}(\mL^p_\omega))$, consider SKE
\begin{align}
\dif u&=[\tr(a_{\lambda_0}\cdot\nabla^2_vu)+v\cdot\nabla_x u+\lambda_0 b\cdot\nabla_v u+f\no\\
&\quad+\tr((a_\lambda-a_{\lambda_0})\cdot\nabla^2_vw)+(\lambda-\lambda_0) b\cdot\nabla_v w]\dif t\nonumber\\
&\quad+[\sigma^k_{\lambda_0}\cdot\nabla_v u+(\sigma^k_\lambda-\sigma^k_{\lambda_0})\cdot\nabla_v w+g^k]\dif W^k_t\label{SPDE0}
\end{align}
with $u(0)=u_0$. Note that by \eqref{24},
\begin{align}
\nor\tr((a_\lambda-a_{\lambda_0})\cdot\nabla^2_vw)\nor_{\mL^\infty_T(\wt\bB^{\beta}_{p;\theta}(\mL^p_\omega))}
&\lesssim |\lambda-\lambda_0|\cdot\nor w\nor_{\mL^\infty_T(\wt\bB^{2+\beta}_{p;\theta}(\mL^p_\omega))},\label{q1}\\
\nor(\lambda-\lambda_0)b\cdot\nabla_vw\nor_{\mL^\infty_T(\wt\bB^{\beta}_{p;\theta}(\mL^p_\omega))}
&\lesssim |\lambda-\lambda_0|\cdot\nor w\nor_{\mL^\infty_T(\wt\bB^{1+\beta}_{p;\theta}(\mL^p_\omega))},\label{q2}
\end{align}
and
\begin{align}
\nor(\sigma_\lambda-\sigma_{\lambda_0})\cdot\nabla_vw\nor_{\mL^\infty_T(\wt\bB^{1+\beta}_{p;\theta}(\mL^p_\omega))}
\lesssim |\lambda-\lambda_0|\cdot\nor w\nor_{\mL^\infty_T(\wt\bB^{2+\beta}_{p;\theta}(\mL^p_\omega))}.\label{q3}
\end{align}
By the assumption, SKE \eqref{SPDE0} admits a unique solution denoted by $Q_\lambda w$. In other words, we obtain a
mapping
$$
\mL^\infty_T(\wt\bB^{2+\beta}_{p;\theta}(\mL^p_\omega))\ni w\to Q_\lambda w\in\mL^\infty_T(\wt\bB^{2+\beta}_{p;\theta}(\mL^p_\omega)).
$$
Moreover, by the linearity of SKE \eqref{SPDE0} and \eqref{ES1}, \eqref{q1}-\eqref{q3}, there is a constant $C_0>0$ 
independent of $\lambda_0$ and $\lambda$ such that
$$
\nor Q_\lambda(w-\bar w)\nor_{\mL^\infty_T(\wt\bB^{2+\beta}_{p;\theta}(\mL^p_\omega))}
\leq C_0 |\lambda-\lambda_0|\nor w-\bar w\nor_{\mL^\infty_T(\wt\bB^{2+\beta}_{p;\theta}(\mL^p_\omega))}.
$$
In particular, if we let $\delta=1/(2C_0)$, we obtain a contraction map $Q_\lambda$. By the fixed point theorem,
for any $\lambda\in[\lambda_0,\lambda_0+\delta]$, there is a unique $u\in  \mL^\infty_T(\wt\bB^{2+\beta}_{p;\theta}(\mL^p_\omega))$ so that
$$
Q_\lambda u=u.
$$
That is, SKE \eqref{SPDE8} is uniquely solvable in $\mL^\infty_T(\wt\bB^{2+\beta}_{p;\theta}(\mL^p_\omega))$.
Since for $\lambda=0$, SKE \eqref{SPDE8} is uniquely solvable (see Theorem \ref{Th326}),
by iteration,  \eqref{SPDE8} is also uniquely solvable for $\lambda=1$.
The proof is complete.
\end{proof}
\subsection{Improvement of the regularity in $x$}
Note that by \eqref{28},
$$
\bB^{2+\beta}_{p;\theta}(\mL^p_\omega)=\bB^{(2+\beta)/3}_{p;x}(\mL^p_\omega)\cap\bB^{2+\beta}_{p;v}(\mL^p_\omega).
$$
For $\beta\in(0,1)$, the regularity of weak solutions in variable $x$ in Theorem \ref{Th41} does not exceed $1$. Naturally,
we may ask under what weakest conditions
on the coefficients, weak solutions $u$ of SKE \eqref{spde2} 
become strong solutions, i.e., $u$ is $C^1$-differentiable in variable $x$.
For this aim, we make the following assumptions:
\begin{enumerate}[{\bf (H$^{\beta,\gamma}_1$)}]
\item We suppose that for some $\gamma,\beta\in(0,1)$,
$$
\|a\|_{\mL^\infty_{T,\omega}(\bC^\gamma_x(\bC^{\beta}_{\theta}))}+\|b\|_{\mL^\infty_{T,\omega}(\bC^\gamma_x(\bC^{\beta}_{\theta}))}
+\|\sigma\|_{\mL^\infty_{T,\omega}(\bC^\gamma_x(\bC^{1+\beta}_{\theta}))}\leq c_3.
$$
\end{enumerate}
We also introduce the following mixed space: for $\beta,\gamma>0$ and $p\in[1,\infty]$,
$$
\|f\|_{\bB^{\gamma,\beta}_{p;x,\theta}(\mB)}:=\sup_{j}2^{j\gamma}\|\cR^x_j f\|_{\bB^{\beta}_{p;\theta}(\mB)}
=\sup_{j,k}2^{j\gamma}2^{k\beta}\|\cR^x_j\cR^\theta_k f\|_{\mL^p_x(\mB)}<\infty,
$$
where $\cR_j^x$ is defined in \eqref{DD9}, and
$$
\nor f\nor_{\wt\bB^{\gamma,\beta}_{p;x,\theta}(\mB)}:=\sup_z\|\chi^\delta_z f\|_{\bB^{\gamma,\beta}_{p;x,\theta}(\mB))}<\infty,
$$
where $\chi^\delta_z$ is defined in \eqref{Cu1}.
We have
\bt
\label{es}
Let $T>0$, $p\in[2,\infty)$, $\beta,\gamma\in(0,1)$ and
$$
u_0\in\widetilde{\bB}^{\gamma,2+\beta}_{p;x,\theta}(\mL^p_\omega),\
f\in\mL^\infty_T(\widetilde{\bB}^{\gamma,\beta}_{p;x,\theta}(\mL^p_\omega)),\
g\in\mL^\infty_T(\widetilde{\bB}^{\gamma,1+\beta}_{p;x,\theta}(\mL^p_\omega(\ell^2))).
$$
Under {\bf (H$_0$)} and {\bf (H$^{\beta,\gamma}_1$)}, 
the unique weak solution $u$ in Theorem \ref{Th41} also satisfies
\begin{align*}
\begin{split}
\nor u\nor_{\mL^\infty_T(\widetilde{\bB}^{\gamma,2+\beta}_{p;x,\theta}(\mL^p_\omega))}
&\lesssim_C\nor u_0\nor_{\widetilde{\bB} ^{\gamma,2+\beta}_{p;x,\theta}(\mL^p_\omega)}+
\nor f\nor_{\mL^\infty_T(\widetilde{\bB} ^{\gamma,\beta}_{p;x,\theta}(\mL^p_\omega))}
+\nor g\nor_{\mL^\infty_T(\widetilde{\bB}^{\gamma,1+\beta}_{p;x,\theta}(\mL^p_\omega(\ell^2)))},
\end{split}
\end{align*}
where $C=C(d,p,\beta,\gamma,T,c_0,c_1,c_3)>0$.
\et

\begin{proof}
Following the proof of Theorem \ref{Th41}, noting that
\begin{align}\label{Eq0}
\begin{split}
\dif \cR^x_j u^\delta_{z_0}&=\big[\tr (a_{z_0}\cdot\nabla^2_v \cR^x_j u^\delta_{z_0})+v\cdot\nabla_x \cR^x_j u^\delta_{z_0}+\cR^x_j\tilde f\big]\dif t\\
&\quad+(\cR^x_j\tilde g^k+\sigma^{k}_{z_0}\cdot\nabla_{v} \cR^x_j u^\delta_{z_0})\dif W^k_t,
\end{split}
\end{align}
for any $q\in[2,\infty]$, by Theorem \ref{Th31} with $\alpha=\beta+2(1-\frac1q)$ and $\kappa=1$, we have 
\begin{align}\label{DM9}
\begin{split}
&\|\cR^x_ju^\delta_{z_0}(t)\|_{\bB^{\beta+2(1-\frac1 q)}_{p; \theta}(\mL^p_\omega)}
\lesssim_C \|\cR^x_j(u_0\chi^\delta_{z_0}(0,\cdot))\|_{\bB^{\beta+2(1-\frac1 q)}_{p;\theta}(\mL^p_\omega)}\\
&\qquad+\left(\int^t_0
\Big[\|\cR^x_j\tilde f(s)\|^q_{\bB^{\beta}_{p;\theta}(\mL^p_\omega)}+\|\cR^x_j\tilde g(s)\|^q_{\bB^{\beta+1}_{p;\theta}(\mL^p_\omega(\ell^2))}\Big]\dif s\right)^{1/q}.
\end{split}
\end{align}
We estimate $\|\cR^x_j\tilde f(s)\|_{\bB^{\beta}_{p;\theta}(\mL^p_\omega)}$ and $\|\cR^x_j\tilde g(s)\|_{\bB^{\beta+1}_{p;\theta}(\mL^p_\omega(\ell^2))}$. 
Note that
$$
\cR^x_j\tr(\bar a\cdot\nabla^2_vu^\delta_{z_0})=\tr(\bar a\cdot\cR^x_j\nabla^2_vu^\delta_{z_0})
+\tr ([\cR^x_j, \bar a]\cdot\nabla^2_vu^\delta_{z_0})
$$
and
$$
\cR^x_j(\bar \sigma\cdot\nabla_vu^\delta_{z_0})=\bar\sigma\cdot\cR^x_j\nabla_vu^\delta_{z_0}
+[\cR^x_j, \bar \sigma]\cdot\nabla_vu^\delta_{z_0}.
$$
By Lemma \ref{Le28},  {\bf (H$^{\beta,\gamma}_1$)} and Young's inequality, it is easy to see that
\begin{align*}
&\|\cR^x_j\tr(\bar a\cdot\nabla^2_vu^\delta_{z_0})\|_{\bB^{\beta}_{p;\theta}(\mL^p_\omega)}
+\|\cR^x_j(\bar b\cdot\nabla_v u^{2\delta}_{z_0})\|_{\bB^{\beta}_{p;\theta}(\mL^p_\omega)}
+\|\cR^x_j(\bar c\cdot u^{2\delta}_{z_0})\|_{\bB^{\beta}_{p;\theta}(\mL^p_\omega)}\\
&\qquad\lesssim  \delta^\beta\|\cR^x_ju^\delta_{z_0}\|_{\bB^{2+\beta}_{p;\theta}(\mL^p_\omega)}
+\|\cR^x_ju^\delta_{z_0}\|_{\mL^p_{z,\omega}}
+2^{-\gamma j}\|u^{2\delta}_{z_0}\|_{\bB^{2+\beta}_{p;\theta}(\mL^p_\omega)},
\end{align*}
and
\begin{align*}
&\| \cR^x_j(\bar\sigma\cdot\nabla_{v} u^\delta_{z_0})\|_{\bB^{1+\beta}_{p;\theta}(\mL^p_\omega(\ell^2))}
+\|\cR^x_j(\bar h\, u^{2\delta}_{z_0})\|_{\bB^{1+\beta}_{p;\theta}(\mL^p_\omega(\ell^2))}\\
&\quad\lesssim\delta\|\cR^x_ju^{\delta}_{z_0}\|_{\bB^{2+\beta}_{p;\theta}(\mL^p_\omega)}
+\| \cR^x_j u^{\delta}_{z_0}\|_{\mL^p_{z,\omega}} +2^{-\gamma j}\|u^{2\delta}_{z_0}\|_{\bB^{1+\beta}_{p;\theta}(\mL^p_\omega)}.
\end{align*}
Recall the definitions of $\tilde f$ and $\tilde g$ in \eqref{FF0} and \eqref{FF1}. 
Substituting these two estimates into \eqref{DM9}, we obtain
\begin{align*}
&\|\cR^x_ju^\delta_{z_0}(t)\|_{\bB^{\beta+2(1-\frac1 q)}_{p; \theta}(\mL^p_\omega)}
\lesssim_C \|\cR^x_j(u_0\chi^\delta_{z_0}(0,\cdot))\|_{\bB^{\beta+2(1-\frac1 q)}_{p;\theta}(\mL^p_\omega)}\\
&\qquad+ \left(\int^t_0\Big[\delta^\beta\|\cR^x_ju^\delta_{z_0}\|_{\bB^{2+\beta}_{p;\theta}(\mL^p_\omega)}+
\|\cR^x_ju^\delta_{z_0}(s)\|^q_{\mL^p_z(\mL^p_\omega)} \Big]\dif s\right)^{1/q}\\
&\qquad\qquad +2^{-\gamma j} \left(\int^t_0\|u^{2\delta}_{z_0}(s)\|^q_{\bB^{2+\beta}_{p;\theta}(\mL^p_\omega)}\dif s\right)^{1/q}\\
&\qquad+\left(\int^t_0
\Big[\|\cR^x_jf^\delta_{z_0}(s)\|^q_{\bB^{\beta}_{p;\theta}(\mL^p_\omega)}+\|\cR^x_jg^\delta_{z_0}(s)\|^q_{\bB^{\beta+1}_{p;\theta}(\mL^p_\omega(\ell^2))}\Big]\dif s\right)^{1/q},
\end{align*}
which further implies by definition and  \eqref{ES1} that for any $q\in[2,\infty]$,
\begin{align}\label{q5}
\begin{split}
&\nor u(t)\nor_{\wt\bB^{\gamma,\beta+2(1-\frac1 q)}_{p;x, \theta}(\mL^p_\omega)}
\lesssim_C \nor u_0\nor_{\wt\bB^{\gamma,\beta+2(1-\frac1 q)}_{p;x,\theta}(\mL^p_\omega)}\\
&\quad+ \left(\int^t_0
\Big[\delta^\beta\nor u(s)\nor^q_{\wt\bB^{\gamma,2+\beta}_{p;x,\theta}(\mL^p_\omega)}+C_\delta\nor u(s)\nor^q_{\wt\bB^{\gamma,1}_{p;x, \theta}(\mL^p_\omega)} 
\Big]\dif s\right)^{1/q}+{\bf K}_T,
\end{split}
\end{align}
where
$$
{\bf K}_T:=\nor f\nor_{\mL^\infty_T(\widetilde{\bB}^{\gamma,\beta}_{p;x,\theta}(\mL^p_\omega)))}
+\nor g\nor_{\mL^\infty_T(\widetilde{\bB}^{\gamma,1+\beta}_{p;x,\theta}(\mL^p_\omega(\ell^2)))}.
$$
In particular, letting $q=\infty$ and $\delta$ be small enough, we get
\begin{align}\label{q6}
\nor u\nor_{\mL^\infty_t(\widetilde{\bB}^{\gamma,2+\beta}_{p;x,\theta}(\mL^p_\omega)))}
&\lesssim_C\nor u_0\nor_{\wt\bB^{\gamma,\beta+2}_{p;x,\theta}(\mL^p_\omega)}+
 \nor u\nor_{\mL^\infty_t(\widetilde{\bB}^{\gamma,1}_{p;x,\theta}(\mL^p_{\omega}))} +  {\bf K}_T.
\end{align}
Substituting this into \eqref{q5} and taking $q=2$, we have
\begin{align*}
\nor u\nor_{\mL^\infty_t(\widetilde{\bB}^{\gamma,1+\beta}_{p;x,\theta}(\mL^p_\omega)))}
&\lesssim_C \nor u_0\nor_{\wt\bB^{\gamma,\beta+2}_{p;x,\theta}(\mL^p_\omega)}+ \left(\int^t_0
\nor u\nor^2_{\mL^\infty_s(\wt\bB^{\gamma,1+\beta}_{p;x,\theta}(\mL^p_\omega)}\dif s\right)^{1/2}+ {\bf K}_T.
\end{align*}
By Gronwall's inequality, we get
$$
\nor u\nor_{\mL^\infty_T(\wt\bB^{\gamma,1+\beta}_{p;x,\theta}(\mL^p_\omega)}\lesssim
 \nor u_0\nor_{\wt\bB^{\gamma,\beta+2}_{p;x,\theta}(\mL^p_\omega)}+ {\bf K}_T.
$$
Substituting this into \eqref{q6}, we obtain the desired estimate.
\end{proof}

\subsection{Solvability of adjoint equation \eqref{spde22}}
In this subsection we study the adjoint equation \eqref{spde22}, and assume
\begin{enumerate}[{\bf (H$^\beta_1$)$'$}]
\item We suppose that for some $\beta\in(0,1)$ and $c_2>0$,
$$
\|a\|_{\mL^\infty_{T,\omega}(\bC^{\beta}_{\theta})}+\|b\|_{\mL^\infty_{T,\omega,z}}
+\|\sigma\|_{\mL^\infty_{T,\omega}(\bC^{\beta}_{\theta}(\ell^2))}\leq c_2,\ T>0.
$$
\end{enumerate}
We also introduce the following notion about the solution of SKE \eqref{spde22}.
\bd\label{Def401}
Let $T>0$, $\beta\in(0,1)$ and $p\in[2,\infty)$.
A predictable process $u\in \mL^\infty_T(\wt\bB^{\beta}_{p;\theta}(\mL^p_\omega))$ 
is called a weak solution of SKE \eqref{spde22} if for any $\varphi\in C^\infty_c(\mR^{2d})$,
\begin{align*}
\<u(t),\varphi\>&=\<u_0,\varphi\>+\int^t_0[\<u,\sL_v\varphi+v\cdot\nabla_x\varphi\>+\<f,\varphi\>]\dif s\\
&\quad+\int^t_0[\<u, \sM^k_v\varphi\>+\<g^k,\varphi\>]\dif W^k_s.
\end{align*}
\ed

We have the following result.
\bt\label{Th45}
Let $T>0$, $p\in[2,\infty)$, $\beta\in(0,1)$ and
$$
u_0\in\widetilde{\bB}^{\beta}_{p;\theta}(\mL^p_\omega),\ \ f\in\mL^\infty_T(\widetilde{\bB}^{\beta-2}_{p;\theta}(\mL^p_\omega)),\ \ g\in\mL^\infty_T(\widetilde{\bB}^{\beta-1}_{p;\theta}(\mL^p_\omega(\ell^2))).
$$
Under {\bf (H$_0$)} and {\bf (H$^\beta_1$)$'$},
there is a unique weak solution $u$ to SKE \eqref{spde22} so that
\begin{align}\label{ES01}
\begin{split}
\nor u\nor_{\mL^\infty_T(\widetilde{\bB}^{\beta}_{p;\theta}(\mL^p_\omega))}
&\lesssim_C
\nor u_0\nor_{\widetilde{\bB} ^{\beta}_{p;\theta}(\mL^p_\omega)}+
\nor f\nor_{\mL^\infty_T(\widetilde{\bB} ^{\beta-2}_{p;\theta}(\mL^p_\omega))}
+\nor g\nor_{\mL^\infty_T(\widetilde{\bB}^{\beta-1}_{p;\theta}(\mL^p_\omega(\ell^2)))},
\end{split}
\end{align}
where the constant $C=C(d,p,q,\beta,T,c_0,c_1,c_2)>0$.
\et
\begin{proof}
We only prove the a priori estimate \eqref{ES01}.
We follow the proof of Theorem \ref{Th41} and use the same notations therein.
It is the same reason as in the proof of Theorem \ref{Th41} that we may assume $u_0\equiv0$.
By Definition \ref{Def401}, one sees that
$u^\delta_{z_0}$ is a weak solution of the following freezing SKE:
$$
\dif u^\delta_{z_0}=\Big[\tr (a_{z_0}\cdot\nabla^2_v u^\delta_{z_0})-v\cdot\nabla_x u^\delta_{z_0}
+\tilde f\Big]\dif t+\Big[\sigma^{k}_{z_0}\cdot\nabla_{v} u^\delta_{z_0}+\tilde g^k\Big]\dif W^k_t,
$$
where
\begin{align*}
\tilde f&:=\p^2_{v_iv_j} ((a^{ij}-a^{ij}_{z_0})u^\delta_{z_0})-2\p_{v_i}(a^{ij}u)\p_{v_j}\chi^\delta_{z_0}+\div_v(b u^\delta_{z_0})\\
&\quad-\big[\tr (a\cdot\nabla^2_v\chi^\delta_{z_0}) +(v-v_0)\cdot\nabla_x\chi^\delta_{z_0}+b\cdot\nabla_v\chi^\delta_{z_0}\big]u+f\chi^\delta_{z_0},
\end{align*}
and
\begin{align*}
\tilde g&:=\div_v((\sigma-\sigma_{z_0})u^\delta_{z_0})-\sigma\cdot\nabla_v\chi^\delta_{z_0} u+g\chi^\delta_{z_0}.
\end{align*}
For simplicity of notations, if we set
$$
\bar a:=(a-a_{z_0})\chi^{2\delta}_{z_0},\ \ \bar\sigma:=(\sigma-\sigma_{z_0})\chi^{2\delta}_{z_0},\ \bar h:=\sigma\cdot\nabla_v\chi^\delta_{z_0},
$$
and
$$
\bar c:=-\big(\tr (a\cdot\nabla^2_v\chi^\delta_{z_0}) +(v-v_0)\cdot\nabla_x\chi^\delta_{z_0}+b\cdot\nabla_v\chi^\delta_{z_0}\big),
$$
then due to $\chi^{2\delta}_{z_0}=1$ on the support of $\chi^\delta_{z_0}$, we can write
\begin{align*}
\tilde f=\p^2_{v_iv_j} (\bar a^{ij}u^\delta_{z_0})-2\p_{v_i}(a^{ij}u^{2\delta}_{z_0})\p_{v_j}\chi^\delta_{z_0}+\div_v(b u^\delta_{z_0})
+\bar c\, u^{2\delta}_{z_0}
+f\chi^\delta_{z_0}
\end{align*}
and
\begin{align*}
\tilde g=\div_v(\bar\sigma u^\delta_{z_0})-\bar h\, u^{2\delta}_{z_0}+g\chi^\delta_{z_0}.
\end{align*}
For any $p\in[2,\infty)$ and $q\in(2,\infty]$, by Theorem \ref{Th31}, there is a constant $C=C(T,d,p,q,\beta,c_0,c_1)>0$
such that for all $\delta\in(0,1)$ and $z_0\in\mR^{2d}$,
\begin{align}\label{DU0}
\begin{split}
\|u^\delta_{z_0}\|_{\mL^\infty_T(\bB^{\beta-\frac2q}_{p; \theta}(\mL^p_\omega))}
&\lesssim_C
\|\tilde f\|_{\mL^q_T(\bB^{\beta-2}_{p;\theta}(\mL^p_\omega))}
+\|\tilde g\|_{\mL^q_T(\bB^{\beta-1}_{p; \theta}(\mL^p_\omega(\ell^2)))}.
\end{split}
\end{align}
We estimate each term of the right hand side. 
By Bernstein's inequality and \eqref{24},
\begin{align*}
\|\p^2_{v_iv_j} (\bar a^{ij}u^\delta_{z_0})\|_{\bB^{\beta-2}_{p;\theta}(\mL^p_\omega)}\lesssim
\|\bar a u^\delta_{z_0}\|_{\bB^{\beta}_{p;\theta}(\mL^p_\omega)}
\lesssim\delta^\beta\|u^\delta_{z_0}\|_{\bB^{\beta}_{p;\theta}(\mL^p_\omega)}
+C_\delta\|u^\delta_{z_0}\|_{\mL^p_z(\mL^p_\omega)},
\end{align*}
and by \eqref{Ea1},
\begin{align*}
\|\p_{v_i}(a^{ij}u^{2\delta}_{z_0})\p_{v_j}\chi^\delta_{z_0}\|_{\bB^{\beta-2}_{p;\theta}(\mL^p_\omega)}
&\lesssim\|\p_{v_i}(a^{ij}u^{2\delta}_{z_0})\|_{\bB^{\beta-2}_{p;\theta}(\mL^p_\omega)}\|\p_{v_j}\chi^\delta_{z_0}\|_{\bC^{2}_{\theta}}\\
&\lesssim\|a^{ij}u^{2\delta}_{z_0}\|_{\mL^p_z(\mL^p_\omega)}\lesssim \|u^{2\delta}_{z_0}\|_{\mL^p_z(\mL^p_\omega)},
\end{align*}
and
$$
\|\div_v(b u^\delta_{z_0})\|_{\bB^{\beta-2}_{p;\theta}(\mL^p_\omega)}\lesssim \|b u^\delta_{z_0}\|_{\bB^{\beta-1}_{p;\theta}(\mL^p_\omega)}
\leq \|b\|_{\mL^\infty_{z,\omega}}\|u^\delta_{z_0}\|_{\mL^p_z(\mL^p_\omega)}.
$$
Since $\bar c$ is bounded, we also have
$$
\|\bar c u^{2\delta}_{z_0}\|_{\bB^{\beta-2}_{p;\theta}(\mL^p_\omega)}\lesssim \|\bar c u^{2\delta}_{z_0}\|_{\mL^p_z(\mL^p_\omega)}
\lesssim \|u^{2\delta}_{z_0}\|_{\mL^p_z(\mL^p_\omega)}.
$$
Therefore,
\begin{align*}
\|\tilde f\|_{\bB^{\beta-2}_{p;\theta}(\mL^p_\omega)}
\lesssim \delta^\beta\nor u\nor_{\wt\bB^{\beta}_{p;\theta}(\mL^p_\omega)}
+C_\delta\nor u\nor_{\wt\mL^p_z(\mL^p_\omega)}+\nor f\nor_{\wt\bB^{\beta-2}_{p;\theta}(\mL^p_\omega)}.
\end{align*}
Similarly, one can show that
\begin{align*}
\|\tilde g\|_{\bB^{\beta-1}_{p;\theta}(\mL^p_\omega(\ell^2))}
\lesssim \delta^\beta\nor u\nor_{\wt\bB^{\beta}_{p;\theta}(\mL^p_\omega)}
+C_\delta\nor u\nor_{\wt\mL^p_z(\mL^p_\omega)}+\nor g\nor_{\wt\bB^{\beta-1}_{p;\theta}(\mL^p_\omega(\ell^2))}.
\end{align*}
Substituting these two estimates into \eqref{DU0}, and taking supremum in $z_0\in\mR^{2d}$, we obtain
that for any $q\in[2,\infty]$,
\begin{align*}
\nor u\nor_{\mL^\infty_T(\wt\bB^{\beta-\frac2q}_{p; \theta}(\mL^p_\omega))}
&\lesssim_C\delta^\beta\nor u\nor_{\mL^q_T(\wt\bB^{\beta}_{p;\theta}(\mL^p_\omega))}
+C_\delta\nor u\nor_{\mL^q_T(\wt\mL^p_z(\mL^p_\omega))}\\
&\quad+\nor f\nor_{\mL^q_T(\wt\bB^{\beta-2}_{p;\theta}(\mL^p_\omega))}+\nor g\nor_{\mL^q_T(\wt\bB^{\beta-1}_{p; \theta}(\mL^p_\omega(\ell^2)))}.
\end{align*}
The remaining proof is completely the same as in Theorem \ref{Th41}.
\end{proof}

\section{Applications in nonlinear filtering problems}

Fix $d,d_1\in\mN$. 
Let $(\tilde B_t)_{t\geq 0}$ and $(B_t)_{t\geq 0}$ be two independent $d$ and $d_1$-dimensional Brownian motions
on some stochastic basis $(\Omega,\mathcal{F},\mathbb{P};(\mathcal{F}_t)_{t\geq 0})$. 
Let
$$
(\tilde b,\tilde\sigma,\sigma,\hat b): \mR_+\times\mR^{2d}\times\mR^{d_1}\to(\mR^d,\mR^d\otimes\mR^d,\mR^d\otimes\mR^{d_1},\mR^{d_1})
$$
and
$$
\hat \sigma: \mR_+\times\mR^{d_1}\to\mR^{d_1}\otimes
\mR^{d_1}
$$
be bounded Borel measurable functions. 
Consider the following SDE of  It\^o's type:
\begin{equation}\label{sde8}
\dif\!
\left(
\begin{array}{c}
X_t\\
V_t\\
Y_t
\end{array}
\right)
=
\left(
\begin{array}{c}
V_t,\\
 \tilde b(t,U_t)\\
\hat b(t,U_t)
\end{array}
\right)\dif t+
\left(
\begin{array}{ccc}
0,&0,&0\\
0,&\tilde\sigma(t,U_t),&\sigma(t,U_t)\\
0, &0,&\hat \sigma(t,Y_t)
\end{array}
\right)\dif\!
\left(
\begin{array}{c}
0\\
\tilde B_t\\
B_t
\end{array}
\right),
\end{equation}
where 
$$
U_t:=(X_t,V_t,Y_t)=:(Z_t, Y_t)\in\mR^{2d}\times\mR^{d_1},
$$
$Z_t$ stands for the {\it unobservable} signal and $Y_t$ denotes the {\it observable} signal,
$U_0$ is an $\mathcal{F}_0$-measurable random variable.
Let $\sF^Y_t$ be the $\mP$-completed $\sigma$-algebra generated by $\{Y_s,s\leq t\}$, which represents the observation information.
In application, we want to predict $Z_t$ in terms of $\sF^Y_t$. More precisely, we want to calculate the conditional 
distribution of $Z_t$ under $\sF^Y_t$:
$$
\Pi_t(\omega,A):=\mP(Z_t\in A|\sF^Y_t).
$$
This is usually called the filtering problem (cf. \cite{roz}).

Throughout this section we suppose that 
\begin{enumerate}[$(\bf {A}_1)$]
\item $\Pi_0(\omega,\dif z)=\pi_0(\omega,z)\dif z$ with
$\pi_0\in \cap_{p\geq 2}\widetilde{\bB}^{\beta}_{p;\theta}(\mL^p_\omega)$ for some $\beta\in(0,1)$.

\item  $\tilde b, \tilde\sigma, \sigma, \hat b, \hat \sigma$ are  bounded and Lipschitz in $(x,v,y)$ and uniformly in $t$. 

\item $\tilde\sigma$ and $\hat \sigma$ are non-singular, and for some $K>0$ and all $t,x,v,y$,
$$
\|\tilde\sigma^{-1}(t,x,v,y)\|_{\mM^{d}_{sym}}+\|\hat \sigma^{-1}(t,y)\|_{\mM^{d_1}_{sym}}\leq K.
$$
\end{enumerate}
Under the above assumptions, it is well-known that SDE \eqref{sde8} admits a unique strong solution.
Let
\begin{align*}
\tilde h(t,z,y)&:=(\hat \sigma^{-1}\hat b)(t,z,y),\ \bar b(t,z,y):=(\tilde b-\sigma\tilde h)(t,z,y)
\end{align*}
and 
$$
W_t=B_t+\int^t_0 \tilde h(s,U_s)\dif s.
$$
Then we can write SDE \eqref{sde8} as
\begin{equation*}
\dif
\left(
\begin{array}{c}
X_t\\
V_t\\
Y_t
\end{array}
\right)
=
\left(
\begin{array}{c}
V_t,\\
 \bar b(t,U_t)\\
0
\end{array}
\right)\dif t+
\left(
\begin{array}{ccc}
0,&0,&0\\
0,&\tilde\sigma(t,U_t),&\sigma(t,U_t)\\
0, &0,&\hat \sigma(t,Y_t)
\end{array}
\right)\dif\!
\left(
\begin{array}{c}
0\\
\tilde B_t\\
W_t
\end{array}
\right).
\end{equation*}
Define
\begin{align}\label{RH0}
 \rho_t:=\exp\left\{\int^t_0\tilde h(s,U_s)\dif  B_s+\frac{1}{2}\int^t_0\|\tilde h (s,U_s) \|^2\dif s\right\}.
\end{align}
Fix $T>0$. Since $\tilde h$ is bounded, by It\^o's formula,  it is easy to see that
\begin{align}\label{HF6}
\mE(\rho^{-1}_T|\sF_t)=\rho^{-1}_t,\ \ t\in [0,T].
\end{align}
Thus,  by Girsanov's theorem, under the new probability measure 
\begin{align}\label{HF5}
\bar{\mathbb{P}}_T(\dif\omega):=\rho^{-1}_T(\omega)\mathbb{P}(\dif \omega),
\end{align} 
$(W_t)_{t\in[0,T]}$ is still a $d_1$-dimensional  Wiener process and independent of $\tilde B$.
Moreover, since $\hat \sigma$ is invertible, it is easy to see that (see \cite[Lemma 6.2]{roz})
\begin{align}\label{DR1}
\sF^Y_t=\sF^{W}_t\vee \sF^Y_0.
\end{align}
This is the crucial point for deriving the density equation of $\Pi_t(\omega,\cdot)$. In other words, by Girsanov's technique, the observation
information is transformed into Brownian filtration information.

Now we let
\begin{align}\label{COE}
\begin{split}
\Sigma(t,\omega,z)&:=\sigma\left(t,z,Y_t(\omega)\right),\ \ \tilde\Sigma(t,\omega,z):=\tilde\sigma\left(t,z,Y_t(\omega)\right),\\
b(t,\omega,z)&:=\bar b\left(t,z,Y_t(\omega)\right), \ \
h(t,\omega,z):=(\hat \sigma^{-1}\hat b)(t,z,Y_t(\omega)),
\end{split}
\end{align}
and
\begin{align}\label{COE0}
\sL_v \varphi:=\tr (a\cdot\nabla^2_v\varphi)+b\cdot\nabla_v \varphi,\ \ \sM_v\varphi:=\Sigma\cdot\nabla_v\varphi+h\varphi,
\end{align}
where $a:=\frac12(\Sigma\Sigma^*+\tilde\Sigma\tilde\Sigma^*)$.
Consider the following linear SKE:
\begin{align}\label{spde4}
\dif u=[\sL^*_vu-v\cdot\nabla_x u]\dif t+ \sM^*_v u\dif W_t,\ \ u(0)=\pi_0,
\end{align}
where $\sL^*_v$ and $\sM^*_v$ stand for the adjoint operators of $\sL_v$ and $\sM_v$, respectively.
Under ${\bf (A_1)}$-${\bf (A_3)}$, by Theorem \ref{Th45}, there exists a unique generalized solution 
\begin{align}\label{re50}
u\in\cap_{p\geq 2}\mL^\infty_T(\widetilde{\bB}^{\beta}_{p;\theta}(\mL^p_\omega)), \quad T>0,
\end{align}
to the above SKE. We first show the following result.
\bl
For each $t\geq 0$, $R>0$ and $\beta'<\beta$, it holds that  for $\mP$-almost all $\omega$, 
\begin{align}\label{CT2}
\chi^R_0(\cdot) u(t,\omega,\cdot)\in\cap_{p\geq 2}\bB^{\beta'}_{p;\theta},
\end{align}
where $\chi^R_0$ is the cutoff function in \eqref{Cu1}, and for any bounded measurable $\varphi$,
\begin{align}
\<u(t,\omega),\varphi\>=\mathbb{E}^{\bar{\mathbb{P}}_T}\left(\varphi(Z_t)\rho_T| \sF^Y_t\right)(\omega),\label{itr}
\end{align}
where $\rho_T$ is defined by \eqref{RH0}.
\el
\begin{proof}
First of all, for each $R, t>0$, by \eqref{re50} and \eqref{CT1}, we have
$$
\chi^R_0u(t)\in \cap_{p\geq 2}{\bB}^{\beta}_{p;\theta}(\mL^p_\omega)\subset \cap_{p\geq 2}\mL^p_\omega({\bB}^{\beta'}_{p;\theta}),
$$
which in turn gives \eqref{CT2}.
Let 
$$
(b_\eps,\tilde \Sigma_\eps, \Sigma_\eps, h_\eps)=
(b,\tilde\Sigma,\Sigma, h)*\varrho_\eps,
$$
where $(\varrho_\eps)_{\eps\in(0,1)}$ is a family of mollifiers in $\mR^{2d}$.
Let $(\sL^{\eps}_v,\sM^\eps_v)$ be defined as in \eqref{COE0} in terms of the above mollifying coefficients.
Let $u_\eps$ be the unique solution of the following SKE:
$$
\dif u_\eps=[(\sL^\eps_v)^* u_\eps-v\cdot\nabla_x u_\eps] \dif t+ (\sM^\eps_v)^* u_\eps\dif W_t,\ \ u_\eps(0)=\pi_0.
$$
By \eqref{ES01}, there is a constant $C>0$ independent of $\eps\in(0,1)$ such that
\begin{align}\label{ES801}
\nor u_\eps\nor_{\mL^\infty_T(\widetilde{\bB}^{\beta}_{p;\theta}(\mL^p_\omega))}\lesssim_C
\nor \pi_0\nor_{\widetilde{\bB}^{\beta}_{p;\theta}(\mL^p_\omega)}.
\end{align}
Let 
$$
U_\eps:=u_\eps-u, \  F_\eps:=(\sL^\eps_v-\sL_v)^* u_\eps,\ G_\eps:=(\sM^\eps_v-\sM_v)^* u_\eps.
$$
Then one sees that
$$
\dif U_\eps=[(\sL_v)^* U_\eps-v\cdot\nabla_xU_\eps+F_\eps] \dif t+ [(\sM_v)^* U_\eps+G_\eps]\dif W_t,\ \ U_\eps(0)=0.
$$
By \eqref{ES01}, we have
\begin{align}\label{ES901}
\nor U_\eps\nor_{\mL^\infty_T(\widetilde{\bB}^{\beta}_{p;\theta}(\mL^p_\omega))}\lesssim
\nor F_\eps\nor_{\mL^\infty_T(\widetilde{\bB} ^{\beta-2}_{p;\theta}(\mL^p_\omega))}
+\nor G_\eps\nor_{\mL^\infty_T(\widetilde{\bB}^{\beta-1}_{p;\theta}(\mL^p_\omega(\ell^2)))}.
\end{align}
Here and below, the implicit constant in $\lesssim$ is independent of $\eps\in(0,1)$.
Noting that
$$
F_\eps=\p^2_{v_iv_j}((a^{ij}_\eps-a^{ij})u_\eps)+\p_{v_i}((b^i_\eps-b^i)u_\eps),
$$
by Bernstein's inequality \eqref{Ber}, \eqref{24} and \eqref{ES801}, we have
\begin{align*}
\nor F_\eps\nor_{\widetilde{\bB} ^{\beta-2}_{p;\theta}(\mL^p_\omega)}
&\lesssim\nor (a_\eps-a)u_\eps\nor_{\wt\bB^{\beta}_{p;\theta}(\mL^p_\omega)}
+\nor (b_\eps-b)u_\eps\nor_{\wt\bB^{\beta-1}_{p;\theta}(\mL^p_\omega)}\\
&\lesssim\| a_\eps-a\|_{\mL^\infty_{t,\omega}(\bC^{\beta}_{\theta})}\nor u_\eps\nor_{\wt\bB^{\beta}_{p;\theta}(\mL^p_\omega)}
+\| b_\eps-b\|_{\mL^\infty_{t,\omega,z}}\nor u_\eps\nor_{\wt\mL^p_{z}(\mL^p_\omega)}\\
&\lesssim(\eps^{1-\beta}+\eps)\nor u_\eps\nor_{\wt\bB^{\beta}_{p;\theta}(\mL^p_\omega)}\lesssim\eps^{1-\beta},
\end{align*}
where in the third inequality we have used that for $|z-z'|\leq 1$,
\begin{align*}
|(a_\eps-a)(t,\omega,z)-(a_\eps-a)(t,\omega,z')|\lesssim \|a_\eps-a\|_\infty^{1-\beta}\|a\|^\beta_{\rm Lip}|z-z'|^\beta
\lesssim\eps^{1-\beta}\ |z-z'|^\beta_\theta.
\end{align*}
Similarly we have
\begin{align*}
\nor G_\eps\nor_{\widetilde{\bB} ^{\beta-1}_{p;\theta}(\mL^p_\omega)}
&\lesssim\nor (\Sigma_\eps-\Sigma)u_\eps\nor_{\wt\bB^{\beta}_{p;\theta}(\mL^p_\omega)}
+\nor (h_\eps-h)u_\eps\nor_{\wt\mL^p_z(\mL^p_\omega)}
\lesssim\eps^{1-\beta}.
\end{align*}
Thus, by \eqref{ES901},
\begin{align}\label{DS1}
\lim_{\eps\to 0}\nor U_\eps\nor_{\mL^\infty_T(\widetilde{\bB}^{\beta}_{p;\theta}(\mL^p_\omega))}=0.
\end{align}
Let $Z^\eps:=(X^\eps, V^\eps)$ solve the following SDE
$$
\dif X^\eps_t=V^\eps_t\dif t,\ \ \dif V^\eps_t=b_\eps(t, Z^\eps_t)\dif t+\tilde\Sigma_\eps(t, Z^\eps_t)\dif \tilde B_t+
\Sigma_\eps(t, Z^\eps_t)\dif B_t,\ Z^\eps_0=Z_0.
$$
Since $b_\eps,\Sigma_\eps,\tilde\Sigma_\eps$ are bounded and Lipschitz continuous in $z$ uniformly in $(t,\omega,\eps)$, 
it is by now standard to show that
\begin{align}\label{DS2}
\lim_{\eps\to 0}\mE |Z^\eps_t-Z_t|^2=0.
\end{align}
Now by \cite[p.201, Theorem 5.6]{roz}, we have
for any $\varphi\in C_b^{\infty}(\mathbb{R}^{2d})$,
\begin{align}
\<u_\eps(t,\omega),\varphi\>=\mathbb{E}^{\bar{\mathbb{P}}_T}\left(\varphi(Z^\eps_t)\rho^\eps_t| \sF^{Y}_t\right)(\omega),
\quad \mP-a.s.\label{itr0}
\end{align}
where
$$
\rho^\eps_t:=\exp\left\{\int^t_0h_\eps(s,Z^\eps_s)\dif  B_s+\frac{1}{2}\int^t_0\|h_\eps (s,Z^\eps_s) \|^2\dif s\right\}.
$$
By \eqref{DS1} and \eqref{DS2} and taking limits for \eqref{itr0}, we get
\begin{align}
\<u(t,\omega),\varphi\>=\mathbb{E}^{\bar{\mathbb{P}}_T}\left(\varphi(Z_t)\rho_t| \sF^{Y}_t\right)(\omega),
\quad \mP-a.s.\label{itr10}
\end{align}
On the other hand, for any $A\in\sF^Y_t$, by  \eqref{HF5} and \eqref{HF6} we have
\begin{align*}
\mathbb{E}^{\bar{\mathbb{P}}_T}\left(\1_A\varphi(Z_t)\rho_t\right)&=\mE\left(\1_A\varphi(Z_t)\rho_t\rho^{-1}_T\right)
=\mE\left(\1_A\varphi(Z_t)\mE(\rho_t\rho^{-1}_T|\sF_t)\right)\\
&=\mE\left(\1_A\varphi(Z_t)\right)=\mE\left(\1_A\varphi(Z_t)\rho_T\rho^{-1}_T\right)
=\mE^{\bar \mP_T}\left(\1_A\varphi(Z_t)\rho_T\right),
\end{align*}
which means that
$$
\mathbb{E}^{\bar{\mathbb{P}}_T}\left(\varphi(Z_t)\rho_t| \sF^Y_t\right)=\mathbb{E}^{\bar{\mathbb{P}}_T}\left(\varphi(Z_t)\rho_T| \sF^Y_t\right).
$$
Substituting this into \eqref{itr10}, we obtain \eqref{itr}.
\end{proof}

We now state the main result of this section.
\bt
Under ${\bf (A_1)}$-${\bf (A_3)}$, for each $t>0$, the conditional distribution $\Pi_t(\omega,\cdot)$ has a continuous density
$\pi_t(\omega,\cdot)$ for $\mP$-almost all $\omega$, and
\begin{align}
\pi_t(\omega,\cdot)=u(t,\omega,\cdot)/\|u(t,\omega,\cdot)\|_{\mL^1_z},\label{den}
\end{align}
where $u$ is given in \eqref{re50}.
Moreover, $\pi_t$ solves the following nonlinear SKE:
\begin{align}\label{GQ1}
\dif \pi_t(\varphi)=\pi_t(\sL_v\varphi+v\cdot\nabla_x\varphi)\dif t+(\pi_t(\sM_v\varphi)-\Pi_t(h_t) (\pi_t(\varphi)))\dif \bar W_t,
\end{align}
for any $\varphi\in C^\infty_c(\mR^{2d})$, where $\bar W_t:=W_t+\int^t_0\Pi_s(h_s)\dif s$.
\et
\begin{proof}
Let $T>0$. Note that by \eqref{RH0},
\begin{align*}
 \rho_t=\exp\left\{\int^t_0\tilde h(s,U_s)\dif W_s-\frac{1}{2}\int^t_0\|\tilde h (s,U_s) \|^2\dif s\right\}
\end{align*}
satisfies
\begin{align*}
\rho_t=1+\int^t_0 \tilde h(s,U_s)\rho_s\dif W_s.
\end{align*}
Taking the conditional expectation with respect to $\sF^Y_t$, by  \cite [p.28, Theorem 1.15]{roz}, we have
\begin{align*}
\mE^{\bar\mP_T}(\rho_t|\sF^Y_t)&=1+\mE^{\bar\mP_T}\left(\int^t_0 \tilde h(s,U_s)\rho_s\dif W_s|\sF^Y_t\right)\\
&=1+\int^t_0\mE^{\bar\mP_T}\left( \tilde h(s,U_s)\rho_s|\sF^Y_s\right)\dif W_s\\
&=1+\int^t_0\mE\left( \tilde h(s,U_s)|\sF^Y_s\right)\mE^{\bar\mP_T}\left(\rho_s|\sF^Y_s\right)\dif W_s,
\end{align*}
which together with $\Pi_t(h_t)= \mathbb{E}\left(h_t(Z_t)\mid\sF^Y_t\right)=\mE( \tilde h(t,U_t)|\sF^Y_t)$ yields that
\begin{align}\label{GD1}
\mathbb{E}^{\bar{\mathbb{P}}_T}\left(\rho_t \mid\mathcal{F}^Y_t\right)&= \exp\left\{\int^t_0 \Pi_s(h_s)\dif W_s
  -\frac{1}{2}\int^t_0\| \Pi_s(h_s) \|^2\dif s\right\}.
\end{align}

Now, by Bayes' formula about the conditional expectations (cf. \cite[p.218, Lemma 6.1]{roz}) and \eqref{itr}, we have
\begin{align*}
\Pi_t(\varphi)=\mathbb{E}\left(\varphi(Z_t)\mid\mathcal{F}^Y_t\right)=\frac{\mathbb{E}^{\bar{\mathbb{P}}}\left(\varphi(Z_t)\rho_T\mid\mathcal{F}^Y_t\right)}{\mathbb{E}^{\bar{\mathbb{P}}}\left(\rho_T\mid\mathcal{F}^Y_t\right)}= \frac{\langle u(t),\varphi\rangle}{\langle u(t),1\rangle}
=\left\langle \frac{u(t)}{\<u(t),1\>},\varphi\right\rangle.
\end{align*}
By the arbitrariness of $\varphi$, we get \eqref{den}. The continuity of $\pi_t(\omega,\cdot)$ follows by \eqref{CT2} and the 
Sobolev embedding.

Moreover, by \eqref{itr0} and \eqref{GD1}, 
$$
\langle u(t),1\rangle=\exp\left\{\int^t_0 \Pi_s(h_s)\dif W_s -\frac{1}{2}\int^t_0\| \Pi_s(h_s) \|^2\dif s\right\}
$$
and
$$
\dif\<u(t),\varphi\>=\<u(t),\sL_t\varphi\>\dif t+\<u(t),\sM_t\varphi\>\dif W_t.
$$
Thus, by It\^o's formula,
\begin{align*}
\dif \pi_t(\varphi)= &\left[\pi_t(\mathscr{L}_t\varphi)+\left(\pi_t(\mathscr{M}_t\varphi)-\pi_t(\varphi)\right) \|\Pi_t(h_t)\|^2\right]\dif t\\
&+ \left(\pi_t(\mathscr{M}_t\varphi)-\pi_t(\varphi)\Pi_t(h_t)\right)\dif  {W}_t\\
=& \pi_t(\mathscr{L}_t\varphi)\dif t+
 \left(\pi_t(\mathscr{M}_t\varphi)-\pi_t(\varphi)\Pi_t(h_t)\right)(\dif {W}_t+\Pi_t(h_t)\dif t).
\end{align*}
From these, we obtain \eqref{GQ1} and complete the proof.
\end{proof}

\br\rm
Under further assumptions on $\tilde b, \tilde\sigma, \sigma, b, \hat \sigma$ and $\pi_0$, by Theorem \ref{es}, one can show the $C^1$ and $C^2$-smoothness of $\pi$ with respect to $x$ variable and $v$ variable so that $\pi$ solves the following SKE in the point-wise sense of $x,v$,
\begin{align}\label{GQ1}
\dif \pi_t=(\sL_v^*\pi_t+v\cdot\nabla_x\pi)\dif t+(\sM_v^*\pi_t-\Pi_t(h_t) \pi_t)\dif \bar W_t.
\end{align}
For instance, we assume that for some $\gamma,\beta\in(0,1)$ with $3\gamma+\beta>1$,
\begin{enumerate}[$(\bf {A}'_1)$]
\item For almost all $\omega$, $\Pi_0(\omega,\dif z)=\pi_0(\omega,z)\dif z$ with
$\pi_0\in \cap_{p\geq 2}\widetilde{\bB}^{\gamma,2+\beta}_{p;x,\theta}(\mL^p_\omega)$.

\item  $\tilde b, \tilde\sigma, \sigma, b, \hat \sigma$ are  bounded and Lipschitz continuous in spatial variables $(x,v,y)$ and uniformly in 
time variable $t$. Moreover, we also suppose that
\begin{align*}
&\sigma,\tilde\sigma,\nabla^2_v\sigma,\nabla^2_v\tilde\sigma, \tilde b,\div_v\tilde b,\nabla_v b,\nabla_v\hat \sigma
\in \mL^\infty_{T,y}(\mathbf{C}^{\gamma}_x(\mathbf{C}^{\beta}_\theta))
\end{align*}									
and
$$
\div_v\sigma\in \mL^\infty_{T,y}(\mathbf{C}^{\gamma}_x(\mathbf{C}^{1+\beta}_\theta)).
$$
\end{enumerate}
Thus, under $(\bf {A}'_1)$, $(\bf {A}'_2)$ and $(\bf {A}_3)$, by Theorem \ref{es}, there exists a unique solution to SKE \eqref{spde4} with regularity
\begin{align}\label{re500}
u\in\cap_{p\geq 2}\mL^\infty_T(\widetilde{\bB}^{\gamma,2+\beta}_{p;x,\theta}(\mL^p_\omega))),\ \  T>0.
\end{align}
Due to $3\gamma+\beta>1$, by \eqref{den}, \eqref{re500}, \eqref{CT1} and Sobolev's embedding, we have
$$
\pi\in\cap_{p\geq 2}\mL^\infty_T(\widetilde{\bB}^{\gamma,2+\beta}_{p;x,\theta}(\mL^p_\omega)))
\subset \cap_{p\geq 2} \mL^\infty_T(\mL^p_\omega(C^{1}_{x}\cap C^{2}_{v})),\ \ T>0,
$$
where $C^1_x\cap C^2_v$ stands for the space of $C^1$ and $C^2$-smooth functions in $x$ and $v$.
\er

\section{SKEs driven by velocity-time white noises}

Let $\{B(t,v), (t,v)\in\mR_+\times\mR\}$ be a Brownian sheet of time and velocity variables, 
whose distribution derivative in $t,v$ is usually considered as velocity-time white noise. Let
$$
(f,g)(t,\omega,x,v,u):\mR_+\times\Omega\times\mR^2\times\mR\to (\mR,\mR)
$$
be $\cP\times\cB(\mR^3)$-measurable functions. We consider the following nonlinear SKE driven by $B$,
\begin{align}\label{CN2}
\dif u=\big[\Delta_vu+v\cdot\nabla_x u+f(u)\big]\dif t+g(u)\dif B(t,v),\ \ u(0)=u_0,
\end{align}
where the stochastic integral is understood in the sense of Walsh \cite{Wa}.
Here we have suppressed the variables $(t,\omega, x,v)$ of $f,g$. We introduce the following notion about the solution to the above SKE.
\bd
Let $p\geq 2$ and $T>0$. We call a predictable process $u\in\mL^\infty_T(\wt\mL^p_z(\mL^p_\omega))$ being
a weak solution of SKE \eqref{CN2} if for any $\varphi\in C^\infty_c(\mR^2)$ and $t\in[0,T]$,
\begin{align}\label{DX7}
\begin{split}
\<u(t),\varphi\>&=\<u_0,\varphi\>+\int^t_0\big[\<u,\Delta_v\varphi-v\cdot\nabla_x\varphi\>+\<f(u),\varphi\>\big]\dif s\\
&\quad+\int^t_0\int_{\mR}\left(\int_{\mR}\varphi(x,v) g(s,x,v,u(s,x,v))\dif x\right)\dif B(s,v)\ \ a.s.
\end{split}
\end{align}
\ed
In order to use our previous results to study SKE \eqref{CN2}, 
we shall use the following representation of $B(t,v)$ (cf. \cite[Section 8.2]{kry99}),
\begin{align}\label{DC8}
B(t,v)=\sum_{k=1}^\infty \left(\int^v_0\eta_k(r)\dif r\right) W^k_t,
\end{align}
where $\{\eta_k,k\in\mN\}\subset C^\infty_c(\mR)$ is an orthogonormal basis of $\mL^2_v(\mR)$.
In particular, for any measurable adapted processes $\xi(s,v)$ with $\int^T_0\mE\|\xi(s,\cdot)\|^2_{\mL^2_v}\dif s<\infty$,
\begin{align}\label{DJ3}
\int^T_0\int_{\mR}\xi(s,v)\dif B(s,v)=\sum_{k=1}^\infty\int^T_0\left(\int_\mR \xi(s,v)\eta_k(v)\dif v\right)\dif W^k_s.
\end{align}

For $N=1,\cdots,\infty$ and a function $h$, we introduce the following notations:
$$
\mG_N(h):=(h\eta_1,\cdots,h\eta_N,0,\cdots),\ \ \mG(h):=\mG_\infty(h).
$$
By this notation and \eqref{DJ3}, we can write SKE \eqref{CN2} as 
\begin{align}\label{CN6}
\dif u=\big[\Delta_vu+v\cdot\nabla_x u+f(u)\big]\dif t+\mG^k(g(u))\dif W^k_t,\ \ u(0)=u_0.
\end{align}
We have the following important lemma.
\bl
For any $p\geq 2$ and $r\in[1,p]$, there is a constant $C=C(p,r)>0$ such that
for all $N=1,\cdots,\infty$ and $h\in \mL^p_\omega(\mL^r_z)$,
\begin{align}\label{DJ0}
\|\mG_N(h)\|_{\bB^{\beta}_{p;\theta}(\mL^p_\omega(\ell^2))}\leq C\|h\|_{\mL^p_\omega(\mL^r_z)},
\end{align}
where $\beta:=4(\frac1p-\frac1r)-\frac12$. Moreover, for any $\beta'<\beta$, we also have
\begin{align}\label{DJ00}
\lim_{N\to\infty}\|\mG_N(h)-\mG(h)\|_{\bB^{\beta'}_{p;\theta}(\mL^p_\omega(\ell^2))}=0.
\end{align}
\el
\begin{proof}
Let $h\in \mL^p_\omega(\mL^r_z)$.
Recall \eqref{Ph0}. By Parsaval's identity, we have for any $j\in\mN_0$,
\begin{align*}
\|\cR^\theta_j\mG_N(h)(x,v)\|^2_{\ell^2}&=\sum_{k=1}^N\left|\int_{\mR^2}\check\phi^\theta_j(x-\bar x,v-\bar v) 
h(\bar x,\bar v)\eta_k(\bar v)\dif\bar x\dif\bar v\right|^2\\
&\leq\sum_{k=1}^\infty\left|\int_{\mR}\left(\int_\mR\check\phi^\theta_j(x-\bar x,v-\bar v) 
h(\bar x,\bar v)\dif\bar x\right)\eta_k(\bar v)\dif\bar v\right|^2\\
&=\int_{\mR}\left|\int_{\mR}\check\phi^\theta_j(x-\bar x,v-\bar v)h(\bar x,\bar v)\dif\bar x\right|^2\dif \bar v\\
&=\int_{\mR}\left|\int_{\mR}\check\phi^\theta_j(\bar x,\bar v)h(x-\bar x,v-\bar v)\dif\bar x\right|^2\dif \bar v.
\end{align*}
Let $q,s\in[1,\infty]$ be defined by
$$
\tfrac1q+\tfrac1r=\tfrac12+\tfrac1p,\ \ \tfrac1s+\tfrac1r=1+\tfrac1p.
$$
Since $p\geq 2$, by Minkowskii's inequality and Young's inequality, we have
\begin{align}
&\|\cR^\theta_j \mG_N(h)\|_{\mL^p_z(\ell^2)}^p
=\int_{\mR}\|\cR^\theta_j (h\eta_\cdot)(\cdot,v)\|^p_{\mL^p_x(\ell^2)}\dif v\no\\
&\leq\int_\mR\left(\int_{\mR}\left\|\int_{\mR}\check\phi^\theta_j(\bar x,\bar v) h(\cdot-\bar x,v-\bar v)
\dif\bar x\right\|^2_{\mL^{p}_x}\dif \bar v\right)^{p/2}\dif v\no\\
&\leq\int_\mR\left(\int_{\mR}\|\check\phi^\theta_j(\cdot,\bar v)\|^2_{\mL^s_x} \|h(\cdot,v-\bar v)\|_{\mL^{r}_x}^2
\dif \bar v\right)^{p/2}\dif v\no\\
&\leq \|\check\phi^\theta_j\|_{\mL^q_v(\mL^s_x)}^{p}\|h\|^{p}_{\mL^r_v(\mL^r_x)}
\lesssim 2^{p(4-\frac3s-\frac1q)j}\|h\|^{p}_{\mL^r_z}=2^{-p\beta j}\|h\|^{p}_{\mL^r_z},\label{DJ7}
\end{align}
where the last inequality is due to the scaling property of $\check\phi^\theta_j$ and $\check\phi^\theta_1\in\cS$. Thus,
by Fubini's theorem,
\begin{align*}
\|\mG_N(h)\|_{\bB^{\beta}_{p;\theta}(\mL^p_\omega(\ell^2))}
&=\sup_{j\geq 0}2^{\beta j}\|\cR^\theta_j \mG_N(h)\|_{\mL^p_z(\mL^p_\omega(\ell^2))}\\
&=\sup_{j\geq 0}2^{\beta j}\|\cR^\theta_j \mG_N(h)\|_{\mL^p_\omega(\mL^p_z(\ell^2))}\lesssim\|h\|_{\mL^p_\omega(\mL^r_z)}.
\end{align*}
Moreover, for $\beta'<\beta$, by the dominated convergence theorem we have
\begin{align*}
\lim_{N\to\infty}\|\mG_N(h)-\mG(h)\|_{\bB^{\beta'}_{p;\theta}(\mL^p_\omega(\ell^2))}
&=\lim_{N\to\infty}\sup_{j\geq 0}2^{\beta' j}\|\cR^\theta_j (\mG_N(h)-\mG(h))\|_{\mL^p_z(\mL^p_\omega(\ell^2))}\\
&\leq\sum_{j\geq 0}2^{\beta' j}\lim_{N\to\infty}\|\cR^\theta_j (\mG_N(h)-\mG(h))\|_{\mL^p_\omega(\mL^p_z(\ell^2))},
\end{align*}
which converges to zero since for each $x,v$,
$$
\lim_{N\to\infty}\|\cR^\theta_j (\mG_N(h)-\mG(h))(x,v)\|_{\ell^2}=0.
$$
The proof is complete.
\end{proof}

\br\rm
Fix $\delta>0$ and $v_0\in\mR$. By \eqref{DC8} and BDG's inequality, it is similar to \eqref{DJ7} that for any $p\geq 2$ and $j\in\mN_0$,
\begin{align*}
\|\cR_j(\chi^\delta_{v_0}\dot B(t,\cdot))\|_{\mL^p_v(\mL^p_\omega)}&
=\left(\int_{\mR}\mE\Big|\sum_k \cR_j(\chi^\delta_{v_0}\eta_k)(v) W^k_t\Big|^p\dif v\right)^{1/p}\\
&\lesssim t^{\frac 12}\left(\int_{\mR}\|\cR_j(\chi^\delta_{v_0}\eta_\cdot)(v)\|^p_{\ell^2}\dif v\right)^{1/p}\\
&\lesssim t^{\frac 12} 2^{\frac j2}\|\chi^\delta_{v_0}\|_{\mL^p_v}=t^{\frac 12} 2^{\frac j2}\|\chi^\delta_{0}\|_{\mL^p_v},
\end{align*}
which implies that
$$
\dot B(t,\cdot)\in \cap_{p\geq 2}\wt\bB^{-1/2}_{p}(\mL^p_\omega).
$$
\er
\subsection{SKEs with Lipschitz coefficients}
In this subsection we make the following Lipschitz assumptions about $f$ and $g$.
\begin{enumerate}[{\bf (H$^s_{f,g}$)}]
\item For some $s\in(8,\infty]$, there is a $\xi\in\cap_{T>0}\mL^\infty_{T,\omega}(\mL^s_z)$ such that
for all $(t,\omega,z)\in\mR_+\times\Omega\times\mR^2$ and $u,u'\in\mR$,
\begin{align}\label{CN1}
|(f,g)(t,\omega,z,u)-(f,g)(t,\omega,z,u')|&\leq \xi(t,\omega,z)|u-u'|.
\end{align}
\end{enumerate}
Since the time variable is not important in the estimates below, we shall drop the time variables in $f$ and $g$.
\bl\label{Le74}
Under {\bf (H$^s_{f,g}$)}, for any $p\geq 2$,
there is a constant $C_0>0$ only depending on $p,s$ and $\sup_{t,\omega}\|\xi(t,\omega,\cdot)\|_{\mL^s_z}$ 
such that for all $u,u'\in\wt\mL^p_z(\mL^p_\omega)$,
\begin{align*}
\nor f(u)-f(u')\nor_{\wt\bB^{\beta}_{p;\theta}(\mL^p_\omega)}+
\nor \mG(g(u)-g(u'))\nor_{\wt\bB^{\beta}_{p;\theta}(\mL^p_\omega(\ell^2))}\lesssim_{C_0}
\nor u-u'\nor_{\wt\mL^p_z(\mL^p_\omega)},
\end{align*}
where $\beta=-\frac12-\frac 4s$, and for all $u,u'\in\mL^p_z(\mL^p_\omega)$,
\begin{align*}
\| f(u)-f(u')\|_{\bB^{\beta}_{p;\theta}(\mL^p_\omega)}+
\| \mG(g(u)-g(u'))\|_{\bB^{\beta}_{p;\theta}(\mL^p_\omega(\ell^2))}\lesssim_{C_0}
\| u-u'\|_{\mL^p_z(\mL^p_\omega)}.
\end{align*}
\el
\begin{proof}
We only prove the first one.
Let $r\in[1,p]$ be defined by $\frac1r=\frac1p+\frac1s$. Fix $\delta>0$ and $z_0\in\mR^2$.
Let $\chi^\delta_{z_0}$ be as in \eqref{Cu1}.
For $\beta=-\frac12-\frac 4s$,  we have
\begin{align*}
&\|\chi^\delta_{z_0}(f(u)-f(u'))\|_{\bB^{\beta}_{p;\theta}(\mL^p_\omega)}
\stackrel{\eqref{CT1}}{\leq}\|\chi^\delta_{z_0}(f(u)-f(u'))\|_{\mL^p_\omega(\bB^{\beta}_{p;\theta})}\\
&\quad\stackrel{\eqref{Ber0}}{\lesssim}\|\chi^\delta_{z_0}(f(u)-f(u'))\|_{\mL^p_\omega(\bB^{-1/2}_{r;\theta})}
\lesssim\|\chi^\delta_{z_0}(f(u)-f(u'))\|_{\mL^p_\omega(\mL^r_z)}\\
&\quad\stackrel{\eqref{CN1}}{\lesssim}\|\chi^\delta_{z_0}\xi|u-u'|\|_{\mL^p_\omega(\mL^r_z)}
\leq\|\chi^{\delta}_{z_0}\xi\|_{\mL^\infty_{T,\omega}(\mL^s_z)}\|\chi^{2\delta}_{z_0}(u-u')\|_{\mL^p_\omega(\mL^p_z)},
\end{align*}
where the last step is due to H\"older's inequality, and also,
\begin{align*}
&\|\chi^\delta_{z_0}\mG(g(u)-g(u'))\|_{\bB^{\beta}_{p;\theta}(\mL^p_\omega(\ell^2))}
\stackrel{\eqref{DJ0}}{\lesssim}\|\chi^\delta_{z_0}(g(u)-g(u'))\|_{\mL^p_\omega(\mL^r_z)}\\
&\quad\stackrel{\eqref{CN1}}{\lesssim}\|\chi^\delta_{z_0}\xi|u-u'|\|_{\mL^p_\omega(\mL^r_z)}
\lesssim\|\chi^{\delta}_{z_0}\xi\|_{\mL^\infty_{T,\omega}(\mL^s_z)}\|\chi^{2\delta}_{z_0}(u-u')\|_{\mL^p_\omega(\mL^p_z)}.
\end{align*}
Combining the above two estimates and by the definition of $\nor\cdot\nor_{{\bB^{\beta}_{p;\theta}(\mL^p_\omega)}}$,
we obtain the first estimate. 
\end{proof}

Now we can state and prove our main result of this section.
\bt\label{Th64}
Let $p\geq 2$, $s\in(8,\infty]$, $\beta=-\frac12-\frac4s$ and $\alpha\in(-1,1+\beta]$.
 Under {\bf (H$^s_{f,g}$)}, for any $\sF_0$-measurable $u_0\in \wt\bB^{\alpha}_{p;\theta}(\mL^p_\omega)$  and $T>0$,
there is a unique weak solution $u\in \mL^\infty_T(\wt\mL^p_z(\mL^p_\omega))$ to SKE \eqref{CN2} with regularity estimate: 
\begin{align}\label{Reg}
\nor u(t)\nor_{\widetilde{\bB}^{1+\beta}_{p;\theta}(\mL^p_\omega)}\lesssim_{C_1} t^{\frac{\alpha-1-\beta}2}
\nor u_0\nor_{\widetilde{\bB}^{\alpha}_{p;\theta}(\mL^p_\omega)}+\tilde \cI_p(f,g),\ \ t\in(0,T],
\end{align}
where $C_1=C_1(C_0,T,p,\alpha,s,\kappa)>0$ and
\begin{align}\label{YY0}
\tilde \cI_p(f,g):=\nor f(\cdot,0)\nor_{\mL^\infty_T(\wt\bB^{\beta}_{p;\theta}(\mL^p_\omega))}
+\nor \mG(g(\cdot,0))\nor_{\mL^\infty_T(\wt\bB^{\beta}_{p;\theta}(\mL^p_\omega(\ell^2)))}.
\end{align}
Moreover, if $u_0\in \bB^{\alpha}_{p;\theta}(\mL^p_\omega)$, then we also have
\begin{align}\label{Reg0}
\| u(t)\|_{{\bB}^{1+\beta}_{p;\theta}(\mL^p_\omega)}\lesssim_{C_1} t^{\frac{\alpha-1-\beta}2}
\| u_0\|_{{\bB}^{\alpha}_{p;\theta}(\mL^p_\omega)}+\cI_p(f,g),\ \ t\in(0,T],
\end{align}
where
\begin{align}\label{YY1}
\cI_p(f,g):=\| f(\cdot,0)\|_{\mL^\infty_T(\bB^{\beta}_{p;\theta}(\mL^p_\omega))}
+\| \mG(g(\cdot,0))\|_{\mL^\infty_T(\bB^{\beta}_{p;\theta}(\mL^p_\omega(\ell^2)))}.
\end{align}
In this case, for $N\in\mN$, let $u_N$ be the unique weak solution of the following SKE driven by finitely many Brownian motions:
$$
\dif u_N=\big[\Delta_vu_N+v\cdot\nabla_x u_N+f(u_N)\big]\dif t+\mG^k_N(g(u_N))\dif W^k_t,\ \ u_N(0)=u_0.
$$
If $p>\frac2{1+\beta}$, then for any $\gamma\in(0,1+\beta-\frac2p)$,
\begin{align}\label{Lim8}
\lim_{N\to\infty}\|u-u_N\|_{\mL^p_\omega(C([0,T]; \bB^\gamma_{p;\theta}))}=0.
\end{align}
\et
\begin{proof}
We divide the proof into four steps.

(i) We use Picard's iteration to show the existence of a weak solution. 
Let 
$$
u_1(t):=P_tu_0,\ \ t>0.
$$ 
For $n\geq 2$, let $u_n$ solve the following linear SKE:
\begin{align}\label{CN8}
\dif u_n=\big[\Delta_vu_n+v\cdot\nabla_x u_n+f(u_{n-1})\big]\dif t+\mG^k(g(u_{n-1}))\dif W^k_t
\end{align}
with $u_n(0)=u_0$.
Since $s\in(8,\infty]$ and $\alpha\in(-1,1+\beta]$, one can choose $\kappa\in(0,\tfrac1{2}-\tfrac4s)$ such that
$$
0<\delta:=\tfrac1{2}-\tfrac4s-\kappa<\alpha+1.
$$
Note that by  Lemma \ref{Le74},
$$
\nor f(u)\nor_{\wt\bB^{\beta}_{p;\theta}(\mL^p_\omega)}
+\nor \mG(g(u))\nor_{\wt\bB^{\beta}_{p;\theta}(\mL^p_\omega)}
\lesssim\nor u\nor_{\wt\mL^p_z(\mL^p_\omega)}+\tilde \cI_p(f,g),
$$
where $\tilde \cI_p(f,g)$ is defined by \eqref{YY0}.
Thus by Theorem \ref{Th326} with therein $\beta=-\frac32-\frac4s$ and $q=2$, there is a constant $C>0$
independent of $n$ such that for all $t\in(0,T]$,
\begin{align*}
&\nor u_n(t)\nor^2_{\widetilde{\bB}^{\delta}_{p;\theta}(\mL^p_\omega)}
\lesssim_C t^{\alpha-\delta}
\nor u_0\nor^2_{\widetilde{\bB}^{\alpha}_{p;\theta}(\mL^p_\omega)}
+\int^t_0(t-s)^{\kappa-1}\nor u_{n-1}(s)\nor^2_{\wt\mL^p_z(\mL^p_\omega)}\dif s+\tilde \cI^2_p(f,g).
\end{align*}
If we let
$$
h_N(t):=\sup_{n=1,\cdots,N}\nor u_n(t)\nor^2_{\widetilde{\bB}^{\delta}_{p;\theta}(\mL^p_\omega)},
$$
then
$$
h_N(t)\lesssim_Ct^{\alpha-\delta}
\nor u_0\nor^2_{\widetilde{\bB}^{\alpha}_{p;\theta}(\mL^p_\omega)}+\int^t_0(t-s)^{\kappa-1}h_N(s)\dif s+\tilde \cI^2_p(f,g).
$$
Since $\alpha-\delta>-1$, by Gronwall's inequality of Volterra's type (see \cite{Zh10}), we get
\begin{align}\label{DX4}
\nor u_n(t)\nor^2_{\widetilde{\bB}^{\delta}_{p;\theta}(\mL^p_\omega)}
\leq h_N(t)\lesssim_C t^{\alpha-\delta}
\nor u_0\nor^2_{\widetilde{\bB} ^{\alpha}_{p;\theta}(\mL^p_\omega)}+\tilde \cI^2_p(f,g).
\end{align}

(ii) Let $U_{n,m}(t):=u_n(t)-u_m(t)$. By Theorem \ref{Th326} with therein $\beta=-\frac32-\frac4s$ and $q=2$,  
and by Lemma \ref{Le74},
there is a constant $C>0$ such that for all $t\in(0,T]$ and $n,m\geq 2$,
\begin{align}\label{DX5}
\nor U_{n,m}(t)\nor^2_{\widetilde{\bB}^{\delta}_{p;\theta}(\mL^p_\omega)}
\lesssim_C\int^t_0(t-s)^{\kappa-1}\nor U_{n-1,m-1}(s)\nor^2_{\wt\mL^p_z(\mL^p_\omega)}\dif s.
\end{align}
Let
$$
\psi(t):=\limsup_{n,m\to\infty}\sup_{s\in(0,t]}\Big(s^{\delta-\alpha}
\nor U_{n,m}(s)\nor^2_{\widetilde{\bB}^{\delta}_{p;\theta}(\mL^p_\omega)}\Big).
$$
By \eqref{DX4} and \eqref{DX5}, it is easy to derive that
$$
\psi(T)\equiv 0.
$$
Thus, there are $u(t)\in \widetilde{\bB}^{\delta}_{p;\theta}(\mL^p_\omega)$ with
\begin{align}\label{DX44}
\nor u(t)\nor_{\widetilde{\bB}^{\delta}_{p;\theta}(\mL^p_\omega)}\lesssim_C t^{\frac{\alpha-\delta} 2}
\nor u_0\nor_{\widetilde{\bB} ^{\alpha}_{p;\theta}(\mL^p_\omega)}+\tilde \cI_p(f,g)
\end{align}
so that
$$
\limsup_{n\to\infty}\sup_{s\in(0,T]}\Big(s^{\delta-\alpha}
\nor u_{n}(s)-u(s)\nor^2_{\widetilde{\bB}^{\delta}_{p;\theta}(\mL^p_\omega)}\Big)=0.
$$
Since for any $\varphi\in C^\infty_c(\mR^2)$,
\begin{align*}
\<u_n(t),\varphi\>&=\<u_0,\varphi\>+\int^t_0\big[\<u_n,\Delta_v\varphi-v\cdot\nabla_x\varphi\>\big]\dif s\\
&\quad+\int^t_0\<f(u_{n-1}),\varphi\>\dif s+\int^t_0\<\mG^k(g(u_{n-1})),\varphi\>\dif W^k_s,
\end{align*}
by the dominated convergence theorem and taking limits $n\to\infty$ for both sides, one sees that the above $u$ is a weak solution of \eqref{CN2}.
The uniqueness of $u\in \mL^\infty_T(\wt\mL^p_z(\mL^p_\omega))$ follows from the same calculations as above.

(iii) Fix $t\in(0,T]$ and $\eps\in(0,t)$. Note that by \eqref{DX44},
$$
\sup_{s\in[\eps,t]}\nor u(s)\nor_{\widetilde{\bB}^{\delta}_{p;\theta}(\mL^p_\omega)}\lesssim_C \eps^{\frac{\alpha-\delta} 2}
\nor u_0\nor_{\widetilde{\bB} ^{\alpha}_{p;\theta}(\mL^p_\omega)}+\tilde \cI_p(f,g).
$$
Starting from the time $\eps$, by Theorem \ref{Th326} with therein $\beta=-\frac32-\frac4s$ and $q=\infty$, 
we can repeat the proof of \eqref{DX4} with $(\delta,\frac12-\frac4s)$ in place of $(\alpha,\delta)$ there, and obtain
\begin{align*}
\nor u(t)\nor_{\widetilde{\bB}^{\frac12-\frac4s}_{p;\theta}(\mL^p_\omega)}
&\lesssim_C (t-\eps)^{\frac{\delta} 2-\frac14+\frac2s}
\nor u(\eps)\nor_{\widetilde{\bB} ^{\delta}_{p;\theta}(\mL^p_\omega)}
+\sup_{s\in[\eps,t]}\nor u(s)\nor_{\wt\mL^p_z(\mL^p_\omega)}+\tilde \cI_p(f,g)\\
&\lesssim_C (t-\eps)^{\frac{\delta} 2-\frac14+\frac2s}\eps^{\frac{\alpha-\delta} 2}
\nor u_0\nor_{\widetilde{\bB} ^{\alpha}_{p;\theta}(\mL^p_\omega)}+\tilde \cI_p(f,g).
\end{align*}
In particular, taking $\eps=\frac t2$ and recalling $\beta=-\frac12-\frac4s$., we obtain \eqref{Reg}. 

(iv) As for \eqref{Reg0}, it is completely the same.
We show the limit \eqref{Lim8}. Let $U_N:=u-u_N$. By definition, we have
\begin{align*}
\dif U_N&=\big[\Delta_vU_N+v\cdot\nabla_x U_N+f(u)-f(u_N)\big]\dif t\\
&\quad+(\mG^k(g(u))-\mG^k_N(g(u_N)))\dif W^k_t.
\end{align*}
Since $p>\frac{2}{1+\beta}=2/(\frac12-\frac4s)$, one can choose $\eps$ small enough so that
$$
\beta+1-\tfrac2p-\eps>0.
$$
By Corollary \ref{Cor38} and Lemma \ref{Le74}, we have
\begin{align*}
\mE\left(\sup_{s\in[0,t]}\|U_N(s)\|^p_{\bB^\gamma_{p;\theta}}\right)
&\lesssim \int^t_0\|f(u)-f(u_N)\|^p_{\bB^{\beta-1-\eps}_{p;\theta}(\mL^p_\omega)}\dif s\\
&\quad+\int^t_0\|\mG(g(u))-\mG_N(g(u_N))\|^p_{\bB^{\beta-\eps}_{p;\theta}(\mL^p_\omega(\ell^2))}\dif s\\
&\lesssim \int^t_0\|f(u)-f(u_N)\|^p_{\bB^{\beta-1-\eps}_{p;\theta}(\mL^p_\omega)}\dif s\\
&\quad+\int^t_0\|\mG_N(g(u)-g(u_N))\|^p_{\bB^{\beta-\eps}_{p;\theta}(\mL^p_\omega(\ell^2))}\dif s\\
&\quad+\int^t_0\|(\mG-\mG_N)(g(u))\|^p_{\bB^{\beta-\eps}_{p;\theta}(\mL^p_\omega(\ell^2))}\dif s\\
&\lesssim \int^t_0\|U_N\|^p_{\mL^p_z(\mL^p_\omega)}\dif s+\int^t_0\|(\mG-\mG_N)(g(u))\|^p_{\bB^{\beta-\eps}_{p;\theta}(\mL^p_\omega(\ell^2))}\dif s,
\end{align*}
which implies \eqref{Lim8} by Gronwall's inequality and \eqref{DJ00}.
\end{proof}

\bc\label{76}
Let $s\in(8,\infty]$, $\beta=-\frac12-\frac4s$, $\gamma\in(0,\frac{1+\beta}3]$, $\alpha\in(-1,1+\beta]$ and
\begin{align}\label{Ex0}
u_0\in\cup_{\alpha>-1}\wt\bB^\alpha_{p;\theta},\ \tilde\cI_p(f,g)<\infty,\ \forall p\geq 2,
\end{align}
where $\tilde\cI_p(f,g)$ is defined by \eqref{YY0}.
 Under {\bf (H$^s_{f,g}$)}, for any $R,T>0$, $p\geq 2$ and $\eps\in(0,T)$, there is a constant $C_2>0$ depending on
$C_1,\alpha,\gamma,\eps,R$, $\tilde\cI_p(f,g)$ and $\nor u_0\nor_{\wt\bB^{\alpha}_{p;\theta}(\mL^p_\omega)}$ 
such that for the solution $u$ in Theorem \ref{Th64} and for all $\eps\leq t_1<t_2\leq T$,
\begin{align}\label{DX6}
\|(u(t_2)- u(t_1))\chi^R_0\|_{\bB^{\gamma}_{p;\theta}(\mL^p_\omega)}\lesssim_{C_2}(t_2-t_1)^{\frac{1+\beta-\gamma}3},
\end{align}
where $\chi^R_0$ is the cutoff function in \eqref{Cu1}. In particular,
for any $\delta\in(0,\frac12-\frac4s)$,
there is a finite random variable $C=C(\eps,T,R,\omega)$ such that for all
$\eps\leq t_1<t_2\leq T$ and $x_1,x_2,v_1,v_2\in B_R$,
$$
|u(t_1,\omega,x_1,v_1)-u(t_2,\omega,x_2,v_2)|\lesssim_C
(|t_1-t_2|^{\frac\delta 3}+|x_1-x_2|^{\frac\delta3}+|v_1-v_2|^{\delta}),\ a.s.
$$
\ec
\begin{proof}
For fixed $R>0$, by definition it is easy to see that
\begin{align*}
\dif (u\chi^R_0)=\big[\Delta_v(u\chi^R_0)+v\cdot\nabla_x(u\chi^R_0)+{\bf F}_R]\dif t+\mG^k(g(u)\chi^R_0)\dif W^k_t,
\end{align*}
where
$$
{\bf F}_R:=F(u)\chi^R_0-2\div_v (u\nabla_v\chi^R_0)+(\Delta_v\chi^R_0-v\cdot\nabla_v\chi^R_0)u.
$$
Thus, for $\gamma\in(0,\frac{1+\beta}3]$, by \eqref{ES9} we have
\begin{align*}
&\| u(t_2)\chi^R_0-\Gamma_{t_2-t_1}(u(t_1)\chi^R_0)\|_{\bB^{\gamma}_{p;\theta}(\mL^p_\omega)}\lesssim_C
(t_2-t_1)^{\frac{1+\beta-\gamma}2}\\
&\quad\times\Big(\|u(t_1)\|_{\bB^{1+\beta}_{p;\theta}(\mL^p_\omega)}+
\|{\bf F}_R\|_{\mL^\infty_{t_1,t_2}(\bB^{\beta-1}_{p;\theta}(\mL^p_\omega))}
+\|\mG(g(u)\chi^R_0)\|_{\mL^\infty_{t_1,t_2}(\bB^{\beta}_{p;\theta}(\mL^p_\omega(\ell^2)))}\Big).
\end{align*}
By Bernstein's inequality \eqref{Ber0}, \eqref{DJ0} and \eqref{Reg}, we clearly have
\begin{align*}
&\|{\bf F}_R\|_{\mL^\infty_{t_1,t_2}(\bB^{\beta-1}_{p;\theta}(\mL^p_\omega))}
+\|\mG(g(u)\chi^R_0)\|_{\mL^\infty_{t_1,t_2}(\bB^{\beta}_{p;\theta}(\mL^p_\omega(\ell^2)))}\\
&\quad\lesssim \|u\chi^{2R}_0\|_{\mL^\infty_{t_1,t_2}(\mL^p_z(\mL^p_\omega))}+1\lesssim t_1^{\frac{\alpha-1-\beta}{2}}
\nor u_0\nor_{\wt\bB^\alpha_{p;\theta}(\mL^p_\omega)}+1\leq C_\eps.
\end{align*}
Thus,
$$
\| u(t_2)\chi^R_0-\Gamma_{t_2-t_1}(u(t_1)\chi^R_0)\|_{\bB^{\gamma}_{p;\theta}(\mL^p_\omega)}\lesssim_C
(t_2-t_1)^{\frac{1+\beta-\gamma}2}.
$$
Moreover, by \eqref{DX11} we also have
$$
\|\Gamma_{t_2-t_1}(u(t_1)\chi^R_0)-u(t_1)\chi^R_0\|_{\bB^{\gamma}_{p;\theta}(\mL^p_\omega)}\lesssim
(t_2-t_1)^{\frac{1+\beta-\gamma}3}\|u(t_1)\chi^R_0\|_{\bB^{1+\beta}_{p;\theta}(\mL^p_\omega)}.
$$
Combining the above two estimates, we obtain \eqref{DX6}.
Finally,
the H\"older continuity of $u(t,\omega,x,v)$ in $t,x,v$ follows by
\eqref{DX6}, \eqref{Reg} and the Kolmogorov continuity theorem.
\end{proof}
\bx\label{Re64}\rm
Here we provide an example for \eqref{Ex0}.
Suppose that 
$$
u_0(x,v)=\rho(x,v)\mu(\dif v),
$$ where $\mu$ is a $\sigma$-finite measure over $\mR$ with $\sup_{v_0}\mu\{v: |v-v_0|\leq 1\}<\infty$ and $\rho(x,v)$ 
is a bounded measurable function on $\mR^2$. Then
\begin{align}\label{Mea}
u_0\in\cap_{p\geq 1}\wt\bB^{1/p-1}_{p;\theta}\subset \cap_{p\geq 2}\cup_{\alpha>-1}\wt\bB^\alpha_{p;\theta}.
\end{align}
Indeed, for $z_0=(x_0,v_0)\in\mR^2$, by definition, we have
\begin{align*}
\cR^\theta_j(\chi^\delta_{z_0}u_0)(x,v)&=\int_{\mR^2}\check\phi^\theta_j(x-\bar x,v-\bar v) 
(\chi^\delta_{z_0}\rho)(\bar x,\bar v)\dif\bar x\mu(\dif\bar v)\\
&=\int_{\mR^2}\check\phi^\theta_j(\bar x,v-\bar v) 
(\chi^\delta_{z_0}\rho)(x-\bar x,\bar v)\dif\bar x\mu(\dif\bar v).
\end{align*}
Thus by Minkowskii's inequality,
\begin{align*}
\|\cR^\theta_j(\chi^\delta_{z_0}u_0)(\cdot,v)\|_{\mL^p_x}
&\leq\int_{\mR^2}|\check\phi^\theta_j(\bar x,v-\bar v) |\dif\bar x
\|(\chi^\delta_{z_0}\rho)(\cdot,\bar v)\|_{\mL^p_x}\mu(\dif\bar v)
\end{align*}
and
\begin{align*}
\|\cR^\theta_j(\chi^\delta_{z_0}u_0)\|_{\mL^p_v(\mL^p_x)}
&\leq\left\|\int_{\mR}|\check\phi^\theta_j(\bar x,\cdot) |\dif\bar x\right\|_{\mL^p_v}\int_{\mR}
\|(\chi^\delta_{z_0}\rho)(\cdot,\bar v)\|_{\mL^p_x}\mu(\dif\bar v)\\
&\lesssim 2^{j(1-\frac1p)}\mu\Big\{\bar v: |\bar v-v_0|\leq 2\delta\Big\} \lesssim  2^{j(1-\frac1p)},
\end{align*}
which in turn gives \eqref{Mea}.
\ex

\subsection{SKEs with super-linear growth coefficients}
In this subsection we consider the following super-linear growth SKE:
$$
\dif u=\big[\Delta_vu+v\cdot\nabla_x u\big]\dif t+|u|^{\gamma+1}\dif B(t,v),\ \ u(0)=u_0,
$$
or equivalently,
\begin{align}\label{CN5}
\dif u=\big[\Delta_vu+v\cdot\nabla_x u\big]\dif t+\mG^k(|u|^{\gamma+1})\dif W^k_t,\ \ u(0)=u_0.
\end{align}

We want to show the following well-posedness result.
\bt\label{Th79}
Let $\gamma\in[0,\frac18)$ and $p>\frac{12}{1-8\gamma}$. 
For any $u_0\in\mL^1_z\cap \bB^{\frac12-4\gamma}_{p;\theta}$, there exists a unique weak solution 
$u\in C_b([0,T]\times\mR^2)$ a.s. to SKE \eqref{CN5}.
\et

Before giving a proof, we first consider the following linear SKE
\begin{align}\label{CN66}
\dif u=\big[\Delta_vu+v\cdot\nabla_x u+f\big]\dif t+\mG^k(\xi u)\dif W^k_t,\ \ u(0)=u_0,
\end{align}
where $\xi:\mR_+\times\Omega\times\mR^2\to\mR$ is a bounded predictable function.
The following a priori estimate is crucial for proving Theorem \ref{Th79}.
\bl\label{78}
Suppose that $\xi$ is a bounded predictable process. Let $p\geq 2$, $u_0\in\mL^p_z$ and
$u\in \mL^\infty_T(\mL^p_z(\mL^p_\omega))$ be a weak solution of SKE \eqref{CN66}. If $u_0^+\in\mL^1_z$, 
then for any bounded stopping time $\tau$,
\begin{align}\label{L1}
\mE\|u^+(\tau)\|_{\mL^1_z}\leq \| u^+_0\|_{\mL^1_z}+\mE\int^\tau_0\|f^+(s)\|_{\mL^1_z}\dif s.
\end{align}
\el
\begin{proof}
By \eqref{Lim8} and Fatou's lemma, it suffices to show \eqref{L1} for the 
unique solution of the following SKE driven by {\it finitely} many Brownian motions:
\begin{align}\label{CN67}
\dif u=\big[\Delta_vu+v\cdot\nabla_x u+f\big]\dif t+\mG^k_N(\xi u)\dif W^k_t,\ \ u(0)=u_0,
\end{align}
where $N<\infty$.
Let $\rho_\eps(x,v):=\eps^{-4}\rho(\eps^{-1}x,\eps^{-3}v)$, where $\rho$ is smooth density function with support in the unit ball.
Define $u_\eps:=u*\rho_\eps$ and $f_\eps:=f*\rho_\eps$. Then
\begin{align*}
\dif u_\eps(t)&=(\Delta_v u_\eps+v\cdot\nabla_x u_\eps+[\rho_\eps*, v\cdot\nabla_x]u+f_\eps)\dif t
+\mG^k_N(\xi u)*\rho_\eps\dif W^k_t,
\end{align*}
where 
$$
[\rho_\eps*, v\cdot\nabla_x]u:=\rho_\eps*(v\cdot\nabla_x u)-v\cdot\nabla_x (u*\rho_\eps).
$$
Let $\psi:\mR\to[0,\infty)$ be a smooth convex smooth function with bounded derivative of all orders greater than $1$.
By It\^o's formula, we have
\begin{align*}
\psi(u_\eps(t))&=\psi(u_\eps(0))+\int^t_0\psi'(u_\eps)(\Delta_v u_\eps+v\cdot\nabla_x u_\eps
+[\rho_\eps*, v\cdot\nabla_x]u+f_\eps)\dif s\\
&+\frac12\int^t_0\psi''(u_\eps)\|\mG_N(\xi u)*\rho_\eps\|^2_{\ell^2}\dif s+\int^t_0\psi'(u_\eps)(\mG^k_N(\xi u)*\rho_\eps)\dif W^k_s.
\end{align*}
Let $\chi$ be a nonnegative smooth function with compact support. Multiplying both sides by $\chi$ and then
integrating on $\mR^{2}$ and noting that
\begin{align*}
\int\psi'(u_\eps)(\Delta_v u_\eps+v\cdot\nabla_x u_\eps)\chi
&=-\int\psi''(u_\eps)|\nabla_v u_\eps|^2\chi+\int (\Delta_v\chi- v\cdot\nabla_x\chi)\psi(u_\eps),
\end{align*}
we obtain that for any bounded stopping time $\tau$,
\begin{align*}
\mE\int\psi(u_\eps(\tau))\chi&\leq\int\psi(u_\eps(0))\chi+\mE\int^\tau_0\!\!\int (\Delta_v\chi- v\cdot\nabla_x\chi)\psi(u_\eps)\\
&\quad+\mE\int^\tau_0\!\!\int\psi'(u_\eps)([\rho_\eps*, v\cdot\nabla_x]u+f_\eps)\chi\\
&\quad+\frac12\mE\int^\tau_0\!\!\int\psi''(u_\eps)\|\mG_N(\xi u)*\rho_\eps\|^2_{\ell^2}\chi.
\end{align*}
Here and below we drop the integral variables $\dif x\dif v\dif s$ for simplicity.
Letting $\eps\to 0$ and by Fatou's lemma and the dominated convergence theorem, we get
\begin{align}\label{LK9}
\begin{split}
\mE\int\psi(u(\tau))\chi&=\mE\int\lim_{\eps\to 0}\psi(u_\eps(\tau))\chi\leq \liminf_{\eps\to 0}\mE\int\psi(u_\eps(\tau))\chi\\
&\leq\int\psi(u(0))\chi+\mE\int^\tau_0\!\!\int (\Delta_v\chi- v\cdot\nabla_x\chi)\psi(u)\\
&+\mE\int^\tau_0\!\!\int\psi'(u)f\chi+\frac12\mE\int^\tau_0\!\!\int\psi''(u)\|\mG_N(\xi u)\|^2_{\ell^2}\chi.
\end{split}
\end{align}
Now we take
$$
\psi(r)=\psi_\delta(r):=(r+\sqrt{r^2+\delta})/2,\ \ \delta>0.
$$
Clearly,
$$
\lim_{\delta\downarrow 0}\psi_\delta(r)=r^+,\ \ \psi''_\delta(r)\geq 0, \ \ r^2\psi''_\delta(r)\leq \sqrt{\delta},\ \ |\psi'_\delta(r)|\leq 1.
$$
Thus by \eqref{LK9},
\begin{align*}
\mE\int\psi_\delta(u(\tau))\chi&\leq\int\psi_\delta(u(0))\chi+\mE\int^\tau_0\!\!\int (\Delta_v\chi- v\cdot\nabla_x\chi)\psi_\delta(u)\\
&\quad+\mE\int^\tau_0\int f^+\chi+\frac{\sqrt\delta}2\mE\int^\tau_0\!\!\int\|\mG_N(\xi)\|^2_{\ell^2}\chi.
\end{align*}
Letting $\delta\downarrow 0$, we get for any bounded stopping time $\tau$,
\begin{align}\label{DX2}
\mE\int u^+(\tau)\chi\leq\int u^+(0)\chi+\mE\int^\tau_0\!\!\int u^+(\Delta_v- v\cdot\nabla_x)\chi
+\mE\int^\tau_0\int f^+\chi.
\end{align}
Now, let $\chi:\mR^2\to[0,1]$ be a nonnegative smooth function with 
$$
\chi(x,v)=
\left\{
\begin{aligned}
1,&\ \ |x|^{1/3}+|v|\leq 1,\\
0, &\ \ |x|^{1/3}+|v|\geq 2.
\end{aligned}
\right.
$$ 
For $R\geq 1$, let 
$$
\chi_R(x,v):=\chi(R^{-3}x,R^{-1}v).
$$
By the chain rule, it is easy to see that
\begin{align}\label{DX3}
\begin{split}
(\Delta_v\chi_R- v\cdot\nabla_x\chi_R)(x,v)&=R^{-2}(\Delta_v\chi- v\cdot\nabla_x\chi)(R^{-3}x,R^{-1}v)\\
&\leq \|\Delta_v\chi- v\cdot\nabla_x\chi\|_\infty R^{-2}\chi_{2R}(x,v).
\end{split}
\end{align}
Thus, for any $n\in\mN_0$, by \eqref{DX2} with $\chi_{2^nR}$ in place of $\chi$, we get for any $t\in[0,T]$,
\begin{align*}
\mE\int u^+(t)\chi_{2^nR} \leq \int u^+_0+C_0(2^nR)^{-2}\mE\int^t_0\int u^+\chi_{2^{n+1}R}+\mE\int^T_0\int f^+,
\end{align*}
where $C_0:=\|\Delta_v\chi- v\cdot\nabla_x\chi\|_\infty$.
Let $H_n(t):=\mE\int u^+(t)\chi_{2^nR}$. The above inequality means that
$$
H_n(t)\leq\int u^+_0+\mE\int^T_0\int f^+\dif s+C_0(2^nR)^{-2}\int^t_0 H_{n+1}(s)\dif s.
$$
By iteration, we get for any $N\in\mN$,
\begin{align*}
H_0(t)&\leq\left(\int u^+_0+\mE\int^T_0\int f^+\dif s\right)\left(1+\sum_{n=1}^{N-1}\frac{(C_0R^{-2} t)^n}{2^{n(n-1)} n!}\right)\\
&\quad+\frac{(C_0R^{-2})^N}{2^{N(N-1)}}\int^t_0\!\!\int^{s_1}_0\!\!\!\cdots\!\!\int^{s_{N-1}}_0H_N(s_N)\dif s_N\cdots\dif s_1.
\end{align*}
Since by H\"older's inequality and the assumption, for $q=\frac{p}{p-1}\in(1,2]$,
$$
H_N(t)\leq\mE\|u(t)\|_p\|\chi_{2^NR}\|_q\leq C_1(2^N R)^{\frac4q},
$$
where $C_1$ does not depend on $N,R,t$, we have
$$
\frac{(C_0R^{-2})^N}{2^{N(N-1)}}\int^t_0\!\!\int^{s_1}_0\!\!\!\cdots\!\!\int^{s_{N-1}}_0H_N(s_N)\dif s_N\cdots\dif s_1
\leq \frac{C_1(2^N R)^{\frac4q} (C_0R^{-2}t)^N}{ 2^{N(N-1)}N!}\to 0
$$
as $N\to\infty$. Hence,
$$
\mE\int u^+(t)\chi_R=H_0(t)\leq\left(\int u^+_0+\mE\int^T_0\int f^+\dif s\right)\e^{C_0R^{-2}t}.
$$
Letting $R\to\infty$, we obtain that for any $t\in[0,T]$,
$$
\mE\int u^+(t)\leq \int u^+_0+\mE\int^T_0\int f^+.
$$
Finally, by this estimate and \eqref{DX2}, \eqref{DX3} again,  
$$
\mE\int u^+(\tau)\chi_R\leq\int u^+(0)\chi_R+\frac{\|\Delta_v\chi- v\cdot\nabla_x\chi\|_\infty}{R^{2}}
\mE\int^\tau_0\!\!\int u^++\mE\int^\tau_0\int f^+\chi_R,
$$
which implies the desired result by taking limits $R\to \infty$.
\end{proof}

\br\rm
Suppose $f^-\equiv 0$. If $u(0)\geq 0$, then $u(t)\geq 0$ a.s. by \eqref{L1}. This in turn implies the nonnegativity of the solution to SKE \eqref{CN5} provided $u_0\geq 0$.
\er
\br\rm\label{Re72}
Estimate \eqref{L1} is also used to derive the local uniqueness of SKE \eqref{CN2} under {\bf (H$^\infty_{f,g}$)}.
More precisely, let $u_1$ and $u_2$ be two solutions of SKE \eqref{CN2} before stopping time $\tau\leq T$
with the same initial values, i.e., for any $\varphi\in C^\infty_c(\mR^2)$,
\begin{align*}
\<u_i(t),\varphi\>&=\<u_0,\varphi\>+\int^t_0\big[\<u_i,\Delta_v\varphi-v\cdot\nabla_x\varphi\>+\<f(u_i),\varphi\>\big]\dif s\\
&\quad+\int^t_0\<\mG^k(g(u_i)),\varphi\>\dif W^k_s,\ \ t\leq\tau,\ \ a.s.,\ \ i=1,2.
\end{align*}
Then it holds that
 $$
 u_1(t)=u_2(t),\ \ \forall t\leq\tau,\ \ a.s.
 $$
Indeed, let $U:=u_1-u_2$. Then $U$ solves the following SKE: for $t\leq\tau$,
$$
\dif U=[\Delta_v U+v\cdot\nabla_x U+\xi_0 U]\dif t+\mG^k(\xi_1 U)\dif W^k_t,\ U(0)=0,
$$
where $\xi_0:=(f(u_1)-f(u_2))/(u_1-u_2)$ and $\xi_1:=(g(u_1)-g(u_2))/(u_1-u_2)$. Here we use the convention $\frac00:=0.$
Under {\bf (H$^\infty_{f,g}$)},  $\xi_0$ and $\xi_1$ are bounded predictable. Using \eqref{L1}, we have
$$
\mE\|U(t\wedge\tau)\|_{\mL^1_z}\leq \mE\int^{t\wedge\tau}_0\|\xi_0 U(s)\|_{\mL^1_z}\dif s\leq \|\xi_0\|_\infty\mE\int^t_0\|U(s\wedge\tau)\|_{\mL^1_z}\dif s.
$$
By Gronwall's inequality, $U(t\wedge\tau)=0$ for any $t\in[0,T]$. We would like to point out that the local uniqueness seems not be an obvious consequence of 
the global uniqueness since for $u^\tau(t):=u(t\wedge\tau)$, we only have
$$
\p_t u^\tau=\1_{t\leq\tau}\big[\Delta_vu^\tau+v\cdot\nabla_x u^\tau+f(u^\tau)\big]\dif t+\mG^k(g(u^\tau))\dif W^k_{t\wedge \tau},
$$
which has a random degenerate coefficient in front of the Laplacian operator.
\er
\br\rm 
By \eqref{L1} and \cite[Theorem III 6.8]{Kry95} or \cite[Lemma 2.7]{Zh-Zh0}, when  $f=0$, we have
\begin{align}\label{LK0}
\mE\left(\sup_{t\in[0,T]}\|u(t)\|_{\mL^1_z}^{1/2}\right)\leq 2\|u_0\|^{1/2}_{\mL^1_z}.
\end{align}
\er

Now we are in a position to give

\begin{proof}[Proof of Theorem \ref{Th79}]
For $m\in\mN$, we define
$$
g_m(u):=(|u|\wedge m)^{\gamma+1}.
$$
Clearly, $g_m$ is global Lipschitz. By Theorem \ref{Th64}, the following SKE has a unique weak solution $u_m\in\mL^\infty_T(\bB^{\frac12-4\gamma}_{p;\theta}(\mL^p_\omega))$
\begin{align}\label{CN6}
\dif u=\big[\Delta_vu+v\cdot\nabla_x u\big]\dif t+\mG^k(g_m(u))\dif W^k_t,\ \ u(0)=u_0.
\end{align}
From the definition, one sees that $u_m$ is also a weak solution of the following linearized equation:
$$
\dif w=\big[\Delta_vw+v\cdot\nabla_x w\big]\dif t+\mG^k(\xi_m w)\dif W^k_t,\ \ w(0)=u_0,
$$
where
$$
\xi_m(t,x,v):=g_m(u_m(t,x,v))/u_m(t,x,v),\ \ \tfrac 00:=0.
$$
Since $|\xi_m|\leq m^\gamma$, by \eqref{LK0}, we have
$$
\mE\left(\sup_{t\in[0,T]}\|u_m(t)\|^{1/2}_{\mL^1_z}\right)\leq 2\|u_0\|^{1/2}_{\mL^1_z}.
$$
In the following we shall use this a priori estimate to derive that for any $T>0$,
\begin{align}\label{LK5}
\lim_{R\to\infty}\sup_m\mP\left(\omega:\sup_{t\in[0,T]}\|u_m(t,\omega,\cdot)\|_{C_b(\mR^2)}\geq R\right)=0.
\end{align}
For $S>0$, define a stopping time
$$
\tau^S_m:=\inf\Big\{t\geq 0: \|u_m(t,\cdot)\|_{\mL^1_z}\geq S\Big\}.
$$
Noting that
$$
|\xi_m(t,\omega,z)|\leq |u_m(t,\omega,z)|^\gamma,
$$
we have for $s=1/\gamma$,
$$
\|\xi_m(t\wedge \tau^S_m,\omega,\cdot)\|_{\mL^s_z}\leq S^\gamma.
$$
On the other hand, by Duhamel's formula, we also have
$$
u_m(t)=P_tu_0+\int^t_0P_{t-s}\mG^k(\xi_mu_m)\dif W^k_s.
$$
Let $\frac1r=\frac1p+\frac1s=\frac1p+\gamma$ and $\beta=4(\frac1p-\frac1r)-\frac12=-\frac12-4\gamma$.
For any $\delta\in(0,\beta+1-\frac2p)$, by \eqref{HH7}, \eqref{323}, \eqref{DJ0} and H\"older's inequality, we have
\begin{align*}
&\mE\left(\sup_{s\in[0,t\wedge\tau^S_m]}\|u_m(s)\|^p_{\bB^{\delta}_{p;\theta}}\right)
\lesssim\|u_0\|^p_{\bB^{\beta+1-\frac2p}_{p;\theta}}
\\&\quad+\int^t_0\|\mG((\xi_m u_m)(\cdot\wedge\tau^S_m))(s)\|^p_{\bB^{\beta}_{p;\theta}(\mL^p_\omega(\ell^2))}\dif s
\\&\quad\lesssim 1+\int^t_0\|(\xi_m u_m)(s\wedge\tau^S_m)\|^p_{\mL^p_\omega(\mL^r_z)}\dif s
\\&\quad\leq 1+S^{\gamma p}\int^t_0\|u_m(s\wedge\tau^S_m)\|^p_{\mL^p_\omega(\mL^p_z)}\dif s,
\end{align*}
which implies by Gronwall's inequality that
$$
\sup_m\mE\left(\sup_{s\in[0,t\wedge\tau^S_m]}\|u_m(s)\|^p_{\bB^{\delta}_{p;\theta}}\right)\leq C_1=C_1(S).
$$
Since $p>\frac{12}{1-8\gamma}$, one can choose $\delta\in(\frac4p,\beta+1-\frac2p)=(\frac4p,\frac12-4\gamma-\frac2p)$ so that by the above estimate and \eqref{Ber1} with $q=\infty$, 
$$
\sup_m\mE\left(\sup_{t\in[0,T\wedge\tau^S_m]}\|u_m(t)\|^p_{C_b(\mR^2)}\right)\leq C_2=C_2(S).
$$
By Chebychev's inequality, we have
\begin{align*}
&\mP\left(\omega:\sup_{t\in[0,T]}\|u_m(t,\omega,\cdot)\|_{C_b(\mR^2)}\geq R\right)\leq 
\mP(\tau^S_m\leq T)\\
&\qquad+\mP\left(\omega:\sup_{t\in[0,T\wedge\tau^S_m]}\|u_m(t,\omega,\cdot)\|_{C_b(\mR^2)}\geq R\right)
\\&\quad\leq\frac{1}{S^{1/2}}\mE\left(\sup_{t\in[0,T]}\|u_m(t,\cdot)\|^{1/2}_{\mL^1_z}\right)
+\frac1R\mE\left(\sup_{t\in[0,T\wedge\tau^S_m]}\|u_m(t)\|^p_{C_b(\mR^2)}\right)\\
&\quad\leq\frac{2\|u_0\|^{1/2}_{\mL^1_z}}{S^{1/2}}+\frac{C_2}{R},
\end{align*}
which yields \eqref{LK5} by firstly letting $R\to\infty$ and then $S\to\infty$.

Finally, by \eqref{LK5}, we can construct a unique solution of \eqref{CN5} as in \cite{kry99}. Indeed, for given $m,n\in\mN$,  we define
$$
\tau_m^n:=\inf\Big\{t>0: \|u_m(t)\|_{C_b(\mR^2)}\geq n\Big\}.
$$
Since for $m>n$, $u_m$ and $u_n$ satisfy the same equation \eqref{CN6} with coefficient $g_n$ before $\tau^n_m$.
By the local uniqueness (see Remark \ref{Re72}), we have $u_m|_{[0,\tau^n_m)}=u_n|_{[0,\tau^n_m)}$. Hence,
$$
\tau^m_m\geq \tau^n_m\geq \tau^n_n,\ \ a.s.,
$$
and we can define without ambiguity
$$
u(t,x,v):=u_m(t,x,v),\ \ t\leq\tau^m_m.
$$
Clearly, $u$ is a unique weak solution of \eqref{CN5} before $\tau^m_m$. By \eqref{LK5}, one sees that
\begin{align*}
\mP(\tau^m_m\leq T)&\leq\mP\left(\sup_{t\in[0,T]}\|u_m(t)\|_{C_b(\mR^2)}\geq m\right)\\
&\leq\sup_n\mP\left(\sup_{t\in[0,T]}\|u_n(t)\|_{C_b(\mR^2)}\geq m\right)\to0
\end{align*}
as $m\to\infty$. The proof is complete.
\end{proof}

{\bf Acknowledgement:} The authors would like to thank Zimo Hao for his quite useful discussions.


\end{document}